\documentclass[12pt,english]{amsart}
\usepackage[T1]{fontenc}
\usepackage[latin1]{inputenc}
\usepackage{geometry}
\usepackage{graphicx}
\usepackage{amssymb}

\makeatletter

\newcommand{\noun}[1]{\textsc{#1}}

 \theoremstyle{plain}    
 \newtheorem{thm}{Theorem}[section]
 \numberwithin{equation}{section} 
 \numberwithin{figure}{section} 
 \theoremstyle{plain}
 \theoremstyle{remark}
 \newtheorem{rem}[thm]{Remark}
 \theoremstyle{plain}    
 \newtheorem{prop}[thm]{Proposition} 
 \theoremstyle{remark}
 \newtheorem*{rem*}{Remark}
 \theoremstyle{plain}    
 \newtheorem{lem}[thm]{Lemma} 
 \theoremstyle{remark}    
 \newtheorem*{conclusion*}{Conclusion} 
 \theoremstyle{definition}
  \newtheorem*{example*}{Example}
 \theoremstyle{remark}    
 \newtheorem*{claim*}{Claim}
 \theoremstyle{plain}    
 \newtheorem{cor}[thm]{Corollary} 

\usepackage{amsmath}
\usepackage{amsfonts}
\usepackage{amssymb}
\input{diagrams}
\input xy
\xyoption{all}
\def\mapiotaj{\buildrel \iota_j \over \to}
\def\mapiota{\buildrel \iota \over \hookrightarrow}

\def\mapbi{\buildrel \mathcal{I} \over \hookrightarrow}
\def\mapbito{\buildrel \mathcal{I} \over \to}

\def\mapcong{\buildrel \cong \over \to}
\def\mapR{\buildrel \mathcal{R} \over \to}
\def\mapRprime{\buildrel ' \mathcal{R} \over \to}
\def\mapb{\buildrel b \over \to}
\def\mapbproj{\buildrel b \over \twoheadrightarrow}
\def\mapbstar{\buildrel b_* \over \twoheadrightarrow}
\def\mapbde{\buildrel b_*^{\de} \over \twoheadrightarrow}
\def\mapaa{\buildrel \iota_1 \over \lInto}
\def\mapbb{\buildrel \iota_2 \over \rInto}
\def\mapcc{\buildrel \mathcal{N}_* \over \to}
\def\mapdd{\buildrel \alpha \over \to}
\def\mapee{\buildrel \Phi \over \to}
\def\mapff{\buildrel \text{Res} \over \to}
\def\mapgg{\buildrel \tilde{\alpha} \over \to}
\def\maphh{\buildrel \text{mod} \, \Lambda \over \equiv}
\def\mapii{\buildrel \text{HC} \over \implies}
\def\mapjj{\buildrel p^*_i \over \to}
\def\mapkk{\buildrel p \over \to}
\def\C{\mathbb{C}}

\def\Z{\mathbb{Z}}
\def\R{\mathbb{R}}
\def\P{\mathbb{P}}
\def\Q{\mathbb{Q}}
\def\c{\mathcal{C}}
\def\z{\mathcal{Z}}
\def\m{\mathcal{M}}
\def\de{\mathcal{D}}
\def\s{\square}
\def\d{\partial}

\def\a{\alpha}
\def\B{\beta}

\def\f{\textbf{f}}

\def\w{\omega}

\def\c{\mathcal{C}}

\def\W{\mathcal{W}}
\def\im{\text{im}}
\def\F{\mathcal{F}}
\def\t{\otimes}

\def\l{\mathcal{L}}
\def\v{\wedge}

\def\dlog{\text{dlog}}

\def\dd{\mathbb{D}}
\def\th{\theta}

\def\di{\text{d}}

\def\i{\iota}
\def\m{\setminus}

\def\DI{\, ' \mathcal{D}}
\def\O{\Omega}
\def\j{\jmath}

\def\o{\mathcal{O}}
\def\ss{\mathcal{S}}
\def\H{\mathcal{H}}

\def\HH{\mathbb{H}}

\def\KK{\mathcal{K}}

\def\EE{\mathcal{E}}

\def\nequiv{/\mspace{-17mu}\equiv}
\def\nsimeq{/\mspace{-18mu}\simeq}
\def\uH{\underline{H}_{\H}}

\def\A{\textbf{a}}
\def\uh{\underline{H}}
\def\sh{\mathsf{H}}
\def\ns{\eta_S}
\def\BI{\mathbb{I}}
\def\WV{\widetilde{V}}
\def\WL{\widetilde{\Lambda}}
\def\ratnequiv{\buildrel \text{rat} \over \nequiv}
\def\homeq{\buildrel \text{hom} \over \equiv}
\def\ratequiv{\buildrel \text{rat} \over \equiv}
\def\mappsi{\buildrel \Psi \over \hookrightarrow}
\def\mappiN{\buildrel \pi_N \over \longrightarrow}
\def\mapq{\buildrel q \over \to}
\def\mapp{\buildrel p \over \to}
\def\mapd{\buildrel \text{d} \over \to}
\def\algequiv{\buildrel \text{alg} \over \equiv}
\def\mapblah{\buildrel \oplus ( \mathcal{R} \circ e^{-1} \circ \Xi^{-1}_i ) \over \longrightarrow}
\def\mapD{\buildrel \text{D} \over \to}
\def\mapphipp{\buildrel ''\phi_{i-1} \over \twoheadrightarrow}
\def\equivmod{\buildrel \text{mod} \, \Lambda \over \equiv}
\def\mapnabla{\buildrel \bar{\nabla} \over \to}
\def\mappia{\buildrel \pi_0 \over \twoheadrightarrow}
\def\mappib{\buildrel \pi \over \twoheadleftarrow}
\def\mappi{\buildrel \pi \over \to}
\def\mapidi{\buildrel \text{id} \times \iota_q \over \hookrightarrow}
\def\maprho{\buildrel \rho \over \to}
\def\mapj{\buildrel \jmath \over \to}

\usepackage{babel}
\makeatother
\begin{document}

\title{Higher Abel-Jacobi maps for 0-cycles}

\author{Matt Kerr}

\begin{abstract}
Starting from the candidate Bloch-Beilinson filtration on $CH_{0}(X_{\C})$
constructed in \cite{L2}, we develop and describe geometrically a
series of Hodge-theoretic invariants $\Psi_{i}$ defined on the graded
pieces. Explicit formulas (in terms of currents and membrane integrals)
are given for certain quotients of the $\Psi_{i}$, with applications
to $0$-cycles on products of curves.
\end{abstract}



\maketitle

\section{\textbf{Introduction}}

A classical question in transcendental algebraic geometry, which has
seen a resurgence of interest (e.g. \cite{Gr1}, \cite{GG3}, \cite{L2},
\cite{RS}) in recent years, is that of $rational$ $equivalence$
of algebraic cycles on a smooth projective variety $X/\C$. For instance,
suppose we take two formal sums $\z_{1}$ and $\z_{2}$ of points
with integer multiplicities (or $0$-cycles) on $X$. When can we
rationally parametrize a {}``path'' between the two collections?
That is, writing $\s:=\P_{\C}^{1}\m\{1\}$ for affine space, when
does there exist an algebraic $1$-cycle $W\in Z_{1}(X\times\s)$
such that \[
\z_{1}-\z_{2}=\pi_{*}^{X}(W\,\cdot\,\, X\times\{0\})-\pi_{*}^{X}(W\,\cdot\,\, X\times\{\infty\})\,\,\,\,?\]
 Or equivalently, when do there exist curves $\W_{j}$ on $X$ (with
desingularizations $\widetilde{\W}_{j}\mapiotaj X$) and functions
$f_{j}\in\C(\widetilde{\W}_{j})^{*}$ such that $\z_{1}-\z_{2}=\sum\i_{j}{}_{_{*}}(f_{j})$
?

For $\dim_{\C}X=1$, this is the familiar question of whether $\z_{1}-\z_{2}$
is the divisor of a function, and it is {}``solved'' by Abel's theorem:
the two invariants $\Psi_{0}:=$degree map and $\Psi_{1}$:=Abel-Jacobi
map completely detect $\ratnequiv0$ (see for example \cite{G}).

Further invariants are required when $\dim X=n>1$. The goal at which
this paper is aimed, is the explicit description of a series of Hodge-theoretically
determined, higher $AJ$-type maps $\Psi_{i}$ which {}``completely
capture'' rational equivalence classes of $0$-cycles. (Each $\Psi_{i}$
should be defined on $\ker(\Psi_{i-1})$.) To make this goal more
attainable we simplify things in a couple of ways. All cycle, Chow,
$K$-groups (and related objects) are taken $\t\Q$ (i.e. modulo torsion);
so $CH^{n}(X)$ means $CH^{n}(X)\t\Q$, while $Hg^{p}(X):=F^{p}H^{2p}(X,\C)\cap H^{2p}(X,\Q)$
and $J^{p}(X):=\frac{H^{2p-1}(X,\C)}{F^{p}\{\text{num}\}+H^{2p-1}(X,\Q)}$.
We also assume throughout that $X$ is defined over $\bar{\Q}$, which
is to say: it is cut out by equations with coefficients in a number
field.

Now restrict to those $0$-cycles defined over some fixed field $k$
finitely generated $/\bar{\Q}$. (By taking a limit over such $k\subset\C$,
one eventually gets all $0$-cycles on $X_{\C}$.) $\Psi_{0}$ and
$\Psi_{1}$ are conveniently summarized (e.g., see \cite{Gr2}) by
the Deligne cycle-class map\[
c_{\de}^{X}:\, CH^{n}(X_{k})\to H_{\de}^{2n}(X_{\C},\Q(n)).\]
 According to the Bloch-Beilinson conjecture (BBC), this is injective
if $k=\bar{\Q}$; otherwise its kernel may be huge, and so one is
led naturally into the arithmetic world in dealing with the extension
$k/\bar{\Q}$. One can in fact exchange this extension for additional
geometry by performing a $\bar{\Q}$-spread as described in $\S2$,
to obtain a cycle $\zeta\in Z^{n}(X\times\eta_{S})$ defined over
$\bar{\Q}$.

A version of the Bloch-Beilinson conjecture for quasi-projective varieties
$X\times U$ (of which $X\times\eta_{S}$ is a limit) now tells us
that a refinement\[
c_{\H}:\, CH^{n}(X\times\eta_{S}{}_{/\bar{\Q}})\to H_{\H}^{2n}(X\times\eta_{S},\Q(n))\]
 of $c_{\de}^{X\times\eta_{S}}$ is injective. ($H_{\H}^{2n}$ is
absolute Hodge cohomology as defined in \cite{L2}.) Placing a Leray
filtration on {[}a subgroup of{]} $H_{\H}^{2n}$ as in \cite{L2},
the corresponding graded pieces of $c_{\H}(\zeta)$ give the desired
invariants\[
\Psi_{i}:\,\l^{i}CH^{n}(X_{k}):=\ker(\Psi_{i-1})\to Gr_{\l}^{i}\underline{H}_{\H}^{2n}(X\times\eta_{S},\Q(n)).\]
 Intutitively, the idea is to use the product structure of $X\times\eta_{S}$
to chop up the cycle and $AJ$-classes of the spread $\zeta$ of $\z$.
This story is reviewed abstractly in Part $1$ ($\S\S2-6$). The difficulty
is in finding explicit formulas for the resulting maps; this occupies
us for much of the remainder of the paper, which is divided into Part
$2$ ($\S\S7-9$), Part $3$ ($\S\S10-12$), and Part $4$ (Applications,
$\S\S13-17$).

Two remarks are in order. First, we wish to emphasize the tenuous
nature of BBC. There are very simple cases where it is not known.
A particularly striking example was pointed out to the author by P.
Griffiths: say we present a $K3$ surface $X$ as a double cover of
$\P^{2}$ branched over a smooth sextic curve with affine equation
$f(x,y)=0$. Then for any $x,y\in\bar{\Q}$ not solving $f$,\[
\z_{x,y}:=(x,y,+\sqrt{f(x,y)})-(x,y,-\sqrt{f(x,y)})\in\ker(c_{\de}^{X})\subseteq CH^{2}(X_{\bar{\Q}})\]
 since $h^{1,0}(X)=0$; but we are not aware of any proof that (some
multiple of) $\z_{x,y}\ratequiv0$. If BBC fails (for quasi-projectives)
then the series $\Psi_{i}$ is still (well-)defined but their kernels
do not exhaust $CH^{n}(X_{k})$. So without the aid of conjectures,
at present the best we can hope for is to detect $\z\ratnequiv0$
when $c_{\de}^{X}(\z)=0$, by showing e.g. $\Psi_{2}(\z)$ or $\Psi_{3}(\z)\neq0$.
Some concrete examples (and techniques) are given in Part $4$.

The second remark, is that in case one replaces $X$ by a relative
variety $(X,Y)$ (as we do in Part $2$) one must also require that
$f_{j}\equiv1$ on $\i_{j}^{-1}(\W_{j}\cap Y)$ in the above definition
of rational equivalence of $0$-cycles.

The present paper is a fusion between the algebraic computation of
the $\Psi_{i}$ target-spaces in \cite{L2}, the geometric point of
view in \cite{GG3}, and the work on regulator maps in \cite{Ke1}.
Its central results are the presentation (in Part $3$) of a quotient
of $\Psi_{i}(\z)$, as a (geometrically defined) vector-valued differential
character on the {}``base'' $\eta_{S}$ of the spread; and (for
$X$ a product of curves) its representation via explicit $(i-1)$-currents
on $S$ (in $\S14$). These {}``reduced invariants'' $\chi_{i}(\z)$
(and the currents) are motivated by a study of analogous invariants
in the degenerate situation $X=(\s^{n},\d\s^{n})$ in Part $2$, which
turn out to be computed by the Milnor-regulator currents of \cite{Ke1}
(on $\eta_{S}$). (Here $(\s^{n},\d\s^{n})$ denotes affine $n$-space
relative its faces; we also study the situation for $X=(\P^{1},\{0,\infty\})^{n}$,
which we think of as a product of degenerate elliptic curves.) Our
general philosophy is that the degenerate invariants, or $\chi_{i}(\z)$,
can be {}``pushed down'' to functionals on topological $(i-1)$-cycles
on the base of the spread.

The geometric construction of the $\chi_{i}$, followed by the proof
that they are well-defined on $\l^{i}CH_{0}$ and are factored by
the $\Psi_{i}$, is contained in $\S11$. Involved in the proof (of
Proposition $11.2$) is an interesting technical observation: K\"unneth
components of Abel-Jacobi classes of (algebraic) cycles $\zeta\homeq0$
on quasiprojective varieties $X\times U$, whose closure $\bar{\zeta}$
on $X\times S$ may not be homologous to zero, can still be partially
computed by straightforward membrane integrals. (In fact, what we
prove is stronger than this.) See Remark $4.1$ and remarks preceding
Prop. $11.2$ for further discussion.

We also offer some description of the full invariants $\Psi_{i}$
and their targets, in $\S\S$ $4$, $5$, $10$, and $12$. That in
$\S12$ leads to a very simple proof (in $\S15$, using $\Psi_{2}$)
of a result of \cite{RS} for certain $0$-cycles on products of curves,
which we then apply to special Kummer ($K3$) surfaces in $\S17$.
Finally, we present evidence that the reduced invariant $\chi_{2}$
may do just as well as $\Psi_{2}$ at detecting such $0$-cycles.
We prove a statement along these lines for products of elliptic curves
in $\S16$ using a result (and in one special case a conjecture) on
algebraic independence of values of elliptic logarithms. This appears
to be the first result of its kind, bringing transcendence theory
to bear on higher Abel-Jacobi maps, and should be interesting to pursue
further.

We should emphasize that this paper is foundational in nature; most
of it was written before \cite{Ke4} and \cite{Ke5}, which rely heavily
on the framework and technical results of the present $\S\S2-5,\,10-14$.
Several original applications to the detection of explicit {}``new''
cycles (not previously known to be nontrivial) can be found in those
papers, e.g. $\S\S6.1,\,6.3,\,8.2$ of \cite{Ke4}. The arguments
in the present $\S\S10-11$ will also be applied to the study of higher
Chow cycles \cite{KL}.

For the orientation of the reader we point out that a review of Lewis's
description of the $\Psi_{i}$-targets $Gr_{\l}^{i}\underline{H}_{\H}^{2n}(X\times\eta_{S},\Q(n))$
may be found in $\S4$. In $\S10$, we give also an explicit {}``Deligne-homology''-type
presentation of the $Gr_{\l}^{i}\underline{H}_{\H}^{2n}$; this just
means a description obtained by resolving a relevant Deligne-Beilinson
cone complex by acyclic sheaves (as in \cite{Ja1}).

We thank Mark Green and James Lewis for helpful conversations, and
UCLA and MPI-Bonn for support while this paper was being written.

\part{Abstract definition of the invariants}

\section{\textbf{Spread philosophy}}

Let $X$ be an $n$-dimensional smooth projective variety defined
$/\bar{\Q}$, and $k\subset\C$ be a finitely generated extension
of $\bar{\Q}$ with $trdeg(k/\bar{\Q})=:t<\infty$. We consider $0$-cycles
$\z\in Z^{n}(X_{k})$, where $X_{k}:=X\times_{(\text{Spec}\,\bar{\Q})}\text{Spec}\, k$;
our aim is to develop invariants to detect their rational equivalence
classes $\left\langle \z\right\rangle \in CH^{n}(X_{k})$. An important
tool will be the {}``$\bar{\Q}$-spread'' $\zeta$ of $\z$; we
give three versions (the first of which we shall not use):

\subsection*{Arithmetic spread}

Let $S/\bar{\Q}$ be a smooth projective variety with $\dim S=t$
and rational function field $\bar{\Q}(S)\begin{array}{c}
_{\cong}\\
\to\\
^{e}\end{array}k$. Then $e^{-1}$ induces base change to \[
\tilde{\zeta}:=\z\times_{(\text{Spec}\, k)}\text{Spec}\,\bar{\Q}(S)\,\in Z^{n}(X\times\text{Spec}\,\bar{\Q}(S)).\]
 Trivially, $\left\langle \tilde{\zeta}\right\rangle \equiv0\,\Leftrightarrow\,\left\langle \z\right\rangle \equiv0$;
we have merely applied $e^{-1}$ to the coefficients of the (homogeneous)
defining equations of $\z$ (in $k[X]$).

\subsection*{Full (geometric) spread}

A choice of hyperplane section now induces a $\bar{\Q}$-algebra homomorphism
$\bar{\Q}[S]\mapbi\bar{\Q}(S)$ of the homogeneous coordinate ring
into the function field. (If $\mathsf{{L}}$ is a linear form corresponding
to the hyperplane, set $\mathcal{{I}}(\mathsf{{F}})=\frac{\mathsf{{F}}}{\mathsf{{L}}^{d}}$
for $\mathsf{{F}}$ homogeneous of degree $d$, and extend linearly.
The induced map $\text{{Spec}}\,\bar{\Q}(S)\,\to\,\text{{Proj}}\,\bar{\Q}[S]$
gives a generic point of $S$, in the Zariski context.) Clearing denominators
from the defining equations of $\tilde{\zeta}$ therefore yields (bihomogeneous)
equations cutting out a cycle $\bar{\zeta}\in Z^{n}(X\times S)$ which
maps to $\tilde{\zeta}$ under $Z^{n}(X\times\text{Proj}\,\bar{\Q}[S])\mapbito Z^{n}(X\times\text{{Spec}}\,\bar{\Q}(S))$.

This {}``full spread'' $\bar{\zeta}$ has defects: neither the choice
of $S$ nor the clearing of denominators is unique; so the behavior
of $\bar{\zeta}$ over codimension-$1$ ($\bar{\Q}$-)subvarieties
of $S$ (where its relative dimension over $S$ may exceed $0$) is
ambiguous. In fact $\left\langle \z\right\rangle \equiv0$ (on $X_{k}$)
$\Longleftrightarrow$ $\bar{\zeta}$ is rationally equivalent to
a cycle supported on $X\times\{\text{codim.\, }1\,\subseteq\, S\}$;
to see this, simply spread or restrict a given rational equivalence.

Its two advantages: that $\bar{\zeta}$ is defined $/\bar{\Q}$ will
be essential for our use of the Bloch-Beilinson conjecture below.
Moreover, a fixed choice of $\bar{\zeta}$ may be complexified (along
with $X\times S$) and considered as a formal sum (with rational coefficients)
of irreducible analytic subvarieties of $X\times S$. This allows
us to work in the realm of Hodge theory, perform the integrals defining
the $AJ$-maps, and so on.

Throughout the paper, we will use the following notation for {}``maps''
related to $\bar{\zeta}$:

\begin{equation*} \xymatrix{& & {\bar{\zeta}} \ar [dll]_{\rho_X} \ar [ddl]^{\rho_S} \ar [dl]_{\mspace{100mu} \iota} \\ X & X \times S \ar [l]^{\pi_X} \ar [d]_{\pi_S} \\ & S } \end{equation*}${}$\\
Here $\rho_{X}$, $\rho_{S}$, and $\iota$ are really formal sums
of the corresponding maps $\rho_{X}^{i}:\,\tilde{\bar{\zeta}}_{i}\to X$
on (desingularizations of) the irreducible components of $\bar{\zeta}=\sum m_{i}\bar{\zeta}_{i}$.

\begin{rem}
These maps are both holomorphic and proper. Hence, given a (closed)
holomorphic form $\omega$ on $X$, one may pull it back along $\rho_{X}^{i}$.
Viewing the resulting holomorphic form as a current of type $(i,0)$,
one may push it forward along $\rho_{S}^{i}$ to obtain (since $\text{{codim}}_{X\times S}(\bar{\zeta})=\dim(X)$)
a closed $(i,0)$ current on $S$. By a standard $\bar{\partial}$-regularity
lemma from Hodge theory (\cite[ch. 10]{L1}), such a current is in
fact a holomorphic form. We write $\rho_{S}^{i}{}_{_{*}}\rho_{X}^{i}{}^{^{*}}\omega$
for the result and define $\rho_{S}{}_{_{*}}\rho_{X}{}^{^{*}}\omega:=\sum_{i}m_{i}\rho_{S}^{i}{}_{_{*}}\rho_{X}^{i}{}^{^{*}}\omega$,
which gives a map $\Omega^{i}(X)\to\Omega^{i}(S)$. This is used in
$\S\S11ff$.
\end{rem}

\subsection*{Final version of spread}

This will share the advantageous properties of $\bar{\zeta}$, while
lacking the ambiguities and being equivalent to $\tilde{\zeta}$ from
an algebraic point of view $/\bar{\Q}$. Namely, take $\zeta$ to
be the image of $\bar{\zeta}$ under the restriction\[
\j:\, Z^{n}(X\times S)\longrightarrow\begin{array}[t]{c}
\varinjlim\\
^{D/\bar{\Q}\,\subset\, S}\\
^{\text{codim.\, }1}\end{array}Z^{n}(X\,\times\,\, S\m D)\,=:\, Z^{n}(X\times\eta_{S}).\]
(Here $\eta_{S}$ is just the limit of the affine Zariski-open sets
$S\m D$; for our purposes, this is a much better definition of {}``generic
point'' than $\text{Spec}\,\bar{\Q}(S)$.) Clearly $\zeta\ratequiv0\Leftrightarrow\z\ratequiv0$
and\[
CH^{n}(X\times\eta_{S}{}_{/\bar{\Q}}):=\begin{array}[t]{c}
\varinjlim\\
^{D\subset S}\end{array}CH^{n}(X\,\times\,\,(S\m D))\cong CH^{n}(X\times\text{Spec}\,\bar{\Q}(S))\cong CH^{n}(X_{k}).\]

To consider $X\times\eta_{S}$ as an analytic variety (possessing
$C^{\infty}$ $(p,q)$-forms/currents, topological chains, etc.),
simply work on a (complexified) arbitrary representative $X\,\times\,\, S\m D$
of the limit. In subsequent sections, this will usually be tacit (as
will complexification). We also note that cohomology groups (singular,
absolute Hodge, etc.) of $\eta_{S}$ or $X\times\eta_{S}$ are implicitly
defined via the obvious direct limits (which are exact in the category
of abelian groups), and homology by inverse limits. 

The (very many) $\C$-points $s\in\eta_{S}(\times_{\bar{\Q}}\C)$
are those which survive in the limit, i.e. belong to the complexification
of no $D/\bar{\Q}\,\subset S$. These are {}``very general'' points
of $S$. They correspond to embeddings $e_{s}:\,\bar{\Q}(S)\mapcong k_{s}\subset\C$,
in the sense that restricting $\zeta_{\C}\left|_{_{X\times\{ s\}}}\right.=:\z_{s}^{\C}$
gives precisely the complexification of $\z_{s}:=\tilde{\zeta}\times_{\bar{\Q}(S)}k_{s}\,\in Z^{n}(X_{k_{s}})$
{[}base change via $e_{s}${]}. We write $s_{0}\in\eta_{S}^{\C}$
for the point corresponding to $e:\,\bar{\Q}(S)\mapcong k$, so that
$\z=\z_{s_{0}}$, $\z\times_{k}\C=\zeta_{\C}\left|_{_{X\times\{ s_{0}\}}}\right.$.

\section{\textbf{Background on Leray filtration}}

Given $\z\in Z^{n}(X_{k})$, the last section yields $\bar{\zeta}\in Z^{n}(X\times S_{/\bar{\Q}})$
and maps \[
Z^{n}(X_{k})\lOnto^{e\circ\mathcal{{I}}}Z^{n}(X\times S_{/\bar{\Q}})\rOnto^{\j}Z^{n}(X\times\eta_{S}{}_{/\bar{\Q}})\]
 sending $\z\lMapsto\bar{\zeta}\rMapsto\zeta$. The restriction $\j$
factors through $e\circ\mathcal{{I}}$, and so the resulting {}``spread
map''\[
Z^{n}(X_{k})\mapcong Z^{n}(X\times\eta_{S}{}_{/\bar{\Q}})\]
 sends $\z\mapsto\zeta$.

Now let $c_{\de}$ and $c_{\H}$ denote cycle-class maps into Deligne
and absolute Hodge cohomology (resp.), e.g. see \cite{L2} or \cite{Ja1}.
(Absolute Hodge keeps track of $F^{\bullet}$ and $W_{\bullet}$,
Deligne only $F^{\bullet}$. They are the same in the {[}smooth{]}
projective case. See $\S2$ of \cite{L2} also for discussion of $\text{BBC}^{q}$.)
Working modulo torsion and assuming a Bloch-Beilinson conjecture ($\text{BBC}^{q}$)
for smooth (quasi-projective) varieties, from the diagram

\begin{equation*} \xymatrix{ & CH^n(X\times S _{/\bar{\Q}})  \ar @{^{(}->} [r]^{\text{(BBC)}}_{c_{\de}} \ar @{>>} [dl] \ar [d] & H^{2n}_{\de}(X\times S, \Q(n)) \ar [d]^{\Phi} \\ CH^n(X_k)  \ar [r]^{\cong} & CH^n(X\times \eta_S {}_{/\bar{\Q}})  \ar @{^{(}->} [r]^{\text{(BBC}^q \text{)}}_{c_{\H}} & H^{2n}_{\H}(X\times \eta_S , \Q(n)) } \end{equation*}
we see that\[
CH^{n}(X_{k})\mappsi\uH^{2n}(X\times\eta_{S},\Q(n)):=\text{im}(\Phi).\]
One should think of $\text{im}(\Phi)$ as the lowest weight part of
$H_{\H}$. Lewis \cite{L2} constructs on this a Leray filtration
$\l^{i}\uH^{2n}$ (with $\l^{0}=\uH^{2n}$, $\l^{i}\supseteq\l^{i+1}$),
which induces a corresponding filtration\[
\l^{i}CH^{n}(X_{k}):=\Psi^{-1}(\l^{i}\uH^{2n})\]
 on the Chow groups. So $\Psi$ breaks naturally into {}``graded
pieces''\[
\Psi_{i}:\,\ker(\Psi_{i-1})=\l^{i}CH^{n}(X_{k})\longrightarrow Gr_{\l}^{i}\uH^{2n}(X\times\eta_{S},\Q(n));\]
$\Psi_{i}(\z)$ is well-defined for $\z\in\l^{i}CH^{n}$, and $\z\in\l^{i}$
if $\Psi_{0}(\z)=\cdots=\Psi_{i-1}(\z)=0$.

Computation of the $Gr_{\l}^{i}\uH^{2n}$ (e.g., eqn. $(4.5)$ below)
shows that they are {}``controlled'' by $H^{2n-i}(X,\begin{array}{c}
\C\\
\Q\end{array})$ (e.g., if the latter groups vanish then so do the $Gr_{\l}^{i}$).
Writing $\mathbf{H}/\bar{\Q}\subset X$ for a hyperplane section,
functoriality properties yield a diagram (see \cite[Prop. 5.0]{L2})

\begin{diagram}
\l^{n+\ell}CH^n(X_k) & \rTo^{\Psi_{n+\ell}} & Gr^{n+\ell}_{\l}\uH^{2n} & ----- & H^{n-\ell }(X) \\
\dTo>{\cdot [\mathbf{H}]^{\ell}} & & \dTo>{\cong} & & \dTo>{\text{Lefschetz}}<{\cong} \\
\begin{array}[t]{c} \underbrace{\l^{n+\ell}CH^{n+\ell}(X_k)} \\ =0 \, (\dim X = n) \end{array} & \rTo^{\Psi_{n+\ell}} & Gr^{n+\ell}_{\l} \uH^{2n+2\ell} & ----- & H^{n+\ell}(X) \\
\end{diagram}
showing $\Psi_{n+\ell}=0$ (for $\ell>1$). That is, $\l^{n+1}CH^{n}(X_{k})=\bigcap_{i}\l^{i}CH^{n}(X_{k})${[}$=0$
if $\text{BBC}^{q}$ holds{]}; so only $\Psi_{0}$ thru $\Psi_{n}$
are nontrivial maps, and one expects they completely capture rational
equivalence.

In fact, if $t\,(=trdeg(k/\bar{\Q})=\dim\, S)<n-1$ then only $\Psi_{0}$
thru $\Psi_{t+1}$ are nonzero. Writing\[
\uh^{*}(\ns,-):=\text{im}\left\{ H^{*}(S,-)\to H^{*}(\ns,-)\right\} \text{\, \, and}\]
\[
\uh_{*}(\ns,-):=\text{coim}\left\{ H_{*}(\ns,-)\to H_{*}(S,-)\right\} ,\]
the cohomology of $S$ involved in $Gr_{\l}^{i}\uH^{2n}$ is $\uh^{i}(\ns)$
and $\uh^{i-1}(\ns)$ (see $\S4$). Since $H^{t+\ell}(S)$ is concentrated
(by Lefschetz again) in codimension $\geq1$ (for $\ell\geq1$), $\uh^{t+\ell}(\ns)=0$;
therefore $Gr_{\l}^{(t+1)+\ell}\uH^{2n}=0$. Summing up this discussion
of {}``termination'' of $\l^{\bullet}$:

\begin{prop}
For $i>\text{{min}}\{ n,\, t+1\}$, we have $Gr_{\l}^{i}CH^{n}(X_{k})=0$;
if $\text{{BBC}}^{q}$ holds, then also $\l^{i}CH^{n}(X_{k})=0$.
\end{prop}
Next we explain how this picture applies to the Chow groups of all
cycles defined $/\C$. To obtain $\Psi_{i}$ for (filtered pieces
of) $CH^{n}(X_{\C})$, we first define (with \cite{L2}) $\l^{i}CH^{n}(X_{\C}):=\varinjlim_{k\subset\C}\l^{i}CH^{n}(X_{k})$
(where the limit is over $k$ finitely generated over $\bar{\Q}$).
Taking a limit over all finite-dimensional smooth projective $S/\bar{\Q}$,
one has\[
\Psi_{i}:\,\l^{i}CH^{n}(X_{\C})\to\begin{array}[t]{c}
\varinjlim\\
^{S/\bar{\Q}}\end{array}Gr_{\l}^{i}\uH^{2n}(X\times\ns,\Q(n)).\]

\begin{rem*}
The maps $Gr_{\l}^{i}\underline{H}_{\H}^{2n}(X\times\eta_{S},\Q(n))\mapjj Gr_{\l}^{i}\underline{H}_{\H}^{2n}(X\times\eta_{S'},\Q(n))$
allowed in defining this limit are only those arising from $dominant$
morphisms $S'\mapkk S$ (which correspond to inclusions $\bar{\Q}(S)\hookrightarrow\bar{\Q}(S')$).
These maps $p_{i}^{*}$ are all injective: since they depend on $S,\, S'$
only up to birational equivalence, it suffices to check this for the
cases where the extension $\bar{\Q}(S'):\bar{\Q}(S)$ is algebraic
or transcendental, which are easy.
\end{rem*}
If $H^{0}(\O_{X}^{\ell})\neq0$ for some $(n\geq)\ell\geq2$, then
the target space of $\Psi_{\ell}$ is $\infty$-dimensional, as one
would expect in connection with $Mumford's$ $theorem$. We briefly
remind the reader what this says, and set a couple of conventions
for later use. For any $N\in\mathbb{N}$ and base point $p\in X$,
there is a map\[
S^{(N)}(X_{\C})\mappiN CH_{hom}^{n}(X_{\C})\]
\[
Sym_{N}(p_{1},\ldots,p_{N})\longmapsto\sum_{j=1}^{N}\left\langle p_{j}\right\rangle \,-N\left\langle p\right\rangle .\]
One says $CH_{[hom]}^{n}(X_{\C})$ is {}``infinite-dimensional''
(over $\C$) if the sequence $d_{N}=\dim_{\C}(\text{im}(\pi_{N}))$
is unbounded (given a suitable definition of dimension for this image,
see \cite[Ch. 15]{L1}). The following consequence of \cite[Prop. 15.3]{L1}
is nontrivial:

\begin{prop}
$\dim(CH^{n}(X_{\C}))=\infty\,\,\Longleftrightarrow\,\,$no $\pi_{N}$
is surjective.
\end{prop}
Mumford's theorem {[}as generalized by Roitman{]} states that $H^{0}(\O_{X}^{\ell})\neq0$
for some $\ell\geq2$ $\implies$ $d_{N}\geq N$ $\implies$ $[Gr_{\l}^{\ell}]CH^{n}(X_{\C})$
is $\infty$-dimensional.

We conclude the section with an important description of the groups
$\underline{H}^{i}(\eta_{S},\Q)$ introduced above. Note that well-definedness
is an issue, since different $S$'s may have the {}``same $\eta_{S}$''.

\begin{prop}
$\underline{H}^{i}(\eta_{S},\Q)$ is a (pure, finite-dimensional)
Hodge structure. In particular, it is well-defined and equal to $W_{i}H^{i}(\eta_{S},\Q):=\varinjlim W_{i}H^{i}(U,\Q)$
(limit over $U\subseteq S$ affine Zariski-open).
\end{prop}
\begin{proof}
First let $U\mapiota S$ be affine Zariski open. By \cite{H}, there
exists an embedded resolution of $S\m U$. This yields $\hat{S}\mapb S$,
with $b$ a finite composition of blow-ups at nonsingular centers
contained in the successive preimages of $S\m U$, such that $b^{-1}(U)\cong U$
and $\hat{S}\m b^{-1}(U)$ is a normal-crossings divisor (NCD). Since
$H^{i}(\hat{S})\to H^{i}(b^{-1}(U))\cong H^{i}(U)$ factors as $\iota^{*}\circ b_{*}$
(or from the description of cohomology of a blow-up), $\text{{im}}\{ H^{i}(S)\to H^{i}(U)\}=\text{{im}}\{ H^{i}(\hat{S})\to H^{i}(U)\}$;
and by \cite[cor. 3.2.17]{De1} (since $(\hat{S}\m b^{-1}(U))=$NCD)
the latter is $W_{i}H^{i}(U)$. Hence $\text{{im}}\{ H^{i}(S)\to H^{i}(U)\}=W_{i}H^{i}(U)$,
and taking a limit gives $\underline{H}^{i}(\eta_{S})=W_{i}H^{i}(\eta_{S})$.

Well-definedness is a trivial consequence of this descrption, but
this uses well-definedness of the $W_{i}H^{i}(U)$. One can argue
more directly as follows: if for some affine $U_{0}$, $S_{1}\mapaa U_{0}\mapbb S_{2}$
are Zariski-open embeddings (that is, if $\eta_{S_{1}}\cong\eta_{S_{2}}$),
then it suffices to show $\text{{im}}(\iota_{1}^{*})=\text{{im}}(\iota_{2}^{*})$
in $H^{i}(U_{0})$. Now simply take the closure (in $S_{1}\times S_{2}$)
of the {}``diagonal'' $(\iota_{1}\times\iota_{2})$-image of $U_{0}$,
and write $S_{0}$ for a desingularization. The obvious morphisms
$S_{0}\twoheadrightarrow S$, $S_{0}\twoheadrightarrow S_{2}$ then
induce surjective push-forwards on cohomology factoring the restriction
$\iota_{0}^{*}$ to $H^{i}(U_{0})$, and so $\text{{im}}(\iota_{1}^{*})=\text{{im}}(\iota_{0}^{*})=\text{{im}}(\iota_{2}^{*})$.
\end{proof}
\begin{rem*}
The reader may also wonder about well-definedness of $\underline{H}_{\H}=\text{{im}}(\Phi)$;
this follows from remark $4.2$ below.
\end{rem*}

\section{\textbf{Description of the invariants}}

Let $\sh_{\Q}$ be a MHS($=$rational mixed Hodge structure), and
write $\sh_{\C}=\sh_{\Q}\t_{\Q}\C,$ $W_{0}\sh_{\C}=(W_{0}\sh_{\Q})\t\C$
and so on. Then (see \cite[sec. 9]{Ja2})\[
Hom_{_{\text{MHS}}}(\Q(0),\sh_{\Q})=F^{0}W_{0}\sh_{\C}\cap W_{0}\sh_{\Q},\,\,\,\, Ext_{_{\text{MHS}}}^{1}(\Q(0),\sh_{\Q})=\frac{W_{0}\sh_{\C}}{W_{0}\sh_{\Q}+F^{0}W_{0}\sh_{\C}}\]
so that

\begin{equation} Hom_{_{\text{MHS}}} (\Q(0), \sh_{\Q} \t \Q(n)) = \left\{ F^nW_{2n} \sh_{\C} \cap W_{2n} \sh_{\Q} \right\} (n), \end{equation} \begin{equation} Ext^1_{_{\text{MHS}}} (\Q(0), \sh_{\Q}\t \Q(n)) = \left\{ \frac {W_{2n} \sh_{\C}} {W_{2n} \sh_{\Q} + F^nW_{2n} \sh_{\C} } \right\} (n). \end{equation}

${}$\\
We will omit the $(n)$ frequently when it is not important.

To relate $\Psi_{i}(\z)$ to geometry we first fix notions of cycle-class
and Abel-Jacobi class for $\zeta$. There is a diagram of short-exact
sequences

\SMALL \begin{equation} \begin{diagram}  \begin{array}[b]{c} J^n(X\times S) \\ \parallel \\ Ext^1_{_{\text{MHS}}} \left( \Q(0), H^{2n-1}(X\times S, \Q (n)) \right) \end{array} & \rTo^q & H^{2n}_{\de}(X\times S ,\Q(n)) & \rTo^p & \begin{array}[b]{c} Hg^n(X\times S) \\ \parallel \\ Hom_{_{\text{MHS}}} \left( \Q(0) , H^{2n} (X\times S ,\Q(n)) \right) \end{array} \\  \dTo>{''\Phi} & & \dTo>{\Phi} & & \dOnto>{' \Phi} \\ Ext^1_{_{\text{MHS}}} \left( \Q(0), H^{2n-1}(X\times \ns ,\Q(n)) \right) & \rTo^q & H^{2n}_{\H} (X \times \ns , \Q(n)) & \rTo^p & Hom_{_{\text{MHS}}} \left( \Q(0) , H^{2n}(X \times \ns ,\Q(n)) \right) \\  \end{diagram} \end{equation} \normalsize
In fact $\text{{im}}({}''\Phi)\mapq\text{{im}}(\Phi)\mapp\text{{im}}({}'\Phi)$
is also short-exact (since $\ker(\Phi)\rOnto\ker({}'\Phi)$, see \cite[p. 317]{L2}),
and equals

\begin{equation} \begin{diagram} \begin{matrix} Hom_{_{\text{MHS}}} \left( \Q(0) , H^{2n} (X \times \ns ,\Q(n)) \right) & =: \underline{Hg}^n(X \times \ns) \\ \uparrow \, p \\ \uH^{2n} (X \times \ns , \Q(n)) \\
\uparrow \, q \\ \frac{ Ext^1_{_{\text{MHS}}} \left( \Q(0), W_{-1} H^{2n-1} (X \times \ns , \Q(n)) \right ) } { Hom_{_{\text{MHS}}} \left( \Q(0), Gr^W_0 H^{2n-1} (X\times \ns ,\Q(n)) \right) } & =: \underline{J}^n (X \times \ns ) . \end{matrix} \\ \end{diagram} \end{equation} 
Recalling that $c_{\H}(\zeta)=\Phi(c_{\de}(\bar{\zeta}))$, we define
the fundamental class of the spread to be\[
[\zeta]:=p(c_{\H}(\zeta))={}'\Phi[p(c_{\de}(\bar{\zeta}))]={}'\Phi[\bar{\zeta}]\,\in\,\underline{Hg}^{n}(X\times\ns).\]
 If this is zero, define the $AJ$-class in $\underline{J}^{n}(X\times\ns)$
by $q[AJ\zeta]=c_{\H}(\zeta)$. 

To compute $[AJ\zeta]$, one observes (using \cite[cor. 8.2.8]{De2})
that\[
[\zeta]=0\,\,\implies\,\,[\bar{\zeta}]\in\ker({}'\Phi)\,=\,\begin{array}[t]{c}
\bigcup\\
^{D/\bar{\Q}\subset S}\\
^{\text{codim.\, }1}\end{array}\text{im}\left\{ Hg^{n-1}(X\times\widetilde{D})\right\} \,;\]
 according to the Hodge conjecture one then expects to be able to
modify $\bar{\zeta}$ {}``over $D$'' by $\eta\in N_{S}^{1}Z^{n}(X\times S):=\bigcup_{D/\bar{\Q}}\text{im}\{ Z^{n-1}(X\times D)\}$
so that $[\bar{\zeta}-\eta]=0$ (of course $\eta\in\ker(\Phi)$).
Then $[AJ\zeta]={}''\Phi[AJ(\bar{\zeta}-\eta)]$, which can be computed
by the usual Griffiths prescription (on $X\times S$). Namely, one
integrates over a topological membrane (i.e., $C^{\infty}$ chain)
bounding on $\bar{\zeta}-\eta$ (of real dimension $2t+1$, written
$\partial^{-1}(\bar{\zeta}-\eta)$). This produces a functional on
$F^{t+1}H^{2t+1}(X\times S,\C)$, which taken modulo periods gives
$AJ(\bar{\zeta}-\eta)$. In particular, if $\dim(S)=1$ and $\dim(X)=2$,
then the HC on $X\times\widetilde{D}$ is just Lefschetz (1,1), so
such an approach is practical. Moreover, there are plenty of cases
(for higher dimensional $X$ and $S$) where one can explicitly write
down a homologically trivial choice of full spread $\bar{\zeta}$
(and avoid the HC).

\begin{rem}
In general there is no circumventing this {}``lift'' to {[}a nullhomologous{]}
$\bar{\zeta}$ for computing $[AJ\,\zeta]$ and its graded pieces
(defined below) by membrane integrals. While from $\zeta\homeq0$
one $can$ obtain a chain $\partial^{-1}\zeta$ bounding on $\zeta$
mod {[}chains w./ support on{]} $X\times D$, the naive idea that
this should lead to a functional on $\varprojlim_{D}F^{t+1}H^{2t+1}(X\times S,X\times D;\,\C)$
{[}modulo periods{]} and ultimately an element of $\underline{J}^{n}(X\times\eta_{S})$,
does not work out for $X$ smooth projective.

The situation is different if $X$ is one of the special $relative$
varieties of Part $2$, due to the occurrence of a Hodge-theoretic
miracle (see $\S\S7-8$). The miracle (or, if one prefers, catastrophe)
produces a very convenient set of test forms (for $\varprojlim_{D}F^{t+1}H^{2t+1}(X\times(S,D),\C)$)
to integrate over (essentially) $\partial^{-1}\zeta$, namely products
of {[}Poincar\'e duals to{]} certain topological cycles on $S$ and
{}``holomorphic'' forms on $X$. Integrating only such forms for
ordinary projective $X$ computes (in contrast) only quotients of
the graded pieces of $[AJ\,\zeta]$; but at least these quotients
are computable by membrane integrals without requiring $[\bar{\zeta}]=0$.
(The integrals take place on $X$ rather than on $S$.) This is the
idea behind the $reduced$ higher $AJ$ maps of $\S11$. 
\end{rem}
Now referring back to eqn. $(4.4)$, the $Ext^{1}/Hom$ quotient (in
$\underline{J}^{n}(X\times\ns)$) is computed by observing that there
is a natural series of maps\[
Hom_{_{\text{MHS}}}\left(\Q(0),\, Gr_{0}^{W}H^{2n-1}(X\times\ns,\Q(n))\right)(-n)\,\cong\]
\[
\ker\left\{ F^{n}Gr_{2n}^{W}H^{2n-1}(X\times\ns,\C)\oplus Gr_{2n}^{W}H^{2n-1}(\ns,\Q)\to Gr_{2n}^{W}H^{2n-1}(X\times\ns,\C)\right\} \]
\small\[
\twoheadrightarrow\frac{\left\{ F^{n}W_{2n}H^{2n-1}(X\times\ns,\C)+W_{2n}H^{2n-1}(X\times\ns,\Q)\right\} \cap\left\{ W_{2n-1}H^{2n-1}(X\times\ns,\C)\right\} }{\text{im}\left\{ F^{n}W_{2n-1}H^{2n-1}(X\times\ns,\C)\oplus W_{2n-1}H^{2n-1}(X\times\ns,\Q)\to W_{2n-1}H^{2n-1}(X\times\ns,\C)\right\} }\]
\normalsize\[
\longrightarrow\frac{W_{2n-1}H^{2n-1}(X\times\ns,\C)}{F^{n}W_{2n-1}H^{2n-1}(X\times\ns,\C)+W_{2n-1}H^{2n-1}(X\times\ns,\Q)}\]
\[
\cong\, Ext_{_{\text{MHS}}}^{1}\left(\Q(0),\, W_{-1}H^{2n-1}(X\times\ns,\Q(n))\right)(-n).\]
(Note that $W_{2n-1}H^{2n-1}(X\times\ns)$ is just $\text{im}\left(H^{2n-1}(X\times S)\to H^{2n-1}(X\times\ns)\right)$.)

Now let $\eta_{0}\in Z_{hom}^{n}(X\times S)$ be a homologically trivial
cycle supported on $X\times D$; $\text{im}(Hom_{_{\text{MHS}}})$
should contain (and thus annihilate) the image of $AJ(\eta_{0})$
under $J^{n}(X\times S)\to Ext_{_{\text{MHS}}}^{1}\left(\Q(0),\, W_{-1}H^{2n-1}(X\times\ns,\Q(n))\right).$
(Intuitively, this is why we quotient it out.) The delta-function
current%
\footnote{See \cite{GH,Ja1} for discussion of currents.%
} $\delta_{\eta_{0}}\in Z_{\di}([F^{n}]\DI_{(X\times S)^{\infty}}^{2n})$
has two primitives: $\omega\in\Gamma(F^{n}\DI_{(X\times S)^{\infty}}^{2n-1})$
and $-\delta_{\KK}$, $\KK\in C_{2t+1}^{top}(X\times S)$ a $C^{\infty}$
chain. Their restrictions ($\j^{*}$) to $X\times\ns$ are $\di$-closed,
and give classes in $F^{n}W_{2n}H^{2n-1}(X\times\ns,\C)$ and $W_{2n}H^{2n-1}(X\times\ns,\Q)$
(resp.); but $\omega+\delta_{\KK}$ is closed on $X\times S$, and
thus has class in $W_{2n-1}H^{2n-1}(X\times\ns,\C).$ So $\j^{*}[\omega+\delta_{\KK}]\in\text{im}(Hom_{_{\text{MHS}}})$,
on the one hand. But the image of $[\omega+\delta_{\KK}]\in H^{2n-1}(X\times S)$
in $J^{n}(X\times S)\cong\left\{ F^{t+1}H^{2t+1}(X\times S,\C)\right\} ^{\vee}\left/\text{(periods)}\right.$
is just $\left[\int_{\KK}(\cdot)\right]$ (where $\d\KK=\eta_{0}$),
which is Griffiths' prescription for $AJ(\eta_{0})$.

Now $\Psi_{i}(\z)$ is a class in the middle term of the short-exact
sequence

\begin{equation} \begin{matrix} Hom_{_{\text{MHS}}} \left( \Q(0) , H^i (\ns ,R^{2n-i}{\pi_S \, _{_*}} \Q(n)) \right) & =Gr^{i}_{\l} \underline{Hg}^n (X\times \ns) \\ \uparrow \, p \\ Gr^i_{\l} \uH^{2n} (X\times \ns ,\Q(n)) \\ \uparrow \, q \\ \frac{Ext^1_{_{\text{MHS}}} \left( \Q(0) , W_{-1} H^{i-1} (\ns , R^{2n-i} {\pi_S \, _{_*} } \Q(n)) \right)} { Hom_{_{\text{MHS}}} \left( \Q(0) , Gr^W_0 H^{i-1} (\ns ,R^{2n-i} {\pi_S \, _{_*}} \Q(n)) \right) } & = Gr^{i-1}_{\l} \underline{J}^n (X\times \ns ) , \end{matrix} \end{equation}

${}$\\
where the $R^{2n-i}\pi_{S}{}_{_{*}}\Q(n)\cong H^{2n-i}(X,\Q(n))$
are constant sheaves. 

\begin{rem*}
\cite{L2} has $Hom_{_{\text{MHS}}}\left(\Q(0),\, Gr_{0}^{W}H^{2n-1}(X\times\ns,\Q(n))\right)$
in the bottom left-hand denominator of $(4.5)$; our simplification
follows easily from the K\"unneth decomposition of $H^{*}(X\times\ns)$
and purity of $H^{*}(X).$
\end{rem*}
Define $[\zeta]_{i}:=p(\Psi_{i}(\z))$, and if this is zero $[AJ\zeta]_{i-1}:=q^{-1}(\Psi_{i}(\z)).$
These should be thought of intuitively as {}``Leray graded pieces''
of the cycle/$AJ$-classes of the spread. One may take this statement
literally for $[\zeta]_{i}$, but not for $[AJ\zeta]_{i}$ since $[AJ\zeta]$
may not be defined (unless all $[\zeta]_{j}$, and thus $[\zeta]$,
are zero). In fact, for $[\zeta]_{i}$ to be defined, one need only
have $[\zeta]_{0}=\cdots=[\zeta]_{i-1}=0$ (not $\Psi_{0}(\z)=\cdots=\Psi_{i-1}(\z)=0$).

One may also think of $\Psi_{i}(\z)$ as a {}``higher normal function'',
namely an $H^{2n-i}(X,\C)$-valued $(i-1)$-current on $S$, with
infinitesimal invariant ($i$-current) $[\zeta]_{i}$. If this is
zero then the resulting closed $(i-1)$-current yields $[AJ\zeta]_{i-1}$;
see $(10.4)$. Moreover, in $\S11$ {[}remark $(11.4)${]} it is shown
that an image of $[\zeta]_{i}$ identifies with the $i^{\text{th}}$
Mumford invariant of S. Saito.

\begin{rem}
Lewis proves $(4.4)$ and $(4.5)$ with $\eta_{S}$ everywhere replaced
by an affine open $U\subseteq S$ with $S\m U$ a NCD, and $\underline{H}_{\H}$
replaced by the image of $\Phi_{S,U}:\, H_{\de}^{2n}(X\times S,\Q(n))\to H_{\H}^{2n}(X\times U,\Q(n))$.
For the $\eta_{S}$-versions, it is enough to show that the NCD condition
can be removed (then take a limit). If $S\m U$ is not a NCD, an embedded
resolution (as in proof of Prop. $3.3$) yields $\hat{S}\mapbproj S$,
$(\hat{S}\m b^{-1}(U))\cong U$ a NCD, and $H^{*}(\hat{S})\mapbstar H^{*}(S)$.
Using the K\"unneth decomposition and the short-exact sequence for
$H_{\de}^{2n}$ in $(4.3)$, $H_{\de}^{2n}(X\times\hat{S},\Q(n))\mapbde H_{\de}^{2n}(X\times S,\Q(n))$;
and so $\Phi_{S,U}\circ b_{*}^{\de}=\Phi_{\hat{S},U}$ $\implies$
$\text{{im}}(\Phi_{S,U})=\text{{im}}(\Phi_{\hat{S},U})$. Everything
else in the $U$-versions of $(4.4-5)$ is intrinsic to $U$, so we
are done.
\end{rem}

\section{\textbf{Some special cases}}

Now we shall compute the terms of $(4.5)$ for $i=0,1$ and determine
the corresponding invariants (to be continued for $i\geq2$ in $\S\S10$ff).
For $i=0$ one has (tacitly twisting by $(-n)$), since $Gr_{\l}^{-1}\underline{J}^{n}=0$,
\[
0\to Gr_{\l}^{0}\uH^{2n}(X\times\ns,\Q(n))\begin{array}{c}
_{\cong}\\
\to\\
^{p}\end{array}\uh^{0}(\ns,Hg^{n}(X)),\]
 where $Hg^{n}(X)$ is regarded as constant sheaf. Now $\Psi_{0}(\z)$
is the function $\ns\to Hg^{n}(X)$ given by $s\mapsto[\z_{s}]$;
since $\Psi_{0}(\z)=p^{-1}[\zeta]_{0}$, it must be constant (this
we owe essentially to rigidity of $H_{\text{DR}}^{*}$) and equal
to $[\z_{s_{0}}]=[\z]$. So $\Psi_{0}(\z)=[\z]$ and $\l^{1}CH^{n}(X_{k})=CH_{hom}^{n}(X_{k})$.

Next consider $i=1$: noting that the $Hom$-denominator in the left-hand
term is zero (as $Gr_{1}^{W}H^{0}(\ns)=0$), one gets\[
\uh^{0}(\ns,J^{n}(X))\mapq Gr_{\l}^{1}\uH^{2n}(X\times\ns,\Q(n))\mspace{250mu}\]
\vspace{-5mm} \[
\mspace{150mu}\mapp F^{n}\left(\uh^{1}(\ns,\C)\t H^{2n-1}(X,\C)\right)\cap\left(\uh^{1}(\ns,\Q)\t H^{2n-1}(X,\Q)\right).\]
Write $V:=H^{2n-1}(X,\C)\supseteq H^{2n-1}(X,\Q)=:\Lambda$, and note
that $V=F^{n-1}V$, $J^{n}(X)=\frac{F^{n-1}V}{F^{n}V+\Lambda}$. As
will be clear from $\S11$ (using $(11.1)$ together with the observation
that $\chi_{1}=\Psi_{1}$), $\Psi_{1}(\z)$ is the holomorphic function
$\sigma:\,\ns\to J^{n}(X)$ given by $s\mapsto AJ(\z_{s})$. This
lifts to a $discontinuous$ function $\ns\to V$, that is, a $V$-valued
$0$-current $\tilde{\sigma}\in\Gamma(\DI_{S}^{0}(V))$. Indeed, $Gr_{\l}^{1}\uH^{2n}$
may be identified with 

\begin{equation} \frac{\ker \left\{ \Gamma(\DI^0_S (V)) \mapd \frac{Z(\DI^1_S (V))} {Z(C^S_{2t-1}(\Lambda)) + Z(F^n \DI^1_S (V))} \right\} } {\Gamma (C^S_{2t}(\Lambda)) + \Gamma (\DI^0_S (F^n V)) } , \end{equation}

${}$\\
where $C_{*}^{S}(\Lambda)$ are (sheafified) $\Lambda$-valued $C^{\infty}$
$*$-chains, $Z(\,\cdot\,)$ means closed sections, and $F^{n}\DI_{S}^{1}(V)=\DI_{S}^{1}(F^{n}V)+F_{S}^{1}\DI_{S}^{1}(V)$.
(Holomorphicity of $\sigma$ $\implies$ $\di[\tilde{\sigma}]$ has
no $\frac{F_{S}^{0}}{F_{S}^{1}}\DI_{S}^{1}\left(\frac{V}{F^{n}V}\right)$
component.) Hence $\di[\tilde{\sigma}]=c+\omega$ gives classes $[c]\in\uh^{1}(\ns,\Q)\t H^{2n-1}(X,\Q)$
and $[-\omega]\in F^{n}\left(\uh^{1}(\ns,\C)\t H^{2n-1}(X,\C)\right)$
which are equivalent, hence the {}``$\cap$'', hence $[\zeta]_{1}=[-\omega]=[c].$
If this class is zero, modifying $\tilde{\sigma}$ by primitives in
the denominator of $(5.1)$ yields $'\tilde{\sigma}$ with $\di[\,'\tilde{\sigma}]=0$,
thus a constant function from $\ns\to V$ $also$ $lifting$ $\sigma$.
So $AJ(\z_{s})=AJ(\z_{s_{0}})$ (for all $s\in\ns$), and $\Psi_{1}(\z)$
lifts to $AJ(\z)\in\uh^{0}(\ns,\, J^{n}(X))$. $\underline{\text{Conclusion}}$:
if $[\zeta]_{1}=0$ then $[AJ\zeta]_{0}$ is defined and $=AJ(\z)$.

On the other hand, $AJ(\z)$ is defined whether or not $[\zeta]_{1}=0$,
and so it is a nontrivial question to ask:\[
\text{if\, }AJ(\z_{(s_{0})})=0\text{,\, must\, }AJ(\z_{s})=0\text{\, for\, all\, }s\text{?}\]
While the $\{\z_{s}\}$ are algebraically isomorphic, they have different
transcendental geometry, and $AJ(\z_{s})$ is an invariant of that
geometry. That the answer is yes, should be considered fortunate;
the same question for $m$-cycles ($m>0$) is open. (The answer is
expected to be affirmative for cycles algebraically equivalent to
zero. This is consistent with the case of $0$-cycles, since all $0$-cycles
are $\algequiv0$.) 

A statement obviously equivalent to the affirmative answer for $0$-cycles,
which we now prove, reads: $AJ(\z)=0$ $\implies$ $\Psi_{1}(\zeta)=0$.

\begin{proof}
We temporarily write {}``$Alb(X)$'' for $J^{n}(X)$ and {}``$Alb$''
for $AJ$, to agree with notation in \cite{M} ($AJ$ will be written
only for curves).

While $Alb$-equivalence does not spread $a$ $priori$, rational
equivalence is algebraic and $does$. For $0$-cycles on curves, $AJ(\W)=0$
$\implies$ $\W\ratequiv0$, so we want to factor a correspondence
giving the identity on $Alb(X)$ through a curve. This is done in
\cite{M}; namely, there exists a smooth curve $C/\bar{\Q}\subseteq X$
and correspondence $^{t}\tilde{E}\in CH^{1}(X\times C_{/\bar{\Q}})$
so that the top (and $\therefore$ bottom) compositions in the following
diagram are isomorphisms:

\begin{diagram}
Alb(X) & \rTo & J(C) & \rTo & Alb(X) \\
\uTo>{Alb} & & \uTo<{\cong}>{AJ} & & \uTo>{Alb} \\
\l^1 CH^n (X_k) & \rTo^{^t \tilde{E}} & \l^1 CH^1 (C_k) & \rTo^{\i_C} & \l^1 CH^n (X_k) \\
\dTo>{\Psi_1} & & \dTo>{\Psi_1} & & \dTo>{\Psi_1} \\
Gr^1_{\l} \uH^{2n}(X\times \ns) & \rTo^{[ \, ^t \tilde{E} \times \Delta_S ]_*} & Gr^1_{\l} \uH^2 (C \times \ns) & \rTo^{[\i_C \times \Delta_S ]_*} & Gr^1_{\l} \uH^{2n} (X\times \ns) \\ 
\end{diagram} \vspace{-12mm}

\begin{diagram}
{} & & & & \rTo_{\cong} & & & & {} \\
\end{diagram}
(where we are identifying $Gr_{\l}^{1}\uH^{2n}(X\times\ns)$ with
{[}holomorphic{]} functions from $\ns\to Alb(X)$). So given $\z\in\l^{1}CH^{n}(X_{k})$
with $Alb(\z)=0$, one traces through the diagram to see $^{t}\tilde{E}\cdot\z\ratequiv0$
$\implies$ $\Psi_{1}({}^{t}\tilde{E}\cdot\z)=0$ $\implies$ $\Psi_{1}(\i_{C}\cdot{}^{t}\tilde{E}\cdot\z)=0$
$\implies$ $\Psi_{1}(\z)=0$.

{[}Geometrically speaking, applying $^{t}\tilde{E}\times\Delta_{S}$
to $\zeta$ turns the family $\{\z_{s}\}$ into a family $\{{}'\z_{s}\}$
of $0$-cycles supported on $C$, with $Alb_{X}({}'\z_{s})=Alb_{X}(\z_{s})$
($\forall\, s$). Now $Alb_{X}(\z_{s_{0}})=0$ $\implies$ $AJ_{C}({}'\z_{s_{0}})=0$
$\implies$ $'\z_{s_{0}}\ratequiv0$ $\implies$ $'\z_{s}\ratequiv0$
$\implies$ $Alb_{X}({}'\z_{s})=0$ $\implies$ $Alb_{X}(\z_{s})=0$.{]}
\end{proof}
\begin{rem*}
(i) $\frac{1}{n}\i_{C}\cdot{}^{t}\tilde{E}\,\in\, CH^{n}(X\times X)$
is called the Albanese projector.

(ii) Owing to the last result we shall call the $\Psi_{i}$ {}``higher
Abel-Jacobi maps''.
\end{rem*}

\section{\textbf{Relation to the work of Griffiths and Green}}

We want to point out the equivalence of the approach taken here (based
on work of Lewis) with that taken in \cite{GG3}, which considers
graded pieces of $c_{\de}(\bar{\zeta})$ (modulo ambiguities) in lieu
of $c_{\H}(\zeta)\,[=\Psi(\z)]$. (This uses a Leray filtration on
$H_{\de}^{2n}(X\times S,\Q(n))$; see the discussion in $\S10$.)
What we do here is briefly explain, in parallel with the discussion
of termination of $\l^{\bullet}$ in $\S3$, why the Griffiths-Green
approach also leads to a filtration $\mathcal{{F}}_{GG}^{\bullet}$
on $CH^{n}$ that terminates at (or before ) the $(n+1)^{\text{{st}}}$
step. A complete proof (modulo HC) of the equivalence of $\mathcal{{F}}_{GG}^{\bullet}$
and $\l^{\bullet}$ has appeared in \cite[sec. 1.6]{mS} for $X$
defined $/\bar{\Q}$. 

A cycle $\z\in CH^{n}(X_{k})$ lies in $\mathcal{{F}}_{GG}^{i}$ if
there $exists$ a choice of complete $\bar{\Q}$-spread $\bar{\zeta}$
such that $c_{\de}(\bar{\zeta})\in\l_{S}^{i}H_{\de}^{2n}(X\times S,\Q(n))$.
In order to exhibit the termination property $\mathcal{{F}}_{GG}^{n+1}=\mathcal{{F}}_{GG}^{n+k}$
for all $k>0$, the following will suffice: given $\bar{\zeta}\in Z^{n}(X\times S_{/\bar{\Q}})$
with $c_{\de}(\bar{\zeta})\in\l_{S}^{n+1}H_{\de}^{2n}$, we must produce
$\Gamma_{\de}\in Z^{n}(X\times D_{/\bar{\Q}})$ (for some divisor
$D_{/\bar{\Q}}\subseteq S$) such that $c_{\de}(\bar{\zeta}-\Gamma_{D})=0$.

The construction of this $\Gamma_{D}$ involves the Hodge conjecture
(HC), which never explicitly appears in Lewis; that is because it
is absorbed into $\text{{BBC}}^{q}$ (in fact, HC$+$BBC$\,\implies\,\text{{BBC}}^{q}$).
That is, in the Lewis approach it is a longer leap from vanishing
of the invariant ($c_{\H}(\zeta)=0$) to rational equivalence $(\zeta\ratequiv0$),
while in the Griffiths-Green approach vanishing of the invariant ($c_{\de}(\bar{\zeta})=0$)
is harder.

One produces $\Gamma_{D}$ one bit at a time (modifying $\bar{\zeta}$
as we go), first picking off $[\bar{\zeta}]_{n+1}$ thru $[\bar{\zeta}]_{2n}$,
then (since now $\bar{\zeta}'\homeq0$) $[AJ(\bar{\zeta}')]_{n}$
thru $[AJ(\bar{\zeta}')]_{2n-1}$. As the idea is the same for each
of these, we just demonstrate it (for one of the $[\bar{\zeta}]_{n+i}$)
in an example.

For this example, we take $X$ to be a $3$-fold ($\z$ a $0$-cycle
as usual), and of course $Gr_{\l_{S}}^{0}c_{\de}(\bar{\zeta})=\cdots=Gr_{\l_{S}}^{3}c_{\de}(\bar{\zeta})=0$.
We take aim at $[\bar{\zeta}]_{4}$, the graded piece of the fundamental
class of $\bar{\zeta}$ in \[
Hg^{3}(X\times S)\cap\{ H^{2}(X)\otimes H^{4}(S)\}.\]
The idea is to take a {}``cross-section'' of $\bar{\zeta}$ by a
hyperplane section $X_{\mathbf{{H}}}$ of $X$ defined $/\bar{\Q}$;
$\bar{\zeta}\cap(X_{\mathbf{{H}}}\times S)$ is then supported on
$X_{(\mathbf{{H}})}\times D$ for some codim.-$1$ $D_{/\bar{\Q}}\subseteq S$.
This cross-section map may be formulated as a correspondence \[
\left\langle \mathbf{H}\times\Delta_{S}\right\rangle \in CH^{*}\left((X\times S)\times(X\times S)\right),\]
 so that the isomorphism of hard Lefschetz is induced by its action:\[
H^{2}(X)\t H^{4}(S)\begin{array}{c}
_{[\mathbf{H}\times\Delta_{S}]_{*}}\\
\longrightarrow\\
^{\cong}\end{array}H^{4}(X)\t H^{4}(S).\]
 According to the Lefschetz standard conjecture ($\Leftarrow$ HC
+ hard Lefschetz), there exists a correspondence $\left\langle \mathbf{H}^{-1}\right\rangle $
(algebraically) inducing its inverse (and inducing the zero map in
other degrees%
\footnote{This is not part of the Lefschetz standard conjecture as such, but
can be arranged without additional assumptions by \cite{Kl}.%
}). That is, $[\mathbf{H}^{-1}]_{*}=\left([\mathbf{H}]_{*}\right)^{-1}$,
so the composition $[\mathbf{H}^{-1}\times\Delta_{S}]_{*}\circ[\mathbf{H}\times\Delta_{S}]_{*}=\text{id}$
on cohomology $H^{2}(X)\otimes H^{4}(S)$ (and zero on $H^{i}(X)\otimes H^{6-i}(S)$
for $i\neq2$). But because it is {}``algebraic'' (i.e. induced
by a correspondence), the composition operates on cycles, and is $not$
the identity on $\bar{\zeta}$. In particular, since already $\left\langle \mathbf{H}\times\Delta_{S}\right\rangle \cdot\bar{\zeta}$
is supported on $X\times D$, \[
\Gamma_{D}^{4}\,:=\,\left\langle \mathbf{H}^{-1}\times\Delta_{S}\right\rangle \left(\left\langle \mathbf{H}\times\Delta_{S}\right\rangle \cdot\bar{\zeta}\right)\,\subseteq\, X\times D.\]
As cohomology classes, though, \[
[\Gamma_{D}^{4}]=[\mathbf{H}^{-1}\times\Delta_{S}]_{*}\circ[\mathbf{H}\times\Delta_{S}]_{*}[\bar{\zeta}]\,=\,[\bar{\zeta}]_{4}\,\in\, H^{2}(X)\t H^{4}(S)\]
 and so the modification $\bar{\zeta}-\Gamma_{D}^{4}$ kills $Gr_{\l_{S}}^{4}$
of the fundamental class (of $\bar{\zeta}$) without affecting either
the other $Gr_{\l_{S}}^{i}$ or the complete spread's restriction
to $\z=\bar{\zeta}\cap(X\times\{ s_{0}\})$ (or even to $\zeta$).
Griffiths and Green say such a class is {}``in the ambiguities''.

\part{Motivation via degenerate cases}

\section{\textbf{Formal computations for relative varieties}}

The results of Part $1$ are stated for $X$ smooth projective $/\bar{\Q}$.
We now propose to apply them formally to $0$-cycles on certain $relative$
$quasiprojective$ $X$. It turns out that the filtrations and maps
they predict (for these special $X$) are correct, explicitly constructible,
and give insight into the case of $X$ projective -- which we resume
in Part $3$.

Our first special $X=(\s_{\bar{\Q}}^{n},\d\s_{\bar{\Q}}^{n})$, the
{}``algebraic $n$-cube'', where $\s_{\bar{\Q}}^{n}:=(\P_{\bar{\Q}}^{1}\m\{1\})$
and $\d\s_{\bar{\Q}}^{n}:=\bigcup_{j}\left\{ (z_{1},\ldots,z_{n})\in\s_{\bar{\Q}}^{n}\,\left|\, z_{j}=0\text{\, or\, }\infty\right.\right\} .$
We consider $0$-cycles $\z\in CH^{n}(\s_{k}^{n},\d\s_{k}^{n})$,
$k\cong\bar{\Q}(S)$. We will write $(\s^{n},\d\s^{n})$ for the corresponding
analytic variety $/\C$, and $\left((\C^{*})^{n},\,\BI^{n}\right)$
for the {}``dual'' variety, where $\BI^{n}=\bigcup_{j}\left\{ (z_{1},\ldots,z_{n})\in(\C^{*})^{n}\,\left|\, z_{j}=1\right.\right\} .$
(One can also think of these two as $(\P^{1}\m\{1\},\,\{0,\infty\})^{n}$
and $(\C^{*},\{1\})^{n}$.)

There is a perfect pairing\[
H^{2n-i}\left((\s^{n},\d\s^{n}),\Q\right)\t H^{i}\left(\left((\C^{*})^{n},\BI^{n}\right),\Q\right)\longrightarrow\Q(-n)\,;\]
while

\begin{equation} H^i \left( \left( (\C^*)^n ,\BI^n \right) ,\Q \right) \cong \left\{ \begin{matrix} \Q (-n) \mspace{20mu} i=n \\ 0 \mspace{30mu} \text{otherwise} \end{matrix} \right\} , \end{equation}

${}$\\
because $\left[\frac{1}{(2\pi\sqrt{-1})^{n}}\dlog z_{1}\v\cdots\v\dlog z_{n}\right]$
has weight $2n$ (due to the $n$ dlog's; see e.g. \cite{GS}). The
$\Q(-n)$ must pair with a $\Q(0)$, and so (writing $\pi_{S}:\, X\times S\to S$
as above)\[
R^{2n-i}\pi_{S}{}_{_{*}}\Q\,\cong\,\left\{ \begin{array}{c}
\Q(0)\mspace{30mu}i=n\\
0\mspace{30mu}\text{otherwise}\end{array}\right.\,\,\,\implies\,\,\, R^{n}\pi_{S}{}_{_{*}}\Q(n)=\Q(n).\]
We emphasize that this is highly peculiar; it has the effect of giving
$H^{*}(\ns)$ where one would expect (were $X$ smooth projective)
$\uh^{*}(\ns)$ in the $\Psi_{n}$-target which we now compute.

The mixed Hodge structure $H^{n-1}(\ns)$ has weights (in general)
from $(n-1)$ to $2(n-1)$; therefore $H^{n-1}(\ns,\Q(n))$ has weights
from $-(n+1)$ to $-2$ so that $W_{-1}$ is everything and $Gr_{0}^{W}=0$.
Similarly, $W_{0}H^{n}(\ns,\Q(n))=H^{n}(\ns,\Q(n))$. Applying this
reasoning to $(7.1)$ with $X=(\s^{n},\d\s^{n})$, one has \[
Gr_{\l}^{i}\uH^{2n}(X\times\ns,\Q(n))=0\text{\, \, for\, \, }i\neq n,\]
\small \[
Gr_{\l}^{n}\underline{Hg}^{n}(X\times\ns)\,\cong\, Hom_{_{\text{MHS}}}\left(\Q(0),\, H^{n}(\ns,\Q(n))\right)\,\cong\, F^{n}H^{n}(\ns,\C)\cap H^{n}(\ns,\Q(n)),\]
and\[
Gr_{\l}^{n-1}\underline{J}^{n}(X\times\ns)\,\cong\, Ext_{_{\text{MHS}}}^{1}\left(\Q(0),\, H^{n-1}(\ns,\Q(n))\right)\,\cong\,\frac{H^{n-1}(\ns,\C)}{F^{n}H^{n-1}(\ns,\C)+H^{n-1}(\ns,\Q(n))}\]
\[
\cong\, H^{n-1}(\ns,\C/\Q(n)).\]
\normalsize So the formalism predicts a map $\Psi$ from $CH^{n}(\s_{k}^{n},\d\s_{k}^{n})$
to the middle term of

\begin{equation} H^{n-1}(\ns ,\C/\Q(n)) \buildrel q \over \to H^n _{\de} (\ns ,\Q(n)) \buildrel p \over \to F^n H^n (\ns ,\C) \cap H^n (\ns ,\Q(n)) , \end{equation}

${}$\\
such that $p(\Psi(\z))=[\zeta]$ (and if this is zero $\Psi(\z)=q[AJ\zeta]$).

Next, in order to get a more interesting Leray filtration, we take
$X=(\P_{\bar{\Q}}^{1},\{0,\infty\})^{n}$ -- a product of {}``degenerate
elliptic curves''. This serves as a transitional example between
$(\s^{n},\d\s^{n})$ and products of smooth projective curves.

For fixed $i$, let $\mathfrak{{S}}_{i}^{n}$ denote the set of strictly
increasing functions $\sigma:\,\{1,\ldots,i\}\hookrightarrow\{1,\ldots,n\}$,
which we may also think of as a choice of $i$ indices \[
1\leq\sigma(1)\leq\cdots\leq\sigma(i)\leq n\,;\]
obviously $|\mathfrak{{S}}_{i}^{n}|=\frac{n!}{(n-i)!\, i!}$. To each
$\sigma\in\mathfrak{{S}}_{i}^{n}$ corresponds a projection\[
\pi_{\sigma}:\,(\C^{*})^{n}\to(\C^{*})^{i},\]
 and these induce maps on cohomology (twisted to have weight $0$)\[
H^{i}\left((\P^{1}\m\{0,\infty\})^{n},\Q(i)\right)\,\cong\,\oplus_{\sigma}\left\langle \v^{\sigma}\dlog\mathbf{{z}}\right\rangle \begin{array}{c}
_{\cong}\\
\longleftarrow\\
^{\oplus_{\sigma}\pi_{\sigma}^{*}}\end{array}\oplus_{\sigma}H^{i}\left((\P^{1}\m\{0,\infty\})^{i},\Q(i)\right)\]
 where $\v^{\sigma}\dlog\mathbf{{z}}:=\dlog z_{\sigma(1)}\v\cdots\v\dlog z_{\sigma(i)}$
and $\oplus_{\sigma}$ is short for $\oplus_{\sigma\in\mathfrak{{S}}_{i}^{n}}$. 

Let $[\v^{\sigma}\dlog\mathbf{{z}}]^{\vee}$ be a functional evaluating
(for each $\sigma'$) to $\delta_{\sigma\sigma'}$ on $[\v^{\sigma'}\dlog\mathbf{{z}}]$;
$\left\langle \v^{\sigma}\dlog\mathbf{{z}}\right\rangle ^{\vee}$
will denote the $\Q$-vector-space it generates. The last map dualizes
to

\begin{equation} H^{2n-i} \left( (\P^1 , \{0 , \infty \} )^n ,\Q(n-i) \right) \cong \oplus_{\sigma} \left< \v^{\sigma} \dlog \mathbf{z} \right> ^{\vee} \buildrel \cong \over \to \oplus_{\sigma} H^i \left( (\P^1 , \{0, \infty \} )^i , \Q \right)  \end{equation}

${}$\\
(also weight $0$) so that {[}for $X=(\P^{1},\{0,\infty\})^{n}${]}\[
R^{2n-i}\pi_{S}{}_{_{*}}\Q(n)\,\cong\,\Q(i)\t R^{2n-i}\pi_{S}{}_{_{*}}\Q(n-i)\,\cong\,\oplus_{\sigma}\Q(i)\t\left\langle \v^{\sigma}\dlog\mathbf{{z}}\right\rangle ^{\vee}\]
\[
\mspace{200mu}\cong\,\oplus_{\sigma}\Q(i).\]
 Twisted by $\Q(i)$, $(7.3)$ induces the top $\cong$ in the following
diagram, in which we are also (formally) applying functorial properties
of the $\l^{i}$ (proved in \cite{L2}):

\begin{equation} \begin{matrix} Gr^i_{\l} \uH^{2n} \left( (\P^1 , \{0,\infty \} )^n \times \ns ,\Q(n) \right) & \buildrel \cong \over \longrightarrow & \oplus_{\sigma} Gr^i_{\l} \uH^{2i} \left( ( \P^1 , \{0, \infty \} )^i \times \ns ,\Q(i) \right) \\ \uparrow \Psi_i (n) & & \uparrow \oplus_{\sigma} \Psi_i (i) \\ \l^i CH^n \left( (\P^1_k ,\{0, \infty \} )^n \right) & \buildrel \oplus (\pi_{\sigma}) {}_{_{*}} \over \longrightarrow & \oplus _{\sigma} \l^i CH^i \left( (\P^1_k , \{0, \infty \} )^i \right) .\end{matrix} \end{equation}

${}$\\
This says that $\Psi_{i}(n)$ factors through $\oplus_{\sigma}(\pi_{\sigma})_{*}$;
combining this with the expectation that $\Psi_{i}(i)$ is injective,
we have formally that\[
\l^{i+1}CH^{n}\left((\P_{k}^{1},\{0,\infty\})^{n}\right)\,=\,\ker(\Psi_{i}(n))\,=\,\ker\{\oplus_{\sigma}(\pi_{\sigma})_{*}\}\]
\[
\mspace{50mu}=\,\left\{ \left\langle \z\right\rangle \,\left|\,(\pi_{\sigma})_{*}\z\ratequiv0\,(\forall\,\sigma\in\mathfrak{{S}}_{i}^{n})\right.\right\} \]
consists of $0$-cycles whose projections to {}``$i$-faces'' $(\P^{1},\{0,\infty\})^{i}$
are rationally equivalent to zero. This obviously $does$ give a filtration,
and we expect maps

\begin{equation} \Psi_i : \, Gr^i_{\ell} CH^n \left( (\P^1_k , \{0,\infty \} )^n \right) \longrightarrow H^i_{\de} (\ns ,\Q(i)) \t \left\{ H^i \left( (\P^1 \m \{0, \infty \} )^n ,\Q(i) \right) \right\} ^{\vee} \end{equation}

${}$\\
(or equivalently, to $\oplus_{\sigma}H_{\de}^{i}(\ns,\Q(i))$) by
a calculation similar to that done above for $(\s^{n},\d\s^{n})$.

\section{\textbf{Enter the Milnor regulator}}

We now explain how maps agreeing with the predictions of $\S7$ are
obtained; we remind the reader that all Chow, cycle, and $K$-groups
are $\t\Q$ (modulo torsion).

Let $\Z\{\ss\}:=$ the free abelian group on the set $\ss$, and $\Q\{\ss\}:=\Z\{\ss\}\t\Q.$
The map\[
\t^{n}\Q\{ k^{*}\}\longrightarrow Z^{n}(\s_{k}^{n},\d\s_{k}^{n})\mspace{50mu}[\{\mathbf{a}\}:=a_{1}\t\cdots\t a_{n}\mapsto(a_{1},\ldots,a_{n})=:(\mathbf{a})]\]
 descends to an isomorphism (taking the quotient by Steinberg relations
$\t\Q$ on the left and rational equivalences $\t\Q$ on the right)\[
K_{n}^{M}(k)\begin{array}{c}
_{\cong}\\
\longrightarrow\\
^{\Xi}\end{array}CH^{n}(\s_{k}^{n},\d\s_{k}^{n})\mspace{50mu}[\{\mathbf{a}\}:=\{ a_{1},\ldots,a_{n}\}\mapsto\left\langle (a_{1},\ldots,a_{n})\right\rangle ].\]

\begin{rem*}
By normalization of Bloch's higher Chow complex, this reduces to $K_{n}^{M}(k)\cong CH^{n}(\text{Spec}k,\, n)$;
proof of the latter (the essential step is Suslin reciprocity) may
be found in \cite{T}.
\end{rem*}
Suppose $(\forall i)$ $f_{i}\mapsto a_{i}$ (i.e. $\f\mapsto\A$)
under the embedding $e:\,\bar{\Q}(S)\mapcong k\,[\subseteq\C]$ given
by evaluation at $0\in\ns,$ and let $\z=(\A)=(e(\f)).$ Then $\bar{\zeta}$
is the graph {[}closure{]} of $f_{1}\t\cdots\t f_{n}$ over $S$,
$\overline{\gamma_{\f}}\subseteq Z^{n}\left((\s_{\bar{\Q}}^{n},\d\s_{\bar{\Q}}^{n})\times S\right).$
Writing $\gamma_{\f}=\zeta$ for its restriction to $(\s_{\bar{\Q}}^{n},\d\s_{\bar{\Q}}^{n})\times\ns$,
the spread map\[
CH^{n}(\s_{k}^{n},\d\s_{k}^{n})\mapcong CH^{n}\left((\s_{\bar{\Q}}^{n},\d\s_{\bar{\Q}}^{n})\times\ns\right)\]
 sends $[\left\langle (\A)\right\rangle =]\,\left\langle (e(\f))\right\rangle \mapsto\left\langle \gamma_{\f}\right\rangle .$

\begin{lem}
\emph{\cite{Ke1}} There exists a regulator map\[
\mathcal{R}:\, K_{n}^{M}(\bar{\Q}(S))\to H_{\de}^{n}(\ns,\Q(n))\]
 such that (referring to $(7.2)$) \emph{$p(\mathcal{R}\{\f\})=[\gamma_{\f}]$}
and (if this is zero) \emph{$q(AJ(\gamma_{\f}))=\mathcal{R}\{\f\}$},
for suitable determinations of \emph{$[\gamma_{\f}]\in F^{n}H^{n}(\ns,\C)\cap H^{n}(\ns,\Q(n))$}
and \emph{$AJ(\gamma_{\f})\in H^{n-1}(\ns,\C/\Q(n)).$}
\end{lem}
If we assume this, taking $\Psi:=\mathcal{R}\circ e^{-1}\circ\Xi^{-1}$
makes good on our prediction from $\S7$, since then (for $\z=(\A)=(e(\f))$)\[
(p\circ\Psi)(\z)=p\left((\mathcal{R}\circ e^{-1})\{\A\}\right)=p(\mathcal{R}\{\f\})=[\gamma_{\f}]=[\zeta],\]
 and so on.

We briefly (and somewhat roughly) summarize the computation of $[\gamma_{\f}]$
and $AJ(\gamma_{\f})$ involved in the construction of $\mathcal{R}$
(from \cite{Ke1}). The idea is that $\gamma_{\f}$ gives a topological
$2t$-cycle essentially on $(\C^{*},\{1\})^{n}\times S_{rel}$, where
$S_{rel}:=\varinjlim_{_{D/\bar{\Q}\subseteq S}}(S,D)$. Moreover each
$(\C^{*},\{1\})$ factor has a fundamental domain $(\dd,\{1\})$ obtained
by cutting along $\R^{-}$, endowed with a $C^{\infty}$ homotopy
$\theta:\,\dd\times[0,1]\to\C^{*}$ contracting $\dd$ to $\{1\}$.
One may lift $\gamma_{\f}$ to $\widetilde{\gamma_{\f}}$ on $(\dd,\{1\})^{n}\times S_{rel}$,
then {}``successively'' contract $(\dd,\{1\})$-factors. Under $this$
{}``total homotopy'' $\widetilde{\gamma_{\f}}$ traces out a $(2t+1)$-chain
$\Gamma_{\f}$ (on $(\C^{*})^{n}\times S$) with $\d\Gamma_{\f}\equiv\gamma_{\f}-(S^{1})^{n}\times T_{\f}$,
where $T_{\f}:=\pi_{S}\left\{ \gamma_{\f}\,\cap\,(\R^{-}\times\cdots\times\R^{-})\times S\right\} $.
(We ignore the rest of the boundary, in particular that supported
on $S_{rel}\times\BI^{n}$.)

Now write $\O_{n}:=\v^{n}\dlog\mathbf{{z}}=\frac{dz_{1}}{z_{1}}\v\cdots\v\frac{dz_{n}}{z_{n}}$;
this generates cohomology of $(\C^{*},\{1\})^{n}$. We may define
the fundamental class of $\gamma_{\f}$ as a functional

\[
[\gamma_{\f}]:=\left[\int_{\gamma_{\f}}\O_{n}\v\pi_{S}^{*}(\cdot)\right]\in\{ H^{2t-n}(S_{rel},\C)\}^{\vee}\cong H^{n}(\ns,\C),\]
 which identifies with {[}integration%
\footnote{(e.g., of forms compactly supported away from some $D$ in the limit)%
} on $S$ against{]} the holomorphic current $\linebreak$ $\v^{n}\dlog\mathbf{{f}}=:\O_{\f}$
or equivalently (using $\Gamma_{\f}$ and Stokes' theorem) $(2\pi\sqrt{-1})^{n}\delta_{T_{\f}}.$
We write \[
[\gamma_{\f}]=[\O_{\f}]=[(2\pi\sqrt{-1})^{n}\delta_{T_{\f}}]\in F^{n}H^{n}(\ns,\C)\cap H^{n}(\ns,\Q(n)).\]
 If this class vanishes then $T_{\f}=\d\mu$ ($\mu$ a $(2t-n+1)$-chain$\t\Q$)
and $\Gamma_{\f}':=\Gamma_{\f}+(S^{1})^{n}\times\mu$ $\implies$
$\d\Gamma_{\f}'\equiv\gamma_{\f}.$ Accounting for ambiguity in choice
of $\mu$, define \[
AJ(\gamma_{\f}):=\left[\int_{\Gamma_{\f}'}\O_{n}\v\pi_{S}^{*}(\cdot)\right]\in\frac{\{ H^{2t-n+1}(S_{rel},\C)\}^{\vee}}{\text{im}\{ H_{2t-n+1}(S_{rel},\Q(n))\}}\cong H^{n-1}(\ns,\C/\Q(n)).\]
 It is computed by the Milnor-regulator current $R_{\f}+(2\pi\sqrt{-1})^{n}\delta_{\mu}=:R_{\f}'\in Z_{\di}(\DI_{\ns}^{n-1}),$
where \[
R_{\f}:=\sum_{j=1}^{n}(\pm2\pi\sqrt{-1})^{j}\log f_{j}\dlog f_{j-1}\v\cdots\v\dlog f_{n}\cdot\delta_{T_{f_{1}}\cap\cdots\cap T_{f_{j-1}}}\]
 (and $T_{f}:=f^{-1}(\R^{-})$, $\pm:=(-1)^{n-1}$). We prefer to
view $R_{\f}'$ as a $\C/\Q(n)$-valued functional on topological
$(n-1)$-cycles ($\t\Q$) avoiding $D$ (rather than forms compactly
supported away from $D$). This allows us to discard the membrane
term. Finally, a technical point: to completely describe $\mathcal{R}$
(and thus $\Psi$) one needs a Deligne-homology description of $H_{\de}^{n}(\ns,\Q(n))$,
or the description below in $\S9$; then (see \cite[Sect. 5.9]{KLM})
we put $\mathcal{R}\{\f\}=\left((2\pi\sqrt{-1})^{n}T_{\f},\O_{\f},R_{\f}\right).$
In terms of our formal analogy to the theory of Part $1$ and the
$\cong$'s of $\S7$, we have the following 

\begin{conclusion*}
If $\z=\left(e(f_{1}),\ldots,e(f_{n})\right)\in Z^{n}(\s_{k}^{n},\d\s_{k}^{n})$,
then \[
[\zeta]_{n}=[\O_{\f}]=[(2\pi\sqrt{-1})^{n}T_{\f}]\in Gr_{\l}^{n}\underline{Hg}^{n}\left((\s^{n},\d\s^{n})\times\ns\right),\]
 and if this vanishes\[
[AJ\zeta]_{n-1}=[R_{\f}]\in Gr_{\l}^{n-1}\underline{J}^{n}\left((\s^{n},\d\s^{n})\times\ns\right);\]
 in other words $\Psi(\z)=\mathcal{R}\{\f\}$.
\end{conclusion*}
\begin{rem}
(i) There is a natural generalization of this to higher-dimensional
cycles on the $n$-cube, to the effect that the one nontrivial higher
$AJ$ map $\Psi\,(=\Psi_{2p-n})$ is computed by the composition\[
CH^{p}(\s_{k}^{n},\d\s_{k}^{n})\begin{array}{c}
_{\cong}\\
\to\\
^{'\Xi^{-1}}\end{array}CH^{p}(\text{Spec}\, k,\, n)\begin{array}{c}
_{\cong}\\
\to\\
^{'e^{-1}}\end{array}CH^{p}(\ns,n)\mapRprime H_{\mathcal{{D}}}^{2p-n}(\eta_{S},\Q(p))\]
 where $'\mathcal{R}$ are the $AJ$ maps of \cite{KLM}. (See \cite{Ke3}
for examples.)

(ii) The isomorphism $\Xi:\, K_{n}^{M}(\C)\mapcong CH^{n}(\s^{n},\d\s^{n})$
(where we have taken the limit over $k\subset\C$) gives a concrete
demonstration of the Mumford(-Roitman) theorem, since the map $\pi_{N}:\, S^{(N)}\left((\C^{*})^{\times n}\right)\to K_{n}^{M}(\C)$
fails to surject for all $N$ (there are always elements of $K_{n}^{M}(\C)$
that cannot be presented as a sum of $\leq N$ symbols). A warning
in this connection: that $K_{n}^{M}(\C)[\t\Q]$ is an $\infty$-dimensional
$\Q$-vector-space is irrelevant (as such), since the $\infty$-dimensionality
which interests us is essentially $/\C$. After all, the Jacobian
of any curve ($g\geq1$) contains $\Q$-vector-spaces of arbitrarily
high dimension; simply choose algebraically independent generators
$\in\C^{g}$ (which are also independent of the period vectors)!
\end{rem}
We give a more terse description of the $\Psi_{i}$ for $X=(\P_{\bar{\Q}}^{1},\{0,\infty\})^{n}$,
starting with $\Psi_{n}$. First recall $\l^{i}CH^{n}\left((\P_{k}^{1},\{0,\infty\})^{n}\right):=\left\{ \left\langle \z\right\rangle \,\left|\,(\pi_{\sigma})_{*{}}\z\ratequiv0\,\,(\forall\,\sigma\in\mathfrak{{S}}_{i-1}^{n})\right.\right\} .$

\begin{lem}
$\l^{n}CH^{n}$ is spanned by $n$-box cycles of the form\[
\z=B(\A)=B(a_{1},\ldots,a_{n}):=((a_{1})-(1))\times\cdots\times((a_{n})-(1)),\]
 $a_{i}\in k^{*}$. The map $\Xi_{n}:\, K_{n}^{M}(k)\to\l^{n}CH^{n}\left((\P_{k}^{1},\{0,\infty\})^{n}\right)$
given by $\{\A\}\mapsto\left\langle B(\A)\right\rangle ,$ is an isomorphism
with inverse sending $\left\langle ({}'\A)\right\rangle \mapsto\{{}'\A\}$
(since all symbols involving a {}``$1$'' are trivial in Milnor
$K$-theory).
\end{lem}
\begin{example*}
The $2$-box cycle $B(a_{1},a_{2})=(a_{1},a_{2})-(1,a_{2})-(a_{1},1)+(1,1).$
Cycles of this form generate the Albanese kernel in $CH^{2}(\P^{1}\times\P^{1},\,\{0,\infty\}\times\P^{1}\cup\P^{1}\times\{0,\infty\}).$
\end{example*}
\begin{proof}
(of Lemma): See \cite{Ke3}; here we just point out why $\Xi_{n}$
surjects. Using the norm in Milnor $K$-theory, one can modify $\z\in\l^{n}CH^{n}\left((\P_{k}^{1},\{0,\infty\})^{n}\right)$
by a rational equivalence so that the coordinates of its points lie
in $k^{*}$ (rather than some algebraic extension). Then one uses
the fact that its face projections are rationally equivalent to zero
to get $\z\ratequiv\z+\sum_{j<n}(-1)^{n-j}\sum_{\sigma\in\mathfrak{{S}}_{j}^{n}}(\tilde{\pi}_{\sigma})_{*}\z=B(\z).$
Here $\tilde{\pi}_{\sigma}:\,(\C^{*})^{n}\to(\C^{*})^{n}$ denotes
the internal projection corresponding to $\pi_{\sigma}$; e.g., if
$\pi_{\sigma}(z_{1},\ldots,z_{5})=(z_{1},z_{4},z_{5})$ then $\tilde{\pi}_{\sigma}(z_{1},\ldots,z_{5})=(z_{1},1,1,z_{4},z_{5})$.
\end{proof}
An $n$-box cycle $\z=B(e(\f))$ spreads to an $n$-box graph $\zeta=B(\gamma_{\f})\in\linebreak Z^{n}\left((\P_{\bar{\Q}}^{1},\{0,\infty\})^{n}\times\ns\right)$;
one checks that $[\zeta]_{0}=\cdots=[\zeta]_{n-1}=0$ and $[\zeta]=[\zeta]_{n}$
is computed by $[\O_{\f}]$. Using homotopies one constructs a topological
chain $\Gamma_{\f}^{B}$ having $\Gamma_{\f}$ as a component, which
satisfies $\d\Gamma_{\f}^{B}=B(\gamma_{\f})-(S^{1})^{n}\times T_{\f}$.
(Here we do $not$ ignore boundary at $\BI^{n}\times S_{rel}$, since
this is $B(\gamma_{\f})\m\gamma_{\f}$!) As before one computes $AJ(B(\gamma_{\f}))$
(if $[\zeta]=0$) by integrating $\int_{\Gamma_{\f}^{B}}\O_{n}\v\pi_{S}^{*}(\cdot)$
(plus membrane term); but as the only nonzero contribution comes from
$\Gamma_{\f}$, $[AJ(B(\gamma_{\f}))]_{n-1}$ is just given by $[R_{\f}]$
after all. See \cite{Ke3} for details.

So for $\left\langle \z\right\rangle =\Xi_{n}(e\{\f\})\in\l^{n}CH^{n}$,
$c_{\de}(\zeta)=\mathcal{R}\{\f\}$; and we write $\Psi_{n}$ for
the composition\[
\l^{n}CH^{n}\left((\P_{k}^{1},\{0,\infty\})^{n}\right)\begin{array}{c}
_{\cong}\\
\to\\
^{\Xi_{n}^{-1}}\end{array}K_{n}^{M}(k)\begin{array}{c}
_{\cong}\\
\to\\
^{e^{-1}}\end{array}K_{n}^{M}(\bar{\Q}(S))\mapR H_{\de}^{n}(\ns,\Q(n)).\]
 To construct the remaining {}``$\Psi_{i}$'', we need the following

\begin{lem}
(i) $\l^{i}CH^{n}\left((\P_{k}^{1},\{0,\infty\})^{n}\right)$ is generated
by {}``$i$-box'' cycles of the form\[
((a_{1})-(1))\times\cdots\times((a_{i})-(1))\times(a_{i+1})\times\cdots\times(a_{n})\]
(written in any order), where the $a_{j}\in k^{*}.$

(ii) $Gr^{i}CH^{n}\left((\P_{k}^{1},\{0,\infty\})^{n}\right)\begin{array}{c}
_{\cong}\\
\to\\
^{\oplus_{\sigma}(\pi_{\sigma})_{*}}\end{array}\oplus_{\sigma\in\mathfrak{{S}}_{i}^{n}}K_{i}^{M}(k).$ \emph{{[}}The above {}``$i$-box'' would be sent to $\{ a_{1},\ldots,a_{i}\}$
in one $\oplus$-factor, $0$ in the others.\emph{{]}}
\end{lem}
\begin{proof}
keys are Suslin reciprocity (in the form proved in \cite{Ke2}), some
combinatorial considerations (see \cite{Ke3}), and a norm argument.
\end{proof}
Now consider the composition\[
Gr_{\l}^{i}CH^{n}\left((\P_{k}^{1},\{0,\infty\})^{n}\right)\begin{array}{c}
_{\oplus(\pi_{i}^{\sigma})_{*}}\\
\to\\
^{\cong}\end{array}\oplus_{\sigma}\l^{i}CH^{i}\left((\P_{k}^{1},\{0,\infty\})^{i}\right)\]
\[
\mapblah\oplus_{\sigma}H_{\de}^{i}(\ns,\Q(i));\]
one can argue (by combining a diagram like $(7.4)$ with the argument
for $\Psi_{n}$) that on cycles of the form described in Lemma $8.4(i)$,
this computes the $i^{\text{th}}$ Leray component of $c_{\de}$ of
the spread. We label this $\Psi_{i}$, with the caveat that one needs
injectivity of $\mathcal{R}$ {[}=BBC in this context{]} to get $\ker(\Psi_{i})=\l^{i+1}CH^{n}$.
Referring to $(7.5)$, we can be more precise and write

\begin{equation} \Psi_i (\z) = \sum_{\sigma \in \mathfrak{S}_i^n} \mathcal{R} \left\{ \f_{\sigma} \right\} \t \left[ \v^{\sigma} \dlog \mathbf{z} \right] ^{\vee} \in H^i_{\de} (\ns ,\Q(i)) \t \left\{ H^i \left( (\P^1 \m {0,\infty } )^n ,\Q(i) \right) \right\} ^{\vee} ,\end{equation}

${}$\\
where $\{\f_{\sigma}\}=(e^{-1}\circ\Xi_{i}^{-1}\circ(\pi_{\sigma})_{*})\left\langle \z\right\rangle $.

\begin{rem}
BBC on $X$ ($not$ to be confused with injectivity of $\mathcal{R}$!)
means that $\l^{2}CH^{n}(X_{k})$ should be zero when $k=\bar{\Q}.$
Lemma $8.4(ii)$ shows that this holds for $X=(\P_{\bar{\Q}}^{1},\{0,\infty\})^{n}$
since $K_{i}^{M}(\bar{\Q})=0$ for $i\geq2$ (see \cite{R}).
\end{rem}

\section{\textbf{Interpretation via differential characters}}

We return to the case $X=(\s_{\bar{\Q}}^{n},\d\s_{\bar{\Q}}^{n})$,
with $\z=(e(\f))\in Z^{n}(X_{k})$ and $\zeta=\gamma_{\f}\in Z^{n}(X\times\eta_{S}{}_{/\bar{\Q}}).$
Up to this point we have presented $\mathcal{R}\{\f\}=\Psi(\z)=c_{\H}(\z)\in H_{\de}^{n}(\ns,\Q(n))$
as an algebraic means of collecting $[\O_{\f}]=[\zeta]$ and $[R_{\f}]=AJ(\zeta)$
into a single object. We now give a geometric description of this
object. In what follows, $D$ refers to some codimension-$1$ subvariety
$D/\bar{\Q}\subset S$ (in the limit defining $\ns$ and $S_{rel}$);
when dealing with $\gamma_{\f}$, $R_{\f}$, etc. it is assumed that
$D\supseteq\cup|(f_{i})|$.

A (holomorphic) differential $(n-1)$-character on $S\m D$ is for
us a functional on cycles $\c\in Z_{n-1}^{top}(S\m D)$ with the following
property: when $\c=\d\mu$ on $S\m D$ (i.e., $\mu$ is compactly
supported away from $D$), the functional reduces to a membrane integral
$\int_{\mu}\O$ (for some holomorphic form $\O$ on $S\m D$). The
group $H_{\de}^{n}(\ns,\Q(n))$ can be described in terms of currents
as

\begin{equation} \frac{ \ker \left\{ \Gamma ( \DI ^{n-1} _{\ns} ) \buildrel \di \over \to \frac{Z( \DI^n_{\ns} )}{Z(F^n \DI ^n_{\ns} ) + Z \left( C^{S_{rel}}_{2t-n} (\Q(n)) \right) } \right\} }{ \di \Gamma ( \DI^{n-2}_{\ns} ) + \Gamma \left( C^{S_{rel}}_{2t-n+1} (\Q(n)) \right) } . \end{equation}

${}$\\
(Here $\DI_{S\m D}^{*}:=\DI_{S}^{*}\left/\{\text{currents\, supported\, on\, }D\}\right.$
; while $C_{*}^{(S,D)}:=$$\linebreak$ $C_{*}^{S}\left/\{\text{chains\, supported\, on\, }D\}\right.$.)
Thus a current $R$ representing a class in $H_{\de}^{n}$ satisfies
$\di[R]\equiv\Omega+(2\pi\sqrt{-1})^{n}\delta_{\KK}$ (on $S\m D$,
i.e. modulo currents supported on $D$). Such an $\int_{(\cdot)}R$
defines a $\C/\Q(n)$-$valued$ $differential$ $character$, since
$\int_{\d\mu}R\,=\,\,\int_{\mu}\Omega\,+\,(2\pi\sqrt{-1})^{n}\#(\mu\cap\KK)\,\equiv\,\int_{\mu}\O$
mod $\Q(n)$.

\begin{prop}
The current \emph{$R_{\f}$} computes the class \emph{$\mathcal{R}\{\f\}\in H_{\de}^{n}(\ns,\Q(n)).$}
\end{prop}
\begin{proof}
That $R_{\f}\in$numerator of $(9.1)$ follows from $\di[R_{\f}]=\O_{\f}-(2\pi\sqrt{-1})^{n}\delta_{T_{\f}}$
(on $\ns$). In the Deligne homology description of $H_{\de}^{n}(\ns,\Q(n))$,
essentially%
\footnote{here $\text{D}(a,b,c)=(\di[c]-b+a)$, and in the denominator $\text{D}(a,c)=(-\d a,0,\di[c]+a).$
One really needs currents with log poles here (and a reduction-to-NCD
argument); e.g. see \cite[Sect. 5.9]{KLM}.%
}\[
\frac{\ker\left\{ Z\left(C_{2t-n}^{S_{rel}}(\Q(n))\right)\oplus Z(F^{n}\DI_{\ns}^{n})\oplus\Gamma(\DI_{\ns}^{n-1})\mapD Z(\DI_{\ns}^{n})\right\} }{\text{im}\left\{ \Gamma\left(C_{2t-n+1}^{S_{rel}}(\Q(n))\right)\oplus\Gamma(\DI_{\ns}^{n-2})\mapD\text{num}\right\} },\]
 $\mathcal{R}\{\f\}$ is taken to be the class defined by $\left((2\pi i)^{n}T_{\f},\O_{\f},R_{\f}\right)$.
But this description maps isomorphically to $(9.1)$ simply by sending
$(T,\O,R)$ {[}in the numerator{]} to $R$.
\end{proof}
There is a more interesting perspective on why $R_{\f}$ gives a differential
character. Given $\c\in Z_{n-1}^{top}(S\m D)$, write $\f(\c)$ for
the $(n-1)$-cycle $\rho_{X}(\rho_{S}^{-1}(\c))$ on $(\C^{*})^{n}$.
(This $\f(\cdot)$ is actually defined on any element of $C_{*}^{top}(S\m D)$,
chains compactly supported away from $D$. One should note in analogy
to Remark $2.1$, that $\rho_{X}\circ(\rho_{S})^{-1}$ really means
$\sum_{i}m_{i}\rho_{X}^{i}\circ(\rho_{S}^{i})^{-1}$.)

\begin{claim*}
$\exists$ $\Gamma_{\c}\in C_{n}^{top}\left((\C^{*})^{n}\right)$
such that $\d\Gamma_{\c}\equiv\f(\c)$ (mod $\BI^{n}$), and $\int_{\Gamma_{\c}}\v^{n}\dlog\mathbf{{z}}\equiv\int_{\c}R_{\f}$
(mod $\Q(n)$). {[}Note that changing the choice of $\Gamma_{\c}$
(by a cycle) only alters $\int_{\Gamma_{\c}}\O_{n}$ by $\Q(n)$.{]}
\end{claim*}
In other words, $\int_{\c}R_{\f}$ can always be computed by a membrane
integral, only on $X$ (rather than $\ns$). The claim is key to understanding
how the form $\v^{n}\dlog\mathbf{{z}}=\O_{n}$ (and thus $H^{n}(X)$)
controls rational equivalence. Assuming it, if $\c=\d\mu$ then setting
$\Gamma_{\c}:=\f(\mu)$ $\implies$ $\int_{\c}R_{\f}\equiv\int_{\Gamma_{\c}}\O_{n}=\int_{\f(\mu)}\O_{n}=\int_{\mu}\O_{\f}$,
as desired. To prove the claim, we recall the {}``total homotopy''
described above as a 

\begin{lem}
$\exists$ map $\Theta:\, C_{*}^{top}\left((\C^{*})^{n}\times S\right)\to C_{*+1}^{top}\left((\C^{*})^{n}\times S\right)$
respecting support over $S$ and satisfying $\d(\Theta(\tau))\,\equiv\,\tau\,-\,\pi_{S}\left\{ \tau\,\cap\,\,(\R^{-})^{n}\times S\right\} $
modulo chains supported on $S\times\BI^{n}$.
\end{lem}
\begin{proof}
(of Claim): We abuse notation, writing $\rho_{S}^{-1}$ for $\iota\circ\rho_{S}^{-1}$.
Now $\rho_{S}^{-1}(\c)=\gamma_{\f}\cap\left\{ (\C^{*})^{n}\times\c\right\} $
is $(n-1)$-dimensional, while $(\R^{-})^{n}\times S$ has codimension
$n$. So $\d\left\{ \Theta(\rho_{S}^{-1{}}(\c))\right\} \equiv\rho_{S}^{-1}(\c)$,
and $\Theta(\rho_{S}^{-1}(\c))$ is (compactly) supported on $(\C^{*})^{n}\times S\m D$.
Set $\Gamma_{\c}:=\pi_{X}\left\{ \Theta(\rho_{S}^{-1}(\c))\right\} $;
clearly $\d\Gamma_{\c}\equiv\pi_{X}\left\{ \rho_{S}^{-1}(\c)\right\} =\f(\c)$
on $(\C^{*},\{1\})^{n}$. That $\int_{\gamma_{\c}}\O_{n}\equiv\int_{\c}R_{\f}$
follows directly from the computation that produces $R_{\f}$ in \cite{Ke1}.
\end{proof}
For $n=2$, we give pictures of $\f(\c)$ {[}$=$boldface loops{]}
together with the {}``canonical'' choice of $\Gamma_{\c}$ (from
the above proof), for two different $\c\in Z_{1}^{top}(S\m D).$ The
$\c$'s in question are just loops around components of $|(f_{1})|$
and $|(f_{2})|$ $\subseteq D$.\\

\vspace{0.3cm}
\begin{center}\includegraphics{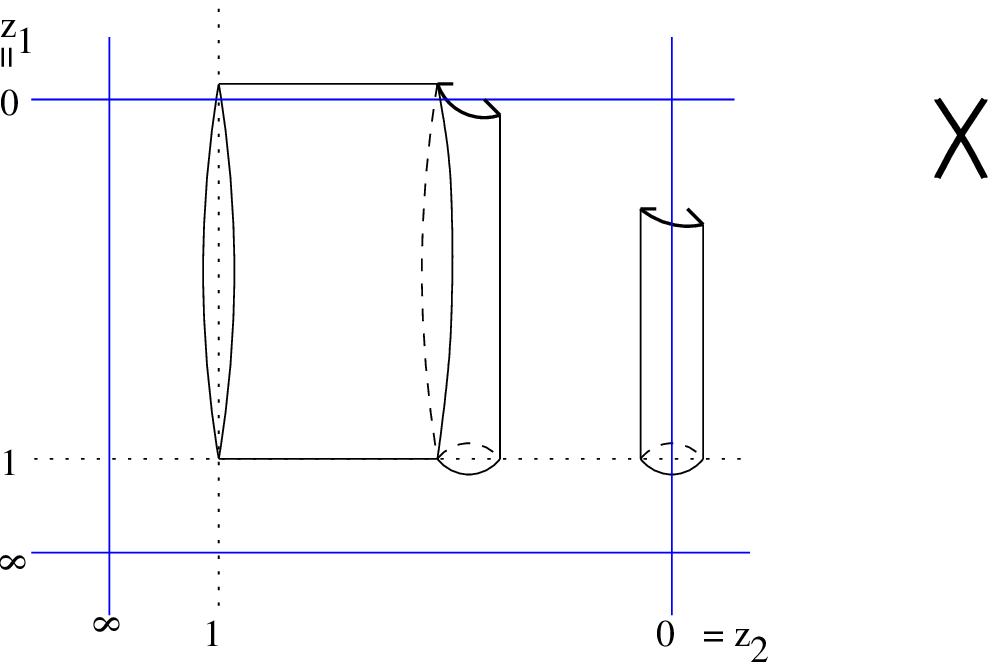}\end{center}
\vspace{0.3cm}

Now take $X=(\P^{1},\{0,\infty\})^{n}$, $\left\langle \z\right\rangle \in\l^{i}CH^{n}(X_{k})$.
In keeping with $(8.1)$, we can see\[
\Psi_{i}(\z)=\sum_{\sigma\in\mathfrak{{S}}_{i}^{n}}\,\int_{(\cdot)}R_{\f_{\sigma}}\,\t\,[\v^{\sigma}\dlog\mathbf{{z}}]^{\vee}\]
 as a formal sum of $\C/\Q(i)$-valued differential $(i-1)$-characters.
Taking $\z$ to be of the form described in Lemma $8.4(i)$, $\v^{\sigma}\dlog\mathbf{{f}}=0$
$(\forall\,\sigma\in\mathfrak{{S}}_{i-1}^{n})$ owing to cancellations;
equivalently $0=${}``$[\zeta]_{i-1}$''$:=\sum_{\sigma\in\mathfrak{{S}}_{i-1}^{n}}[\v^{\sigma}\dlog\mathbf{{f}}]\t[\v^{\sigma}\dlog\mathbf{{z}}]^{\vee}$.
So $[\zeta]_{i-1}{}^{^{*}}:=(\rho_{S})_{*}\circ(\rho_{X})^{*}:\,[F^{i-1}]H^{i-1}\left((\C^{*})^{n}\right)\to H^{i-1}(\ns)$
is zero; hence $[\zeta]_{i-1}{}_{_{*}}:=(\rho_{X}\circ(\rho_{S})^{-1})_{*}:\, H_{i-1}(\ns,\Q)\to H_{i-1}((\C^{*})^{n},\Q)$
is zero also: given $\c\in Z_{i-1}^{top}(S\m D)$, there exists $\Gamma_{\c}\in C_{i}^{top}((\C^{*})^{n})$
bounding on $\mathbf{{f}}(\c)$. A straightforward argument then shows
$\int_{\c}R_{\f_{\sigma}}=\int_{\Theta(\pi_{\sigma}(\rho^{-1}(\c)))}\O_{i}\equiv\int_{\pi_{\sigma}(\Gamma_{\c})}\O_{i}=\int_{\Gamma_{\c}}\v^{\sigma}\dlog\mathbf{{z}}\,\,\,\,(\text{mod\, }\Q(i)).$
Consequently\[
\Psi_{i}(\z)(\c)\,=\,\sum_{\sigma\in\mathfrak{{S}}_{i}^{n}}\,\int_{\Gamma_{\c}}\v^{\sigma}\dlog\mathbf{{z}}\,\t\,[\v^{\sigma}\dlog\mathbf{{z}}]^{\vee}\mspace{100mu}\]
 $\vspace{-3mm}$\[
\mspace{80mu}=\,\int_{\Gamma_{\c}}(\cdot)\,\,\in\,\frac{\left\{ [F^{i}]H^{i}((\C^{*})^{n},\C)\right\} ^{\vee}}{\text{im}\left\{ H_{i}((\C^{*})^{n},\Q(i))\right\} }\,\cong\,\oplus_{\sigma}\C/\Q(i)\,,\]
and $\Psi_{i}(\z)$ may be viewed as $one$ differential character
taking values in a vector space modulo a $\Q$-lattice. This is the
point of view we will generalize to $X$ smooth projective.

\part{Concrete description of the invariants}

\section{\textbf{Reduction of the target spaces}}

Once more take $X/\bar{\Q}$ smooth projective ($n=\dim X$), $\z\in Z^{n}(X_{k})$
such that $\left\langle \z\right\rangle \in\l^{i}CH^{n}(X_{k})$,
and spread $\z$ out to get $\zeta\in Z^{n}(X\times\eta_{S}{}_{/\bar{\Q}})$.
We want to compute $\Psi_{i}(\z)=Gr_{\l}^{i}c_{\H}(\zeta)$; but in
general the target space $Gr_{\l}^{i}\uH^{2n}(X\times\ns,\Q(n))$
can be complicated, and $\Psi_{i}$ difficult to describe geometrically.

That is in stark contrast to the relative cases considered in Part
$2$, for instance $X'=(\P^{1},\{0,\infty\})^{n}$ {[}with dual $\hat{X}'=(\C^{*})^{n}${]}.
Since $F^{n-i+1}H^{2n-i}(X',\C)=0$, rational equivalence was {}``controlled
by'' the rather simple object $\frac{[F^{n-i}]}{F^{n-i+1}}H^{2n-i}(X')$
{[}or dually by $[F^{i}]H^{i}(\hat{X}')${]}.

One naturally wonders whether neglecting $F^{n-i+1}H^{2n-i}(X,\C)[\neq0]$
in $Gr_{\l}^{i}\uH^{2n}$ markedly simplifies its presentation. The
answer is yes, and the resulting invariants $\chi_{i}(\z)$ (quotients
of $\Psi_{i}(\z)$) turn out to be differential characters. These
are computed by membrane integrals on $X$ (different from Griffiths's
prescription; see Remark $4.1$). The caveat is that they only capture
rational inequivalences {}``regulated'' by $holomorphic$ forms
on $X$.

The following notation is good for the rest of the paper. Set\[
\widetilde{\Lambda}:=H^{2n-i}(X,\Q)\,,\,\,\,\widetilde{V}:=\widetilde{\Lambda}\t\C=H^{2n-i}(X,\C)\,,\]
and note that $\widetilde{V}=F^{n-i}\widetilde{V}$. Also let\[
V:=\widetilde{V}\left/F^{n-i+1}\widetilde{V}\right.=\left\{ F^{i}H^{i}(X,\C)\right\} ^{\vee},\,\,\,\,\Lambda:=\text{im}\left\{ \widetilde{\Lambda}\mapj V\right\} \cong\widetilde{\Lambda}\left/\ker(\j)\right.,\]
 and remark that $\dim_{\Q}\Lambda\geq2\dim_{\C}V$ (usually $>$
holds). We now explain how to take the appropriate quotients of $Gr_{\l}^{i}\underline{Hg}^{n}(X\times\ns)$
and $Gr_{\l}^{i-1}\underline{J}^{n}(X\times\ns)$, then finally $Gr_{\l}^{i}\uH^{2n}$.

\subsection*{Reduction of $Gr_{\l}^{i}\underline{Hg}^{n}(X\times\ns)$}

Set \[
H_{\Q}^{2n}:=H^{i}(\ns,\Q)\t H^{2n-i}(X,\Q),\,\,\,\uh_{\Q}^{2n}:=W_{2n}H_{\Q}^{2n}=\uh^{i}(\ns,\Q)\t H^{2n-i}(X,\Q).\]
 Referring to $(4.5)$ and $(4.1)$ (omitting the twist),

\begin{equation} Gr^i_{\l} \underline{Hg}^n(X\times \ns) = Hom_{_{\text{MHS}}} \left( \Q(0) , H^{2n}_{\Q} \t \Q(n) \right) \cong F^n \uh^{2n}_{\C} \cap \uh ^{2n}_{\Q} . \end{equation}

${}$\\
It is worth noting that this maps $injectively$ to\[
\uh_{\Q}^{2n}\cong Hom_{\Q}\left(\uh_{i}(\ns,\Q),\, H_{i}(X,\Q)\right),\]
 where one thinks of $\uh_{i}(\ns,\Q)\cong\text{coim}\left\{ H_{i}(\ns,\Q)\to H_{i}(S,\Q)\right\} $
as topological cycles on $S$, which can be moved (in their class
on $S$) to avoid an arbitrary {[}not necc. irreducible{]} divisor
$D/\bar{\Q}\subset S$.

Now observe that $F^{n}\uh_{\C}^{2n}=\sum_{k=0}^{i}F^{k}\uh^{i}(\ns,\C)\t F^{n-k}H^{2n-i}(X,\C)$;
writing\[
F_{X}^{n-i+1}H_{\C}^{2n}:=\uh^{i}(\ns,\C)\t F^{n-i+1}H^{2n-i}(X,\C)\,,\,\,\,\, W_{\C}:=(F_{X}^{n-i+1}\uh_{\C}^{2n})\cap(F^{n}\uh_{\C}^{2n}),\]
 we have\[
\frac{F^{n}\uh_{\C}^{2n}}{W_{\C}}=F^{i}\uh^{i}(\ns,\C)\t\frac{[F^{n-i}]}{F^{n-i+1}}H^{2n-i}(X,\C)=:F^{i}\uh^{i}(\ns,V).\]
 (Note that $F^{i}\uh^{i}(\ns,\C)\cong F^{i}H^{i}(S,\C).$) The reduction
is then achieved via the composition

\begin{diagram}
\ker \left\{ F^n \uh^{2n}_{\C} \oplus \uh^{2n}_{\Q} \to \uh_{\C}^{2n} \right\} & \rEq : & F^n \uh^{2n}_{\C} \cap \uh ^{2n}_{\Q} \\
\dOnto>{(*)} \\
\ker \left\{ \frac{F^n \uh^{2n}_{\C} }{W_{\C}} \oplus \frac{\uh^{2n}_{\Q}}{\uh^{2n}_{\Q}\cap W_{\C}} \to \frac{\uh_{\C}^{2n}}{W_{\C}}\right\} & & \dTo>{'\phi_i} \\
\dInto \\
\begin{array}[b]{c} \ker \left\{ \frac{F^n \uh^{2n}_{\C}}{W_{\C}} \oplus \frac{\uh_{\Q}^{2n}}{\uh_{\Q}^{2n} \cap F^{n-i+1}_X \uh^{2n}_{\C}} \to \frac{\uh^{2n}_{\C}}{F^{n-i+1}_X \uh^{2n}_{\C}} \right\} \\ \downarrow \cong \\ \ker \left\{ F^i \uh^i (\ns ,V) \oplus \uh^i (\ns ,\Lambda )\to \uh^i (\ns ,V) \right\} \end{array} & \rEq : & F^i \uh^i (\ns ,V) \cap \uh^i(\ns ,\Lambda) , \\
\end{diagram}
which may also be expressed as a map%
\footnote{here $\overline{H}_{i}(X,\Q):=$Poincar\'e duals of $\Lambda$. To
get a (sort of) feeling for what this looks like: in a situation where
the $unamended$ GHC (in general false, see \cite{L1}) holds, this
is $\varinjlim_{Y\subseteq X\text{(codim.\, }1\text{)}}\text{coker}\left\{ H_{i}(Y,\Q)\to H_{i}(X,\Q)\right\} =\text{coim}\left\{ H_{i}(X,\Q)\to H_{i}(X_{rel},\Q)\right\} $.%
}\small\[
'\phi_{i}:\, Gr_{\l}^{i}\underline{Hg}^{n}(X\times\ns)\to Hom_{\C}\left\{ F^{i}H^{i}(X,\C),F^{i}\uh^{i}(\ns,\C)\right\} \cap Hom_{\Q}\left\{ \uh_{i}(\ns,\Q),\overline{H}_{i}(X,\Q)\right\} .\]
 \normalsize This time injectivity of the map (from the r.h.s.) to
the $Hom_{\C}$ piece is what shall be useful.

\begin{rem}
(i) Write $F_{h}^{1}\uh^{i}(\ns,\Q)$ for the largest sub-Hodge-structure
of $\linebreak$ $\uh^{i}(\ns,\Q)\cap F^{1}\uh^{i}(\ns,\C)$; this
is zero if {[}Grothendieck's{]} GHC holds. Since \[
\ker(*)\,\,\subseteq\,\,\left\{ \sum_{k=1}^{i-1}F^{k}\uh^{i}(\ns,\C)\t F^{n-k}H^{2n-i}(X,\C)\right\} \cap\uh_{\Q}^{2n}\]
 \[
\subseteq\, F_{h}^{1}\uh^{i}(\ns,\Q)\t H^{2n-i}(X,\Q),\]
we expect $(*)$ is an isomorphism; but we won't assume this.

(ii) By contrast (provided $V\neq0$), $\left\{ F^{i}\uh^{i}(\ns,\C)\t V\right\} \cap\left\{ \uh^{i}(\ns,\Q)\t\Lambda\right\} $
is $not$ contained in $\left\{ F^{i}\uh^{i}(\ns,\C)\cap\uh^{i}(\ns,\Q)\right\} \t\Lambda$
(which is zero), essentially because $\dim_{\Q}\Lambda>\dim_{\C}V$
(and so similar reasoning does not apply). So the target of $'\phi_{i}$
is certainly nonzero.
\end{rem}

\subsection*{Reduction of $Gr_{\l}^{i-1}\underline{J}^{n}(X\times\ns)$}

Now let\[
H_{\Q}^{2n-1}:=H^{i-1}(\ns,\Q)\t H^{2n-i}(X,\Q)\,,\,\,\,\uh_{\Q}^{2n-1}:=W_{2n-1}H_{\Q}^{2n-1}.\]
 From $(4.5)$, $(4.2)$ (omitting the twist), and the description
of the $Ext^{1}/Hom$ quotient (in $\S4$), \[
Gr_{\l}^{i-1}\underline{J}^{n}(X\times\ns)\,=\,\frac{Ext_{_{\text{MHS}}}^{1}\left(\Q(0),\, W_{-1}(H_{\Q}^{2n-1}\t\Q(n))\right)}{\text{im}\left\{ \text{Hom}_{_{\text{MHS}}}\left(\Q(0),\, Gr_{0}^{W}(H_{\Q}^{2n-1}\t\Q(n))\right)\right\} }\]
\[
=\,\frac{Ext_{_{\text{MHS}}}^{1}\left(\Q(0),\,\uh_{\Q}^{2n-1}\t\Q(n)\right)}{\text{im}\left\{ Hom_{_{\text{MHS}}}\left(\Q(0),\,(Gr_{2n}^{W}H_{\Q}^{2n-1})\t\Q(n)\right)\right\} }\]

\begin{equation} \cong \, \frac{\uh^{2n-1}_{\C}}{\uh^{2n-1}_{\Q} + F^n \uh^{2n-1}_{\C} + \left\{ W_{2n}H^{2n-1}_{\Q} + F^n W_{2n} H^{2n-1}_{\C} \right\}\cap \uh^{2n-1}_{\C}} \end{equation}

${}$\\
(since $W_{2n}\uh_{\Q}^{2n-1}=\uh_{\Q}^{2n-1}$).

Next one notices that\[
F^{n}\uh_{\C}^{2n-1}\,=\,\sum_{k=0}^{i-1}F^{k}\uh^{i-1}(\ns,\C)\t F^{n-k}H^{2n-i}(X,\C)\]
\[
\subseteq\,\uh^{i-1}(\ns,\C)\t F^{n-i+1}H^{2n-i}(X,\C)\,=:\, F_{X}^{n-i+1}\uh_{\C}^{2n-1}.\]
 Similarly, \[
F^{n}W_{2n}H_{\C}^{2n-1}\,=\,\sum_{k=0}^{i-1}F^{k}W_{i}H^{i-1}(\eta_{S},\C)\t_{\C}F^{n-k}H^{2n-i}(X,\C)\]
\[
\mspace{70mu}\subseteq\, W_{i}H^{i-1}(\eta_{S},\C)\t_{\C}F^{n-i+1}H^{2n-i}(X,\C)\]
\[
\mspace{70mu}=\, W_{i}H^{i-1}(\eta_{S},\Q)\t_{\Q}F^{n-i+1}H^{2n-i}(X,\C),\]
hence \[
\{ W_{2n}H_{\Q}^{2n-1}+F^{n}W_{2n}H_{\C}^{2n-1}\}\cap\underline{H}_{\C}^{2n-1}\,\subseteq\]
\[
W_{i}H^{i-1}(\eta_{S},\Q)\t_{\Q}\left\{ H^{2n-i}(X,\Q)+F^{n-i+1}H^{2n-i}(X,\C)\right\} \mspace{200mu}\]
\[
\mspace{300mu}\cap\,\, W_{i-1}H^{i-1}(\eta_{S},\Q)\t_{\Q}H^{2n-i}(X,\C)\]
\[
=\, W_{i-1}H^{i-1}(\eta_{S},\Q)\t_{\Q}\left\{ H^{2n-i}(X,\Q)+F^{n-i+1}H^{2n-i}(X,\C)\right\} \]
\[
\subseteq\,\underline{H}_{\Q}^{2n-1}+F_{X}^{n-i+1}\underline{H}_{\C}^{2n-1}.\]

We write our conclusion as a projection\[
Gr_{\l}^{i-1}\underline{J}^{n}(X\times\ns)\mapphipp Gr_{\l}^{i-1}\underline{J}^{n}(X\times\ns)\left/F_{X}^{n-i+1}\uh_{\C}^{2n-1}\right.=\frac{\uh_{\C}^{2n-1}}{\uh_{\Q}^{2n-1}+F_{X}^{n-i+1}\uh_{\C}^{2n-1}}\]
\[
=\,\frac{\uh^{i-1}(\ns,\Q)\t H^{2n-i}(X,\C)}{\uh^{i-1}(\ns,\Q)\t H^{2n-i}(X,\Q)\,+\,\uh^{i-1}(\ns,\Q)\t F^{n-i+1}H^{2n-i}(X,\C)}\]
\[
\cong\, Hom\left(\uh_{i-1}(\ns,\Q),\,\frac{H^{2n-i}(X,\C)}{H^{2n-i}(X,\Q)+F^{n-i+1}H^{2n-i}(X,\C)}\right)\]
\[
=\, Hom\left(\uh_{i-1}(\ns,\Q),\,\{ F^{i}H^{i}(X,\C)\}^{\vee}\left/\text{im}\{ H_{i}(X,\Q)\}\right.\right)\]
\[
=:\,\uh^{i-1}(\ns,\, V\left/\Lambda\right.).\]

\subsection*{Reduction of $Gr_{\l}^{i}\uH^{2n}(X\times\ns,\Q(n))$}

Define a complex of sheaves in the analytic topology (on $S_{\C}$)\[
\Lambda\{ i\}_{\de}:=\,\Lambda\hookrightarrow\o_{S}\t V\to\O_{S}^{1}\t V\to\cdots\to\O_{S}^{i-1}\t V\to0\]
 (where $\Lambda=$ constant sheaf placed in degree $0$), and set\[
H_{\de}^{*}(S,\Lambda\{ i\}):=\HH^{*}(S,\Lambda\{ i\}_{\de}).\]
 One may define similar groups for $\ns$ via a Deligne-Beilinson
approach (as in \cite{EV}); and we write \[
\uh_{\de}^{i}(\ns,\Lambda\{ i\}):=\text{im}\left\{ H_{\de}^{i}(S,\Lambda\{ i\})\to H_{\de}^{i}(\ns,\Lambda\{ i\})\right\} .\]
 By a standard argument one has a diagram (with exact rows)

\begin{diagram}
0 & \rTo & H^{i-1}(S, V/\Lambda ) & \rTo & H^i_{\de} (S,\Lambda \{ i \}) & \rTo & \ker \left\{ \left( \begin{matrix} F^iH^i(S,V) \\ \oplus \, H^i(S,\Lambda) \end{matrix} \right) \to H^i (S,V) \right\} & \rTo & 0 \\
& & \dTo & & \dTo & & \dTo>{\alpha} \\
0 & \rTo & H^{i-1}(\ns ,V/\Lambda ) & \rTo & H^i_{\de} (\ns ,\Lambda \{ i \} ) & \rTo & \ker \left\{ \left( \begin{matrix} F^iH^i(\ns ,V) \\ \oplus \, H^i (\ns , \Lambda ) \end{matrix} \right) \to H^i (\ns ,V) \right\} & \rTo & 0 \, \, . \\
\end{diagram}
Here $V$ is treated simply as an abelian group (and ignored by the
Hodge filtration). Since $\ker\{ H^{i}(S,\Lambda)\to H^{i}(S,V)\}$
and $\ker\{ F^{i}H^{i}(S,V)\to F^{i}H^{i}(\ns,V)\}$ are both zero,
$\alpha$ is injective. It follows immediately that the top row of
the following (where we have not yet defined $\phi_{i}$) is exact:

\small \begin{equation} \begin{matrix} 0 \to \uh^{i-1} (\ns ,V/\Lambda ) & \rTo^q & \uh^i_{\de} (\ns ,\Lambda \{i \} ) & \rTo^p & F^i \uh ^i (\ns ,V) \cap \uh^i (\ns ,\Lambda ) \to 0 \\ '' \phi_{i-1} \uparrow & & \phi_i \uparrow & & '\phi_i \uparrow \\ 0 \to Gr^{i-1}_{\l} \underline{J}^n (X\times \ns ) & \rTo^q & Gr^i_{\l} \uH^{2n} (X\times \ns ,\Q(n)) & \rTo^p & Gr^i_{\l} \underline{Hg}^n (X\times \ns ) \to 0 \, \, . \end{matrix} \end{equation} \normalsize

${}$\\
One notices a formal similarity between the top row, and the situation
{[}$X=(\P^{1},\{0,\infty\})^{n}$, $V=\oplus_{\sigma}\C$, $\Lambda=\oplus_{\sigma}(2\pi i)^{i}\Q${]}
in Part $2$; but the restriction on the weights (e.g., $\uh^{i-1}=W_{i-1}H^{i-1}$)
constitutes a substantial difference.

The bottom row of $(10.3)$ admits a Deligne-homology-type presentation,
namely

\begin{equation} \begin{matrix} \underline{\ker \left\{ \frac{Z(C^S_{2t-i}(\WL ))}{\di \Gamma (C^S_{2t-i+1}(\WL ))} \oplus \frac{Z(F^n\DI^i_S(\WV ))}{\di \Gamma (F^n \DI^{i-1}_S (\WV ))} \to \frac{Z(\DI^i_S(\WV))}{\di \Gamma ( \DI^{i-1}_S(\WV)) } \right\} } \\ {N^1 (\text{num}) } \\ \uparrow \, p \\ {\ker \left\{ Z(C^S_{2t-i} (\WL )) \oplus Z(F^n \DI^i_S (\WV )) \oplus \Gamma (\DI^{i-1}_S (\WV )) \buildrel \text{D} \over \to Z(\DI^i_S (\WV )) \right\} } \\ \overline{\widetilde{N}^1(\text{num}) + \text{im} \left\{ \Gamma (C^S_{2t-i+1} (\WL )) \oplus \Gamma (F^n \DI^{i-1}_S (\WV)) \oplus \Gamma (\DI^{i-2}_S (\WV)) \to \text{num} \right\}} \\ \uparrow \, q \\ {Z(\DI^{i-1}_S(\widetilde{V}))} \\ \overline{ \di \Gamma (\DI^{i-2}_S (\widetilde{V})) + \left\{ N^1 \Gamma (\DI^{i-1}_S (\WV )) + W^S_i \Gamma (C_{2t-i+1} (\WL )) + W^S_i \Gamma (F^n \DI^{i-1}_S (\WV ))\right\} \cap \text{num}} \end{matrix} \end{equation}

${}$\\
(see notation%
\footnote{\textbf{Notation:} here $\DI_{S}^{*}(\WV)$ denotes $\WV$-valued
currents (and so forth), $Z(\cdots)$ is closed sections, and {}``num''
means numerator. When working with $\WV$-valued currents, the Hodge
filtration $is$ conscious of $\tilde{V}$, e.g. $F^{n}\DI_{S}^{i}(\WV)=\sum_{k}(F_{(S)}^{k}\DI_{S}^{i})\otimes(F_{(X)}^{n-k}\WV)$;
whereas in $(10.5)$ $V$ is treated (as above) as a $\C$-vector-space,
$F_{[S]}^{i}\DI_{S}^{i}(V)=(F^{i}\DI_{S}^{i})\otimes V$. 

The $N^{1}$'s work as follows: let $K_{S}$, $M_{S}$ be groups of
global sections of currents or chains on $S$ (including direct sums
and kernels). Then $N^{1}(K_{S})$ denotes sections supported on any
codimension $1$ $D_{/\bar{\Q}}\subseteq S$, so that in particular
$N^{1}\{\ker(K_{S}\to M_{S})\}=\ker\{ N^{1}(K_{S})\to M_{S}\}$. In
contrast, define $\widetilde{N}^{1}\{\ker(K_{S}\to M_{S})\}:=\varinjlim_{D/\bar{\Q}}\text{{im}}\left(\ker\{ K_{\widetilde{D}}\to M_{\widetilde{D}}\}\mapcc K_{S}\right),$
where $K_{\widetilde{D}}$ and $M_{\widetilde{D}}$ are corresponding
objects on a normalization of $D$. (Strictly speaking, since an arbitrary
$D$ has non-normal crossings, we must actually take $\widetilde{D}$
to be the normalization of the \emph{preimage} of $D$ under an embedded
resolution $b:\,\hat{S}\to S$ as in the proof of Prop. $3.2$; $\mathcal{{N}}$
is then the composition $\widetilde{D}\to\hat{S}\to S$, which still
has image in codimension $1$.) 

Finally, $W_{i}^{S}\Gamma(\cdots)$ means sections closed on $\eta_{S}$
(i.e., mod $N^{1}$) such that the image in $H^{*}(\eta_{S},\C\text{{\, or\,}}\Q)$
lies in $W_{i}$.%
} below). We explain how to derive the middle term below.

On the other hand, a (less horrible) differential-character-type presentation
of $(10.3)$'s top row will serve to motivate and guide the discussions
in $\S11$ and $\S14$:

\begin{equation} \begin{matrix} \ker \left\{ \frac{Z(C^S_{2t-i}(\Lambda))}{\di \Gamma (C^S_{2t-i+1}(\Lambda )) + N^1 (\text{num})} \oplus Z(F^i _S \DI^i_S (V)) \to \frac{Z(\DI^i_S (V))}{\di \Gamma (\DI^{i-1}_S (V)) + N^1 Z ( C^S_{2t-i} (\Lambda))} \right\} \\ \uparrow \, p \\ {\ker \left\{ \Gamma (\DI^{i-1}_S (V)) \buildrel \di \over \to \frac{Z(\DI^i_S(V))}{Z(C^S_{2t-i}(\Lambda))+Z(F^i_S D^i_S (V))} \right\} } \\ \overline{ \di \Gamma (\DI^{i-2}_S (V)) + N^1 \Gamma (\DI^{i-1}_S (V))\cap \text{num} + \Gamma (C^S_{2t-i+1} (\Lambda )) } \\ \uparrow \, q \\ Z(\DI^{i-1}_S (V)) \\ \overline{ \di \Gamma (\DI^{i-2}_S (V)) + N^1 Z(\DI^{i-1}_S (V)) + Z(C^S _{2t-i+1} (\Lambda) ) } \end{matrix} \end{equation}

${}$\\
where $p$ is obtained by taking $\di$. To obtain $\phi_{i}$ in
$(10.3)$, simply send $(\overline{\c,\O,R})$ in the middle term
of $(10.4)$, to $(\overline{R})$ in the middle term above.

We sketch an argument for the middle terms of $(10.4)$ and $(10.5)$
above; the outer terms (and proving exactness) are more straightforward
and left to the reader. First recall that $H_{\de}^{2n}(X\times S,\Q(n))$
may be computed as the $2n^{\text{{th}}}$ cohomology of the {[}Deligne
homology($=\de H$){]} complex \[
C_{\de}^{\bullet}(X\times S,\Q(n)):=\left\{ \Gamma(C_{2(n+t)-\bullet}^{X\times S}\otimes\Q(n))\oplus\Gamma(F^{n}\DI_{X\times S}^{\bullet})\oplus\Gamma(\DI_{X\times S}^{\bullet-1})\right\} \]
of abelian groups, with differential $\text{{D}}(a,b,c)=\left(-\d a,-\di[b],\di[c]-b+\delta_{a}\right)$
and where $C$ (resp. $\DI$) denotes $C^{\infty}$ chains (resp.
currents). Multiplying $a,\, b,\, c$ by $\frac{1}{(2\pi\sqrt{-1})^{n}}$,
we may omit the $\Q(n)$ twist in the first $\oplus$-factor (for
computational purposes).

Denote by $C_{\de,\times}^{\bullet}$ the subcomplex consisting of
products of chains (resp. currents) on $S$ and $X$ (for example,
$\Gamma(\DI_{X\times S}^{\bullet-1})$ is replaced by $\sum_{j+k=\bullet-1}\Gamma(\DI_{S}^{j})\otimes\Gamma(\DI_{X}^{k})$).
The inclusion $C_{\de,\times}^{\bullet}\hookrightarrow C_{\de}^{\bullet}$
is a quasi-isomorphism, and $C_{\de,\times}^{\bullet}$ admits an
obvious filtration $\l_{S}^{\bullet}$ on each term, e.g.\[
\l_{S}^{i}C_{\de,\times}^{2n}\,=\,\left\{ \left(\sum_{j\geq i}\Gamma(C_{2t-j}^{S})\otimes\Gamma(C_{j}^{X})\right)\oplus\left(\sum_{j\geq i}\sum_{\ell}\Gamma(F^{\ell}\DI_{S}^{j})\otimes\Gamma(F^{n-\ell}\DI_{X}^{2n-j})\right)\right.\]
\[
\mspace{100mu}\left.\oplus\left(\sum_{j\geq i}\Gamma(\DI_{S}^{j-1})\otimes\Gamma(\DI_{X}^{2n-j})\right)\right\} .\]
If we set 

\begin{equation} \l^i_S H^{2n}_{\de} \left( X\times S, \Q(n) \right) \, := \, \frac{\left\{ \ker (\text{D}) \subseteq \l^i_S C^{2n}_{\de, \times}  \right\}}{\text{im}\left\{ \l^{i-1}_S C^{2n-1}_{\de, \times} \to C^{2n}_{\de , \times} \right\} \cap \l^i_S C^{2n}_{\de, \times} } , \end{equation}${}$\\
then $\Phi(\l_{S}^{i}H_{\de}^{2n})=\l^{i}\underline{H}_{\H}^{2n}(X\times\eta_{S},\Q(n))$.
(One has {}``$=$'' and not just {}``$\subseteq$'', because the
maps of graded pieces are surjective.) In words, the Leray filtration
is represented by $\de H$-triples of the form $(\l^{i},\l^{i},\l^{i-1})$
in degree $2n$ modulo coboundaries coming from $(\l^{i-1},\l^{i-1},\l^{i-2})$
in degree $2n-1$. Straightforward (but tedious) arguments with $\de H$-triples
then show that the middle term of $(10.4)$ $without$ the $\widetilde{N}^{1}$,
is equal to $Gr_{\l_{S}}^{i}H_{\de}^{2n}(X\times S,\Q(n))$.

Now arguing for the normal-crossings case (replacing $\eta_{S}$ by
$S\m D$, $D/\bar{\Q}$ a NCD), one has an exact sequence (e.g., see
\cite{mS})\[
H_{\de,X\times D}^{2n}(X\times S,\Q(n))\mapdd H_{\de}^{2n}(X\times S,\Q(n))\mapee\underline{H}_{\H}^{2n}(X\times(S\m D),\Q(n))\to0\]
(where $\underline{H}_{\H}$ is $\text{{im}}\{ H_{\de}\to H_{\H}\}$).
It follows from work of Deligne%
\footnote{Use the exact sequence $W_{i}H^{i-1}(S\m D)\mapff W_{i}H_{D}^{i}(S)\to H^{i}(S)\to W_{i}H^{i}(S\m D)$
together with Deligne's description (in \cite[sec. 3.2]{De1}) of
$W_{i}H^{i-1}(S\m D)$ and the fact (\cite[Cor. 8.2.8]{De2}) that
$\text{{im}}\{ H_{D}^{i}(S)\to H^{i}(S)\}=\text{{im}}\{ H^{i-2}(\widetilde{D})\to H^{i}(S)\}$,
to get $H^{i-2}(\widetilde{D})\twoheadrightarrow W_{i}H_{D}^{i}(S)$. %
} that $H^{2n-2}(X\times\widetilde{D})\twoheadrightarrow W_{2n}H_{X\times D}^{2n}(X\times S)$,
hence that one may replace $\alpha$ above by $H_{\de}^{2n-2}(X\times\widetilde{D},\Q(n-1))\mapgg H_{\de}^{2n}(X\times S,\Q(n)).$
Since $\tilde{\alpha}$ and $\Phi$ are strictly compatible with Leray,%
\footnote{i.e., for the obvious $\l_{\widetilde{D}}^{\bullet}$ on $H_{\de}^{2n-2}(X\times\widetilde{D},\Q(n-1))$
one can show $\tilde{\alpha}(\l_{\widetilde{D}}^{i-2})=\l_{S}^{i}\cap\text{{im}}(\tilde{\alpha})$
(e.g., by arguing explicitly with $\de H$-triples). We have already
noted this property for $\Phi$.%
} one concludes that $Gr_{\l}^{i}\{\underline{H}_{\H}^{2n}\}\cong\text{{coker}}\{ Gr_{\l}^{i}(\tilde{\alpha})\}.$
To transfer this result to the non-NCD case, assume $D$($=$NCD)
is the preimage of $D_{0}$($=$non-NCD) under an embedded resolution
$b:\, S\to S_{0}$, so that $(S\m D)=(S_{0}\m D_{0})$. The surjective
map $H_{\de}^{2n}(X\times S,\Q(n))\mapbde H_{\de}^{2n}(X\times S_{0},\Q(n))$
is then (since $Gr_{\l}^{i}b_{*}^{\de}$ are surjective) strictly
compatible with Leray and factors $\Phi$; hence $Gr_{\l}^{i}\underline{H}_{\H}^{2n}(X\times(S_{0}\m D_{0}),\Q(n))\cong\frac{Gr_{\l}^{i}H_{\de}^{2n}(X\times S_{0},\Q(n))}{(\overline{b_{*}^{\de}\circ\tilde{\alpha}})\left\{ Gr_{\l}^{i-1}H_{\de}^{2n-2}(X\times\widetilde{D},\Q(n-1))\right\} }.$
Combining this with $(10.6)$ and taking a limit gives the middle
term of $(10.4)$.

Moving on to $(10.5)$, it will suffice to show that its middle term
without the $N^{1}\Gamma(\cdots)\cap\text{{num}}$ in the denominator,
is isomorphic to $H_{\de}^{i}(S,\Lambda\{ i\})$. The latter has a
standard Deligne-homology presentation as the $i^{\text{{th}}}$ cohomology
of $C_{\de}^{\bullet}(S,\Lambda\{ i\})=\Gamma(C_{2t-\bullet}^{S}(\Lambda))\oplus\Gamma(F_{S}^{i}\DI_{S}^{\bullet}(V))\oplus\Gamma(\DI_{S}^{\bullet-1}(V)).$
Moreover, $\ker\{\text{{D}}:\, C_{\de}^{i}\to C_{\de}^{i+1}\}\subseteq$

\begin{equation} Z\left( C^S_{2t-i}(\Lambda ) \right) \oplus Z \left( F^i_S \DI^i(V) \right) \oplus \ker \left\{ \Gamma \left( \DI^{i-1}_S (V) \right) \to \frac{Z(\DI^i(V))}{Z(C^S_{2t-i}(\Lambda))+Z(F^i \DI^i_S(V))} \right\} \end{equation}${}$\\
since (for $(a,b,c)\in C_{\de}^{i}$) $0=\text{{D}}(a,b,c)$ $\implies$
$\d a,\,\di b=0$ while $\di c=b-a$. If we project $\ker\{\text{{D}}:\, C_{\de}^{i}\to C_{\de}^{i+1}\}$
to the third $\oplus$-factor of $(10.7)$, then $\text{{im}}\{\text{{D}}:\, C_{\de}^{i-1}\to C_{\de}^{i}\}$
is sent to the subgroup $\di\Gamma(\DI_{S}^{i-2}(V))+\Gamma(C_{2t-i+1}^{S}(V))$;
explicitly for $(A,B,C)\in C_{\de}^{i-1}$ one has $F_{(S)}^{i}\DI_{S}^{i-1}=0$
$\implies$ $B=0$, and $D(A,0,C)=(-\d A,0,\di C+A)$ is sent to $\di C+A$.
Hence we get the desired isomorphism.

\section{\textbf{Reduced higher $AJ$-maps}}

Let $\l^{i}Z^{n}(X_{k})$ denote cycles with class in $\l^{i}CH^{n}(X_{k})$;
we now define maps\[
\chi_{i}:\,\l^{i}Z^{n}(X_{k})\to\uh_{\de}^{i}(\ns,\Lambda\{ i\}).\]
Considering the middle term of $(10.5)$ {[}$\cong\uh_{\de}^{i}(\ns,\Lambda\{ i\})${]}
and the discussion in $\S9$, one might expect $\chi_{i}(\z)$ to
be presented as a $V/\Lambda$-valued differential $(i-1)$-character.
(Recall this is a functional on topological $(i-1)$-cycles with a
special property on bounding cycles.) Accordingly, here is the initial
$ad$ $hoc$ construction of $\chi_{i}$:

Take $\z\in\l^{i}Z^{n}(X_{k})$ with any (fixed) choice of complete
spread $\bar{\zeta},$and let $D\subset S$ (defined $/\bar{\Q}$,
codim. $1$) be such that the restriction of $\bar{\zeta}$ to $X\times(S\m D)$
has relative dimension $0$ (over $S\m D$). Now pick any $\c\in Z_{i-1}^{top}(S\m D).$
Since $\z\in\l^{i},$ $[\zeta]_{i-1}=0\in Hom_{\Q}\left(\uh_{i-1}(\ns,\Q),\, H_{i-1}(X,\Q)\right).$
Therefore $\rho_{X}(\rho_{S}^{-1}(\c))$ has trivial class in $H_{i-1}(X,\Q)$,
and so may be written $\d\Gamma_{\c}$ (for $\Gamma_{\c}\in C_{i}^{top}(X)$
defined up to an $i$-cycle). If $\{\omega_{\ell}\}\subseteq\Omega^{i}(X)$
is a basis for $F^{i}H^{i}(X,\C)$ then the membrane integrals

\begin{equation} \chi_i (\z ) \c := \left\{ \int_{\Gamma_{\c}} \omega_{\ell} \right\} \end{equation}

${}$\\
give a {}``value'' in $V/\Lambda$. If $\c\homeq0$ on $S$, then
there exists $\KK\in C_{i}^{top}(S)$ such that $\c=\d\KK$ and $X\times\KK$
intersects $\bar{\zeta}$ properly (i.e. in the right real codimension).
Again modulo $i$-cycles on $X$, $\Gamma_{\c}=\rho_{X}(\rho_{S}^{-1}(\KK))$
and so

\begin{equation} \chi_i (\z) \d \KK \buildrel \text{mod} \, \Lambda \over \equiv \left\{ \int_{\KK} \rho_S \, _* \rho_X \, ^* \omega_{\ell} \right\} \end{equation}

${}$\\
amounts to integrals of holomorphic $i$-forms (or one $V$-valued
holomorphic $i$-form) over $\KK$. So $(11.1)$ defines a differential
character.

To produce a class in $\uh_{\de}^{i}(\ns,\Lambda\{ i\})$ from this,
we need a corresponding $V$-valued $(i-1)$-current $R$, with $\int_{\c}R\equivmod\chi_{i}(\z)\c;$
that is, one current for each $\omega_{\ell}$. There are in fact
many possible choices of $R$, but (as we will show) their differences
lie in the denominator of the middle term of $(10.5)$. We sketch
one explicit construction (of an $R$) here, and another (for $X$
a product of curves) in $\S14$ which resembles the regulator currents
of $\S\S8-9$ more closely.

Let $\{\gamma^{k}\},\,\{\gamma_{k}^{*}\}$ be dual bases for $H_{2t-i}(S,\Q)$
and $H_{i}(S,\Q)$ (resp.) and $\{\sigma^{m}\}$ be a basis for $H_{2t-i+1}(S,\Q)\,[\cong H_{i-1}(S,\Q)^{\vee}]$.
(The $\gamma$'s and $\sigma$'s are all topological cycles.) For
some $\alpha_{\ell}^{k}\in\C$, $\left[\sum\alpha_{\ell}^{k}\gamma^{k}\right]=\left[\rho_{S}{}_{*}\rho_{X}{}^{*}\omega_{\ell}\right]\in H^{i}(S,\C)$;
and so there exist $R_{\ell}^{0}\in\Gamma(\DI_{S}^{i-1})$ with $\di[R_{\ell}^{0}]=\rho_{S}{}_{*}\rho_{X}{}^{*}\omega_{\ell}-\sum_{k}\alpha_{\ell}^{k}\delta_{\gamma^{k}}.$
Collecting these together, and writing $[\omega_{\ell}]^{\vee}\in V$
for the functional evaluating (for each $m$) to $\delta_{\ell m}$
on $[\omega_{m}]$, $R^{0}:=\{ R_{\ell}^{0}\}=\sum_{\ell}R_{\ell}^{0}\t[\omega_{\ell}]^{\vee}\in\Gamma(\DI_{S}^{i-1}(V))$
has\[
\di[R^{0}]=\{\rho_{S}{}_{*}\rho_{X}{}^{*}\omega_{\ell}\}-\{\sum_{k}\alpha_{\ell}^{k}\delta_{\gamma^{k}}\}\]
\[
=\,\sum_{\ell}(\rho_{S}{}_{*}\rho_{X}{}^{*}\omega_{\ell})\t[\omega_{\ell}]^{\vee}-\sum_{k}\gamma_{k}\t(\sum_{\ell}\alpha_{\ell}^{k}[\omega_{\ell}]^{\vee})\,\,=:\,\O^{0}-\mathsf{C}^{0}.\]
Since $\alpha_{\ell}^{k}=\int_{\gamma_{k}^{*}}\rho_{S}{}_{*}\rho_{X}{}^{*}\omega_{\ell}\,=\,\int_{\rho_{X}(\rho_{S}^{-1}(\gamma_{k}^{*}))}\omega_{\ell}$,
\[
\sum_{\ell}\alpha_{\ell}^{k}[\omega_{\ell}]^{\vee}=\int_{\rho_{X}(\rho_{S}^{-1}(\gamma_{k}^{*}))}(\cdot)\,\in\,\text{im}\{ H_{i}(X,\Q)\to F^{i}H^{i}(X,\C)^{\vee}\}=\Lambda.\]
Therefore $\mathsf{C}^{0}\in Z(C_{2t-i}^{S}(\Lambda))$, and $R^{0}$
has the property $\int_{\c=\d\KK}R^{0}\equivmod\int_{\KK}\O^{0}.$
However, if $\c$ is $not$ a boundary the value $\int_{\c}R^{0}$
needs correcting (in order to be $\equiv\chi_{i}(\z)\c$); an \emph{ad
hoc} modification $R=R^{0}+\left\{ \sum_{m}\beta_{\ell}^{m}\delta_{\sigma^{m}}\right\} $
($\beta_{\ell}^{m}\in\C$) of the current $R^{0}$ does the trick
(for all $\c$). Clearly $\di[R]=\O^{0}-\mathsf{C}^{0}\in Z(C_{2t-i}^{S}(\Lambda))+Z(F_{S}^{i}\DI_{S}^{i}(V))$,
so $R$ gives a class in the middle term of $(10.5)$.

To check well-definedness, let $R'$ be another $(i-1)$-current in$\linebreak$
$\ker\left\{ \Gamma(\DI_{S}^{i-1}(V))\to\frac{Z(\DI_{S}^{i}(V)}{Z(C_{2t-i}^{S}(\Lambda)+Z(F_{S}^{i}\DI_{S}^{i}(V))}\right\} $,
writing $\di R'=\Omega_{0}'-\mathsf{{C}}_{0}'$; and assume that $\int_{\c}R'\maphh\int_{\c}R$
(\emph{a fortiori} $\int_{\d\KK}R'\maphh\int_{\KK}\Omega_{0}$) for
all $\c\in Z_{i-1}^{top}(S\m D)$ (resp. admissible $\KK\in C_{i}^{top}(S)$).
By a $\bar{\d}$-regularity lemma from Hodge theory, a $\di$-closed
current of pure type $(i,0)$ is in fact a holomorphic form; hence
$\Omega_{0},\,\Omega_{0}'\in\Omega^{i}(S)\otimes V$. But $0\maphh\int_{\d\KK}(R-R')=\int_{\KK}\di[R-R']=\int_{\KK}\{(\Omega_{0}-\Omega_{0}')-(\mathsf{{C}}_{0}-\mathsf{{C}}_{0}')\}$
is then incompatible with the fact that continuous deformations of
$\KK$ should continuously change the value of $\int_{\KK}(\Omega_{0}-\Omega_{0}')$,
unless $\Omega_{0}-\Omega_{0}'=0$. Hence $\mathsf{{C}}_{0}-\mathsf{{C}}_{0}'=\di[R-R']$,
and must be $\di$ of some $\mathsf{{C}}_{1}\in\Gamma(C_{2t-i+1}^{S}(\Lambda))$.
This makes $R-R'-\mathsf{{C}}_{1}$ $\di$-closed on $S$ \emph{a
fortiori} on $S\m D$, where $\int_{\c}R-R'-\mathsf{{C}}_{1}\maphh\int_{\c}R-R'\maphh0$
shows its class is {}``$\Lambda$-integral'', thus shared with some
$\mathsf{{C}}_{2}\in\Gamma(C_{2t-i+1}^{S}(\Lambda))$ (with $\d\mathsf{{C}}_{2}$
supported on $D$). We have therefore that $\{ R-R'-(\mathsf{{C}}_{1}+\mathsf{{C}}_{2})\}\in\{\di\Gamma(\DI_{S}^{i-2}(V))+N^{1}\Gamma(\DI_{S}^{i-1}(V))\}$,
and so $R-R'\,\in$ denominator of $(10.5)$'s middle term.

We emphasize an important point: that $\underline{H}_{\de}^{i}(\eta_{S},\Lambda\{ i\})$
can be put into $1-1$ correspondence with {[}the limit over $D$
of{]} a certain class of ($V/\Lambda$-valued) functionals%
\footnote{namely, those that evaluate (mod $\Lambda$) on all $\d\KK$ to $\int_{\KK}\Omega$
for some fixed $\Omega\in\Omega^{i}(S)\otimes V$. These are different
from those encountered in $\S9$ since here $(a)$ $\KK$ need not
avoid $D$ and $(b)$ $\Omega$ cannot have poles along $D$.%
} on topological $(i-1)$-chains on $S\m D$.

Now that we have a well-defined class $\chi_{i}(\z)=R\in\uh_{\de}^{i}(\ns,\Lambda\{ i\})$,
set\[
\overline{[\zeta]_{i}}:=p(\chi_{i}(\z)).\]
 This is essentially $\di[R]$ (in the last term of $(10.5)$), whose
image under $\linebreak$ $\ker\left\{ \uh^{i}(\ns,\Lambda)\oplus F^{i}\uh^{i}(\ns,V)\to\uh^{i}(\ns,V)\right\} \hookrightarrow F^{i}\uh^{i}(\ns,V)$
is just $\O^{0}$, i.e.\[
\sum_{\ell}[\omega_{\ell}]^{\vee}\t[\rho_{S}{}_{*}\rho_{X}{}^{*}\omega_{\ell}]\in Hom_{\C}\left(F^{i}H^{i}(X,\C),\, F^{i}\uh^{i}(\ns,\C)\right).\]

If $\overline{[\zeta]_{i}}$ vanishes then $(\forall\ell)$ $[\rho_{S}{}_{*}\rho_{X}{}^{*}\omega_{\ell}]=0\in F^{i}\uh^{i}(\ns,\C)$
$\implies$ the forms $\rho_{S}{}_{*}\rho_{X}{}^{*}\omega_{\ell}$
are identically zero. For $\c=\d\KK$, $(11.2)$ now becomes $\chi_{i}(\z)\d\KK\equiv0$
(mod $\Lambda$) and $\chi_{i}(\z)$ gives a $V/\Lambda$-valued functional
on $homology$ $classes$, i.e. an element of $\uh^{i-1}(\ns,V/\Lambda)$.
On the level of the exact sequence $(10.5)$, $\O^{0}=0$ $\implies$
$-\di[R]=\mathsf{C}^{0}$ $\implies$ $[\mathsf{C}^{0}]=0\in\uh^{i}(\ns,\Lambda)$
$\implies$ $\mathsf{C}^{0}$ has a primitive $\mathsf{K}^{0}\in\Gamma(C_{2t-i+1}^{S}(\Lambda))$.
We may modify $R\mapsto R+\mathsf{K}^{0}=:\tilde{R}$ (for free) to
get $\tilde{R}\in Z(\DI_{S}^{i-1}(V))$; this slides back along $q$
to give a class in $\uh^{i-1}(\ns,V/\Lambda)$. Of course, $\tilde{R}$
and $R$ do not yield distinct differential characters; $\tilde{R}$
is just the right choice of current.

So if $\overline{[\zeta]_{i}}=0$, define $\overline{[AJ\zeta]_{i-1}}$
by\[
q\overline{[AJ\zeta]_{i-1}}=\chi_{i}(\z)\,;\]
 and use either approach above to compute it.

\begin{rem}
In some situations it may be easier to show $\overline{[\zeta]_{i}}=0$
than $[\zeta]_{i}=0$ (equivalent only if GHC holds). However, for
$i>t=\dim S$ one already knows $[\zeta]_{i}=0$ --- a Lefschetz argument
as in $\S3$ shows $\uh^{i}(\ns,\Q)=0$.
\end{rem}
We now say how all this relates to the invariants of Part $1$: the
present invariants are {}``quotients'' of them. For instance, in
case $[\bar{\zeta}]=0$ (so that $[AJ(\zeta)]_{i-1}$ is computed
by $\int_{\d^{-1}\bar{\zeta}}(\cdot)$) for $\z\in\l^{i}CH_{0}$,
one can easily check the following (see \cite[sec. 4.3]{Ke4}): write
$\underline{H}^{*}(S_{rel}):=\text{{coim}}\{ H^{*}(S_{rel})\to H^{*}(S)\}\cong\varprojlim_{D}\text{{im}}\{ H^{*}(S,D)\to H^{*}(S)\};$
then restricting functionals via$\linebreak$ $\{ F^{t+1}(H^{2t-i+1}(S,\C)\otimes H^{i}(X,\C))\}^{\vee}\twoheadrightarrow\{\underline{H}^{2t-i+1}(S_{rel},\C)\otimes F^{i}H^{i}(X,\C)\}^{\vee}$
induces a (well-defined) map $Gr_{\l}^{i-1}\underline{J}^{n}(X\times\eta_{S})\twoheadrightarrow\text{{Hom}}_{\Q}(\underline{H}_{i-1}(\eta_{S},\Q),\frac{\Omega^{i}(X)^{\vee}}{\text{{im}}\{ H_{i}(X,\Q)\}})$
sending $[AJ(\zeta)]_{i-1}\mapsto\overline{[AJ(\zeta)]_{i-1}}$. However,
we need a stronger result for $\S\S13-14$, where $\bar{\zeta}\homeq0$
cannot necessarily be arranged. Most importantly, we want $\chi_{i}$
to be well-defined on $Gr_{\l}^{i}CH^{n}(X_{k})$ (i.e., modulo $\ratequiv$
and $\l^{i+1}$); and proving this requires precisely the characterization
$(10.6)$. Here, then, is what one can say:

\begin{prop}
On $\l^{i}Z^{n}(X_{k}),$ $\phi_{i}\circ\Psi_{i}=\chi_{i}$.
\end{prop}
\begin{cor}
(i) $\chi_{i}$ kills $\l^{i+1}Z^{n}$ and cycles $\ratequiv0$, i.e.
is defined on $Gr_{\l}^{i}CH^{n}(X_{k})$.

(ii) $\overline{[\zeta]_{i}}={}'\phi_{i}([\zeta]_{i})$, and if $[\zeta]_{i}$
(not $\overline{[\zeta]_{i}}$) vanishes (so that $[AJ\zeta]_{i-1}$
is defined), then $\overline{[AJ\zeta]_{i-1}}={}''\phi_{i-1}([AJ\zeta]_{i-1}).$

(iii) If either $\overline{[\zeta]_{i}}$ or $\overline{[AJ\zeta]_{i-1}}$
is nonzero, then $\left\langle \z\right\rangle \notin\l^{i+1}CH^{n}(X_{k})$.
In particular, $\z\ratnequiv0$.
\end{cor}
\begin{proof}
\emph{(of Proposition)} Given $\z\in\l^{i}Z^{n}(X_{k})$, let $\bar{\zeta}$
be any choice of complete spread; this has Deligne class $c_{\de}(\bar{\zeta})=(\bar{\zeta},\delta_{\bar{\zeta}},0)$
in

\begin{equation} H^{2n}_{\de} (X \times S ,\Q(n)) \cong \frac{\ker \left\{ Z(C^{X\times S}_{2t} \t \Q(n) ) \oplus Z(F^n \DI^{2n}_{X\times S} ) \oplus \Gamma (\DI^{2n-1}_{X\times S} ) \buildrel \text{D} \over \to Z(\DI^{2n}_{X\times S} ) \right\} }{\text{im} \left\{ \Gamma (C^{X\times S}_{2t+1} \t \Q(n) ) \oplus \Gamma ( F^n \DI^{2n-1}_{X\times S} ) \oplus \Gamma (\DI^{2n-2}_{X\times S} ) \buildrel \text{D} \over \to \text{num} \right\} } , \end{equation}

${}$\\
where $\text{{D}}(A,B,C)=(-\d A,-\di[B],\di[C]-B+\delta_{A})$ {[}$\delta_{A}=$
delta-function current via integration along $A${]}. This $c_{\de}(\bar{\zeta})$
maps to $c_{\H}(\zeta)\in\l^{i}\underline{H}_{\H}^{2n}(X\times\eta_{S},\Q(n))$
under $\Phi$. Since $\Phi(\l_{S}^{i}H_{\de}^{2n})=\l^{i}\underline{H}_{\H}^{2n}$
and $\ker(\Phi)=\varinjlim_{D}\{\text{{im}}H_{\de,X\times D}^{2n}(X\times S,\Q(n))\}$
{[}limit over divisors $D_{/\bar{\Q}}\subseteq S${]}, we can modify
$c_{\de}(\bar{\zeta})$ by a D-closed triple $(N_{top},N_{F},N_{\KK})$
supported on some $X\times D$, to get a class in $\l_{S}^{i}H_{\de}^{2n}$
(with the same $\Phi$-image). This means that we can add a D-coboundary
(triple) to the modified class, to get a triple $(\mathsf{{C}},\Omega,\mathsf{{R}})$
in the numerator of $(10.6)$ (i.e. $\mathsf{{C}},\Omega\in\l_{S}^{i}$,
$\mathsf{{R}}\in\l_{S}^{i-1}$, and all three K\"unneth-decomposed).
Explicitly,

\begin{equation} \begin{matrix} \text{D} (\Xi_{top}, \Xi_F , \KK ) + (\bar{\zeta}, \delta_{\bar{\zeta}}, 0)+(N_{top},N_F,N_{\KK}) \, = \\ \sum_{j\geq i} \left( \sum_k \c^{S,k}_{2t-j}\times \c^{X,k}_j , \sum_{\ell} \Omega^j_{S,\ell} \v \Omega^{2n-j}_{X,\ell} , \sum_m R^{j-1}_{S,m} \v R^{2n-j}_{X,m} \right) \\ =: (\mathsf{C}, \Omega, \mathsf{R} ) . \end{matrix} \end{equation}

${}$\\
(where $k,\,\ell,\, m$ are dummy indices and the others indicate
degree). We may assume that the $\c_{*}^{S,k},\,\c_{*}^{X,k}$ resp.
$\Omega_{S,\ell}^{*},\,\Omega_{X,\ell}^{*}$ are $\d$- resp. $\di$-closed,
and consequently that the $R_{X,m}^{2n-i}$'s (but not all $R_{X,m}^{*}$'s)
are also d-closed. Hence, writing $[\cdot]$ for (co)homology class
on $X$,\[
\left(\sum_{k}\c_{2t-i}^{S,k}\t[\c_{i}^{X,k}],\,\sum_{\ell}\Omega_{S,\ell}^{i}\t[\Omega_{X,\ell}^{2n-i}],\,\sum_{m}R_{S,m}^{i-1}\t[R_{X,m}^{2n-i}]\right)\]
gives a class in the middle term of $(10.4)$, and this is $\Psi_{i}(\z)$.

Next, $\phi_{i}$ sends this to the class of $\sum_{m}R_{S,m}^{i-1}\t\overline{[R_{X,m}^{2n-i}]}\in\Gamma(S,\DI_{S}^{i-1}(V))$
in the middle term of $(10.5)$. (Here the bar indicates the quotient
$H^{2n-i}(X,\C)=\widetilde{V}\twoheadrightarrow V=\{ F^{i}H^{i}(X,\C)\}^{\vee}$.)
We must check that this $\phi_{i}(\Psi_{i}(\z))$ gives the same $V/\Lambda$-valued
differential character as $\chi_{i}(\z)$ by integrating it over topological
$(i-1)$-cycles $\c\in Z_{i-1}^{top}(S\m D)$. Let $\{\omega_{\ell}\}\subseteq\Omega^{i}(X)$
be a basis. Since $\left(\sum_{j\geq i+1}\sum_{m}R_{S,m}^{j-1}\wedge R_{X,m}^{2n-j}\right)\wedge\pi_{X}^{*}\omega_{\ell}$
integrates to $0$ on $\c\times X$ (by degree), this amounts to checking
that the vectors $\left\{ \int_{\c\times X}R\wedge\pi_{X}^{*}\omega_{\ell}\right\} $
and $\chi_{i}(\z)\c$ are (for each such $\c$) equal in $V/\Lambda$.
It is thus more than sufficient to show for $any$ $\omega\in\Omega^{i}(X)$
and $\c\in Z_{i-1}^{top}(S\m D)$, that $\int_{\c\times X}R\wedge\pi_{X}^{*}\omega=\int_{\Gamma}\omega$
where $\d\Gamma=\rho_{X}(\rho_{S}^{-1}(\c))$.

From $(11.4)$, $\mathsf{{R}}=\di\KK-\Xi_{F}+\delta_{\Xi_{top}}+N_{\KK}$.
Since $\Xi_{F}\in\Gamma(F^{n}\DI_{X\times S}^{2n-1})$ and $\pi_{X}^{*}\omega\in\Gamma(F^{i}\Omega_{X\times S}^{i})$,
$\Xi_{F}\wedge\pi_{X}^{*}\omega$ has less than $n-1$ $d\bar{z}$'s
and so $\pi_{S}{}_{_{*}}$ (which involves fiberwise integration along
$X$) kills it. Hence, $\int_{\c\times X}\Xi_{F}\wedge\pi_{X}^{*}\omega=\int_{\c}\pi_{S}{}_{_{*}}\{\Xi_{F}\wedge\pi_{X}^{*}\omega\}=0$.
Next, $\c\times X$ does not intersect the support (in $D\times X$)
of $N_{\KK}$, so $\int_{\c\times X}N_{\KK}\wedge\pi_{X}^{*}\omega=0$;
and $\int_{\c\times X}\di\KK\wedge\pi_{X}^{*}\omega=\int_{\c\times X}\di\{\KK\wedge\pi_{X}^{*}\omega\}=0$
since $\omega$ is closed. We have therefore\[
\int_{\c\times X}\mathsf{{R}}\wedge\pi_{X}^{*}\omega\,=\,\int_{(\c\times X)\cap\Xi_{top}}\pi_{X}^{*}\omega\,=\,\int_{\pi_{X}\{(\c\times X)\cap\Xi_{top}\}}\omega.\]
Now from $(11.4)$ again, $\d\Xi_{top}=\bar{\zeta}+N_{top}-\sum_{j\geq i}\sum_{k}\c_{2t-j}^{S,k}\times\c_{j}^{X,k}$;
clearly $N_{top}\cap(\c\times X)=0$ and the same may be assumed for
$\sum_{j\geq i}(\cdots)\cap(\c\times X)$ also (by degree on $S$).
Hence, $\d\pi_{X}\{(\c\times X)\cap\Xi_{top}\}=\pi_{X}\{(\c\times X)\cap\d\Xi_{top}\}=\pi_{X}\{(\c\times X)\cap\bar{\zeta}\}=\rho_{X}(\rho_{S}^{-1}(\c))$.
\end{proof}
\begin{rem}
The reader familiar with the work of S. Saito (e.g., \cite{sS}) may
expect that the $\overline{[\zeta]_{i}}$ identify with so-called
{}``Mumford invariants'' (specialized to the case where $X$ is
defined $/\bar{\Q}$). A version of these comes from filtering the
image of the dlog cycle-class map $CH^{n}(X_{k})\to H^{n}(X_{k},\O_{X_{k}/\bar{\Q}}^{n})[\t\C]$.
(Note that \cite{EP} defined similar invariants. Strictly speaking,
Saito's target is $\HH^{n}(X,\O_{X_{k}/\bar{\Q}}^{\bullet\geq n})$.)
Writing $Gr_{\l}^{i}\O_{X_{k}/\bar{\Q}}^{n}:=\O_{k/\bar{\Q}}^{i}\t\O_{X_{k}/k}^{n-i}$,
one defines a spectral sequence $E_{1}^{i,j}=H^{i+j}(X_{k},Gr_{\l}^{i}\O_{X_{k}/\bar{\Q}}^{n})\,\,\implies\,\, H^{i+j}(X_{k},\O_{X_{k}/\bar{\Q}}^{n}).$
Since $d_{1}$ is just the (graded) arithmetic Gauss-Manin connection
$\bar{\nabla}$, $Gr_{\l}^{i}H^{n}(X_{k},\O_{X_{k}/\bar{\Q}}^{n})=E_{\infty}^{i,n-i}=E_{2}^{i,n-i}$
is the middle cohomology of the complex\[
\to\O_{k/\bar{\Q}}^{i-1}\t H^{n-1}(\O_{X_{k}/k}^{n-i+1})\mapnabla\O_{k/\bar{\Q}}^{i}\t H^{n}(\O_{X_{k}/k}^{n-i})\mapnabla\O_{k/\bar{\Q}}^{i+1}\t H^{n+1}(\O_{X_{k}/k}^{n-i-1})\to.\]
For us $X_{k}=X_{\bar{\Q}}\t k$, so $\bar{\nabla}=0$ and that cohomology
is simply $\O_{k/\bar{\Q}}^{i}\t H^{n}(\Omega_{X_{k}/k}^{n-i})$;
tensored with $\C$, this becomes $\Omega_{k/\bar{\Q}}^{i}\t Gr_{F}^{n-i}H^{2n-i}(X_{\C}^{(an)},\C)$
by GAGA. For $\z\in\l^{i}CH^{n}(X_{k})$, the resulting invariant
is just the image of $\overline{[\zeta]_{i}}$ under the composition\[
F^{i}\uh^{i}(\ns,V)\cap\uh^{i}(\ns,\Lambda)\hookrightarrow F^{i}H^{i}(S^{(an)},\C)\t V\cong H^{0}(\O_{S/\bar{\Q}}^{i})\t V\hookrightarrow\O_{k/\bar{\Q}}^{i}\t V,\]
 where the last map is given by evaluation at the {}``very general''
point $s_{0}\in S_{(\C)}$.
\end{rem}

\section{\textbf{Finer invariants: $[AJ(\zeta)]_{1}^{tr}$ and coniveau}}

Recall from $\S4$ that when $[\bar{\zeta}]$ vanishes, $[AJ\zeta]_{i-1}$
(for any $i$) is computed by Griffiths' prescription $\int_{\d^{-1}\bar{\zeta}}(\cdot)$,
restricted to give a functional on$\linebreak$ $F^{t+1}\{\uh^{2t-i+1}(S_{rel},\C)\t H^{i}(X,\C)\}$
{[}$\subseteq F^{t+1}\{ H^{2t-i+1}(S,\C)\t H^{i}(X,\C)\}${]}. (This
function identifies with an element of $\uh_{\C}^{2n-1}\left/F^{n}\uh_{\C}^{2n-1}\right.$
in the notation of $(10.2)$). In an ideal world (where the Hodge
conjecture is a theorem), this completely describes $[AJ\zeta]_{n-1}$
in particular (or $[AJ\zeta]_{t}$ if $t<n-1$), since $[\zeta]_{0}$
thru $[\zeta]_{n}$ must be zero ($\implies$ $[\zeta]=0$ $\mapii$
$[\bar{\zeta}]=0$) for it to be defined. What we lack, is a precise
description of $Gr_{\l}^{i-1}\underline{J}^{n}(X\times\ns),$ and
towards some $[AJ\zeta]_{1}$ computations in $\S\S15-17$ we now
remedy this for $i=2$.

Specializing the notation of $\S10$ to $i=2$, $H_{\Q}^{2n-1}=H^{1}(\ns,\Q)\t H^{2n-2}(X,\Q)$
and $\uh_{\Q}^{2n-1}=W_{2n-1}H_{\Q}^{2n-1}$; also note that $W_{2n}H_{\Q}^{2n-1}=H_{\Q}^{2n-1}$.
According to $(10.2)$,\[
Gr_{\l}^{1}\underline{J}^{n}(X\times\ns)=\frac{Ext_{_{\text{MHS}}}^{1}(\Q,\uh_{\Q}^{2n-1}(n))}{\im\left\{ Hom_{_{\text{MHS}}}(\Q,\frac{H_{\Q}^{2n-1}}{\uh_{\Q}^{2n-1}}(n))\right\} }=\frac{Ext_{_{\text{MHS}}}^{1}(\Q,\uh_{\Q}^{2n-1}(n))}{\left\{ H_{\Q}^{2n-1}\cap F^{n}H_{\C}^{2n-1}\right\} \cap\uh_{\C}^{2n-1}}.\]
It is the denominator that we want to better compute, namely the image
of

\begin{equation} Hom_{\text{MHS}} \left( \Q(-n) , \, Gr^W_2 H^1(\eta_S,\Q)\t H^{2n-2}(X,\Q) \right) . \end{equation}${}$\\
Since $Gr_{2}^{W}H^{1}(\eta_{S},\Q)$ is of pure type $(1,1)$, recalling
$Hg^{n-1}(X)\,=\, H^{2n-2}(X,\Q)\cap F^{n-1}H^{2n-2}(X,\C)\,=\, Hom_{\text{{MHS}}}(\Q(-n+1),\, H^{2n-2}(X,\Q))$
we have that $(12.1)\,\cong$

\begin{equation} Hom_{\text{MHS}} \left( \Q(-n) , \, Gr^W_2 H^1(\eta_S,\Q)\t Hg^{n-1}(X) \right) . \end{equation}${}$\\
Using now purity of $Hg^{n-1}(X)$, we easily see that the image of
$(12.2)$ {[}hence $(12.1)${]} in $Ext_{\text{{MHS}}}^{1}$ is the
same as the image of

\begin{equation} \left( \left\{ F^1 [W_2] H^1 (\ns ,\C) + [W_2] H^1(\ns ,\Q) \right\} \cap  \uh^1 (\ns ,\C) \right) \t Hg^{n-1} (X) . \end{equation}${}$\\
Since $\dim(X)=n$, the (dualized) Lefschetz $(1,1)$ theorem allows
us to replace $Hg^{n-1}(X)$ by $H_{alg}^{2n-2}(X,\Q)$; similarly
$Hg^{1}(X)=H_{alg}^{2}(X,\Q)$ and we write $H_{tr}^{2}(X,\Q)$ and
$H_{tr}^{2n-2}(X,\Q)$ (resp.) for their orthogonal complements under
the cup product. Of course $\cup$ makes $H_{alg}^{2}$ and $H_{alg}^{2n-2}$,
hence $H_{tr}^{2}$ and $H_{tr}^{2n-2}$, duals.

Now in $\uh_{\C}^{2n-1}\left/F^{n}\uh_{\C}^{2n-1}\right.\,\cong$
\[
\left\{ F^{t}\uh^{2t-1}(S_{rel},\C)\t F^{1}H^{2}(X,\C)\,+\,[F^{t-1}]\uh^{2t-1}(S_{rel},\C)\t F^{2}H^{2}(X,\C)\right\} ^{\vee},\]
the image of $(12.3)$ may be seen as functionals arising from integrals
over \[
\{(2t-1)\text{-chains\, bounding\, on\, divisors\, defined\, }/\bar{\Q}\,\,\subset S\}\times\{\text{curves\, on\, }X\}.\]
 However for our purposes it will be acceptable to quotient out by
$\uh^{1}(\ns,\C)\t H_{alg}^{2n-2}(X,\C)$ (which contains $(12.3)$).
Writing $\uh_{\Q,\, tr}^{2n-1}=\uh^{1}(\ns,\Q)\t H_{tr}^{2n-2}(X,\Q)$
and so on, set\[
\left(Gr_{\l}^{1}\underline{J}^{n}(X\times\ns)\right)_{tr}:=\frac{Gr_{\l}^{1}\underline{J}^{n}(X\times\ns)}{\uh^{1}(\ns)\t H_{alg}^{2n-2}(X)}\cong\frac{\uh_{\C,\, tr}^{2n-1}}{\uh_{\C,\, tr}^{2n-1}+F^{n}\uh_{\C,\, tr}^{2n-1}}\]

\begin{equation} \cong \begin{matrix} \underline{ \left\{ F^t \uh ^{2t-1} (S_{rel} ,\C) \t F^1 H^2_{tr} (X,\C) + \uh^{2t-1} (S_{rel} ,\C) \t F^2 H^2 (X,\C) \right\} ^{\vee} } \\ \text{im} \left\{ \uh _{2t-1} (S_{rel} ,\Q) \t H^{tr}_2 (X,\Q) \right\} \end{matrix} \end{equation}

${}$\\
where $H_{2}^{tr}(X,\Q)$ is Poincar\'e  duals of $H_{tr}^{2n-2}(X,\Q)$
and $\uh_{*}(S_{rel},\Q):=\text{im}\{ H_{*}(S)\to H_{*}(S_{rel})\}.$
This is a much less drastic quotient than \[
\uh^{1}(\ns,V/\Lambda)\cong\left\{ \uh^{2t-1}(S_{rel},\C)\t F^{2}H^{2}(X,\C)\right\} ^{\vee}\left/\text{im}\left\{ \uh_{2t-1}(S_{rel},\Q)\t H_{2}^{tr}(X,\Q)\right\} \right..\]
We write $[AJ(\zeta)]_{1}^{tr}$ for the (projected) image of $[AJ(\zeta)]_{1}$
in $(12.4)$.

In the event that $\boxed{t\,(=\dim S)=1}$, $\uh^{1}(S_{rel})=\text{coim}\{ H^{1}(S_{rel})\to H^{1}(S)\}\cong H^{1}(S)$
and $(12.4)$ simplifies to\[
\left(Gr_{\l}^{1}\underline{J}^{n}(X\times\ns)\right)_{tr}=\begin{array}{c}
\underline{\left\{ F^{1}H^{1}(S,\C)\t F^{1}H_{tr}^{2}(X,\C)\,+\, H^{1}(S,\C)\t F^{2}H^{2}(X,\C)\right\} ^{\vee}}\\
\text{im}\left\{ H_{1}(S,\Q)\t H_{2}^{tr}(X,\Q)\right\} \end{array}\]

\begin{equation} = \begin{matrix} \underline{ \left\{ H^{1,0} (S,\C) \t H^{1,1}_{tr} (X,\C) \oplus H^1 (S,\C) \t H^{2,0} (X,\C) \right\} ^{\vee} } \\  \text{periods}  \end{matrix} . \end{equation}

${}$\\
Again, $[AJ(\zeta)]_{1}^{tr}\in(12.5)$ is the image of $[AJ(\zeta)]_{1}$
(which in this case $=\Psi_{2}(\z)$ since $Gr_{\l}^{2}\underline{Hg}^{n}(X\times\ns)=\{0\}$
for $t=1$). Note that the $\Q$-dimension of the period lattice equals
the $\R$-dimension of the numerator.

We summarize the picture for $i>2$. Write $F_{h}^{n-i+1}H^{2n-i}(X,\Q)$
for the largest sub-Hodge-structure of $H^{2n-i}(X,\Q)$ contained
in $H^{2n-i}(X,\Q)\cap F^{n-i+1}H^{2n-i}(X,\C)$, and $F_{h}^{n-i+1}H^{2n-i}(X,\C)$
for its complexification. The analogue of $(12.1)$ is then

\[
Hom_{\text{{MHS}}}\left(\Q(-n),\, Gr_{i}^{W}H^{i-1}(\eta_{S},\Q)\t H^{2n-i}(X,\Q)\right)\cong\]

\begin{equation} Hom_{\text{MHS}} \left( [ Gr^W_{2t-i} H^{2t-i+1} (S_{rel},\Q)] \t \Q(t-n), H^{2n-i} (X,\Q ) \right) . \end{equation}${}$\\
Since $F^{i}H^{i-1}(\eta_{S},\C)=\{0\}$, dually $H^{2t-i+1}(S_{rel},\C)=F^{t-i+1}H^{2t-i+1}(S_{rel},\C)$;
hence picking any $\theta\in(12.6)$, $\text{{im}}(\theta)\subseteq F_{h}^{n-i+1}H^{2n-i}(X,\Q)$.
It follows that $(12.6)\cong$\[
Hom_{\text{{MHS}}}\left([Gr_{2t-i}^{W}H^{2t-i+1}(S_{rel},\Q)]\t\Q(t-n),\, F_{h}^{n-i+1}H^{2n-i}(X,\Q)\right)\cong\]

\begin{equation} Hom_{\text{MHS}} \left( \Q(-n), \,  Gr^W_i H^{i-1} (\eta_S,\Q)\t F^{n-i+1}_h H^{2n-i}(X,\Q)  \right) , \end{equation}${}$\\
the analogue of $(12.2)$. Writing $\mathcal{{H}}_{\Q}:=W_{i}H^{i-1}(\eta_{S},\Q)\otimes F_{h}^{n-i+1}H^{2n-i}(X,\Q)$,
the image of $(12.7)$ {[}hence $(12.6)${]} in $Ext_{\text{{MHS}}}^{1}(\Q(-n),\,\underline{H}^{i-1}(\eta_{S},\Q)\t H^{2n-i}(X,\Q))$
is then given by\[
\left\{ F^{n}\mathcal{{H}}_{\C}+\mathcal{{H}}_{\Q}\right\} \,\,\,\cap\,\,\,\underline{H}^{i-1}(\eta_{S},\Q)\t F_{h}^{n-i+1}H^{2n-i}(X,\Q)\]
 (compare to $(12.3)$), and so $im\{ Hom\}$ will be absorbed if
we quotient $Gr_{\l}^{i-1}\underline{J}^{n}(X\times\ns)$ by $\uh^{i-1}(\ns,\C)\t F_{h}^{n-i+1}H^{2n-i}(X,\C).$
It is a very important point, that this is usually (substantially)
smaller than the $\uh^{i-1}(\ns,\C)\t F^{n-i+1}H^{2n-i}(X,\C)$ by
which we quotiented in $\S10$ (to get $\uh^{i-1}(\ns,V/\Lambda)$).
In fact, provided GHC holds, it is precisely $\uh^{i-1}(\ns,\C)\t N^{n-i+1}H^{2n-i}(X,\C)$;
while {[}hard Lefschetz $\implies${]} $H^{2n-i}(X,\C)=N^{n-i}H^{2n-i}(X,\C).$

Writing $Gr_{N_{X}}^{n-i}\uh_{\Q}^{2n-1}:=\uh^{i-1}(\ns,\Q)\t\frac{N^{n-i}}{N^{n-i+1}}H^{2n-i}(X,\C)$,
GHC $\implies$ $\exists$ of a quotient 

\begin{equation} \left( Gr^{i-1}_{\l} \underline{J} ^n (X\times \ns) \right) _N := \begin{matrix} Gr_{N_X}^{n-i} \uh^{2n-1}_{\C} \\ \overline{Gr_{N_X}^{n-i} \uh^{2n-1}_{\Q} + F^n \text{(num)}} \end{matrix} \end{equation}

${}$\\
of $Gr_{\l}^{i-1}\underline{J}^{n}(X\times\ns).$ This is one way
to understand the notion that $Gr_{\l}^{i}CH^{n}(X_{k})$ should be
{}``controlled'' by $Gr_{N}^{n-i}H^{2n-i}(X)$.

\part{Applications to products of curves}

\section{\textbf{General observations}}

We now specialize to $X=C_{1}\times\cdots\times C_{n}$, where $C_{j}$
are smooth projective curves (of positive genus) defined $/\bar{\Q}$,
with fixed base points $o\in C_{j}(\bar{\Q})$. 

Two observations are in order. 

The first is that the product structure is reflected in the Leray
filtration on $CH^{n}(X_{k})$. As in $\S7$ one uses the $\frac{n!}{i!(n-i)!}$
projections (one for each $\sigma\in\mathfrak{{S}}_{i}^{n}$) \[
\pi_{\sigma}:\, X\longrightarrow\, C_{\sigma(1)}\times\cdots\times C_{\sigma(i)}=:X^{\sigma}\]
 to products of $i$ curves, to induce isomorphisms\[
H^{i}(X,\Q)\begin{array}{c}
_{\oplus_{\sigma}(\pi_{\sigma})^{*}}\\
\longleftarrow\\
^{\cong}\end{array}\oplus_{\sigma}H^{i}(X^{\sigma},\Q)\]
$\vspace{-5mm}$\[
H^{2n-i}(X,\Q)\begin{array}{c}
_{\oplus_{\sigma}(\pi_{\sigma})_{*}}\\
\longrightarrow\\
^{\cong}\end{array}\oplus_{\sigma}H^{i}(X^{\sigma},\Q)\]
and (using functoriality and the formulas for the $Gr_{\l}^{*}\uH^{*}$)
a diagram

\begin{diagram}
\l^i CH^n(X_k) & \rTo^{\oplus (\pi_{\sigma})_* } & \oplus_{\sigma} \l^i CH^i (X^{\sigma_i}_k ) \\
\dTo>{\Psi_i (n)} & & \dTo>{\oplus \Psi_i(i)} \\
Gr^i_{\l} \uH^{2n} (X\times \ns ,\Q(n) ) & \rTo^{\cong} & \oplus _{\sigma} Gr^i_{\l} \uH^{2i} (X^{\sigma} \times \ns ,\Q(i) ) . \\
\end{diagram}${}$\\
Setting (for $0\leq i\leq n$)\[
F_{\times}^{i+1}CH^{n}(X_{k}):=\ker\left\{ \oplus_{\sigma\in\mathfrak{{S}}_{i}^{n}}(\pi_{\sigma})_{*}\right\} \subseteq CH^{n}(X_{k})\]
 (so $F_{\times}^{n+1}CH^{n}=\{0\}$), one has from this diagram\[
F_{\times}^{i}CH^{n}(X_{k})\subseteq\l^{i}CH^{n}(X_{k})\]
 and equality if $\text{BBC}^{q}$ holds (or without $\text{BBC}^{q}$
for $i\leq2$), since the $\Psi_{i}(i)$ are injective in that case.
Here the $\l^{i}$ are defined, as usual, via kernels of successive
$\Psi_{i}$.

In particular, we shall concentrate on\[
F_{\times}^{n}CH^{n}(X_{k}):=\left.\left\{ \begin{array}{c}
0\text{-cycles\, whose\, projections\, to}\\
\text{products\, of\, }n-1\text{\, curves,\, are\, }\ratequiv0\end{array}\right\} \right/\ratequiv\,\,,\]
which obviously contains all {}``$n$-box'' cycles\[
\z=B(p_{1},\ldots,p_{n}):=((p_{1})-(o))\times\cdots\times((p_{n})-(o))\]
 where $p_{j}\in C_{j}(k)$. We note that $\chi_{n}$ is (well-)defined
on $F_{\times}^{n}CH^{n}(X_{k})$ since it is already defined on $\l^{n}$
by Corollary $11.3$ (more generally, $\chi_{i}$ is defined on $F_{\times}^{i}$).

\begin{rem*}
The above $n$-box cycles $/k$ do not necessarily generate$\linebreak$$F_{\times}^{n}CH^{n}(X_{k})$;
the latter also contains certain sums of $n$-box cycles $with$ $p_{j}\in C_{j}(L)$,
$L/k$ $algebraic$, which are essentially {}``norms'' (relative
to $k$) of $n$-box cycles in $F_{\times}^{n}CH^{n}(X_{L})$. (We
could ignore this possibility in Part $2$ because of the existence
of a norm $K_{n}^{M}(L)\to K_{n}^{M}(k)$: such cycles could be replaced
by ones with coordinates $\in k$.) If $k$ was algebraically closed
then of course this problem would disappear, but we have restricted
in this paper to $k/\bar{\Q}$ finitely generated (in order to have
$S$ available; one may then deal with $\bar{k}$ via a limit).
\end{rem*}
The second observation is this: when $\z\left/\{ k\cong\Q(S)\}\right.$
is an $n$-box we can explicitly compute $\chi_{n}(\z)$ as a $V/\Lambda$-valued
$(n-1)$-current {}``$R$'' on $S$. This seems particularly intriguing
when the spread $\zeta_{0}$ of the leading corner $(p_{1},\ldots,p_{n})$
of $\z$, projects to a nontrivial correspondence $\pi_{X}(\zeta_{0})\subseteq C_{1}\times\cdots\times C_{n}.$

The idea is to extend the homotopy technique described in $\S8$ for
$(\C^{*})^{n}$, to $C_{1}\times\cdots\times C_{n}$ by using the
fact that each $C_{j}$ can be cut into a fundamental domain $\dd_{j}$.
One then lifts the {}``$n$-box graph'' $\zeta=B(\zeta_{0})$ to
$(\dd_{1}\times\cdots\times\dd_{n})\times S$ and {}``successively''
contracts $\dd_{j}$-factors (to $o$, in a sense that will be clear
to readers of \cite{Ke1}). The lift traces out a $(2t+1)$-chain
$\Gamma_{\zeta}\subseteq S\times X$ bounding on $\zeta$ and $\{$topological
$(2t-n)$-cycles on $S\}\times\{$topological $n$-cycles on $X\}$.
The function of this chain is essentially to standardize the membranes
$\Gamma_{\c}$ used (see $(11.1)$) in constructing the differential
character $\chi_{n}$, or alternatively the choice of current $R\in\Gamma(\DI_{S}^{n-1}(V))$
representing it. (We are certainly $not$ claiming that the topological
$(2t-n)$-cycles lie in any divisor $D\subseteq S$.)

Specifically, the currents we write down in $\S14$ (without further
proof) will correspond essentially to $\Gamma_{\c}:=\pi_{X}\{(X\times\c)\cap\Gamma_{\zeta}\}.$
They will give classes in the middle term of $(10.5)$, which is to
say $\di[R]$ will be $(i)$ a $\Lambda$-valued $(2t-n)$-cycle plus
$(ii)$ a $V$-valued holomorphic current {[}or form{]}. Either $(i)$
or $(ii)$ can be considered to represent $\overline{[\zeta]_{n}}$,
but we shall write $\overline{[\zeta]_{n}}$ for $(ii)$, and moreover
$\chi_{n}(\z)$ in lieu of $R$. So\[
\di[\chi_{n}(\z)]=\overline{[\zeta]_{n}}+T_{\zeta}\]
 now expresses the fact that $\chi_{n}(\z)$ is a differential character.

When $\overline{[\zeta]_{n}}=0$ we instead write $\overline{[AJ\zeta]_{n-1}}$
for $R$, so that \[
\di\overline{[AJ\zeta]_{n-1}}=T_{\zeta}.\]
 It is unnecessary to correct $R$'s discrepancy from closedness here,
since integrating it over topological cycles$(\t\Q)$ already gives
a cohomology class:\[
\overline{[AJ\zeta]_{n-1}}\in Hom\left(\uh_{n-1}(\ns,\Q),\,\frac{\{ F^{n}H^{n}(X,\C)\}^{\vee}}{\text{im}\{ H_{n}(X,\Q)\}}\right)\cong\uh^{n-1}(\ns,V/\Lambda).\]
 Recall that $\uh_{n-1}(\ns,\Q)$ means cycles which may be moved
to avoid any $D/\bar{\Q}$. We just note that if $S$ is a curve,
then $\uh_{1}(\ns,\Q)\cong H_{1}(S,\Q)$, while if $S$ is a surface
then $\uh_{2}(\ns,\Q)\cong H_{2}^{tr}(S,\Q)$.

Before embarking on computations, we indicate how to see Mumford's
theorem easily using $\overline{[\zeta]_{n}}$, provided one believes
Proposition $(3.2)$. We only do $X=C_{1}\times C_{2}$ (i.e., $n=2$),
and leave the generalization to the reader. To prove no $\pi_{N}$
surjects, it is enough to show\[
\mathcal{B}_{N}:\, S^{(N)}(C_{1}\times C_{2})\longrightarrow\l^{2}CH^{2}(X_{\C})\]
\[
Sym_{N}\left((p_{1},q_{1}),\ldots,(p_{N},q_{N})\right)\longmapsto\sum_{j=1}^{N}B(p_{j},q_{j})\]
 is surjective for no $N$.%
\footnote{or equivalently: for no $N$ is $\text{im}(\mathcal{B}_{N})$ closed
under addition. The smallest subgroup of $\l^{2}CH^{2}(X_{\C})$ containing
$\text{im}(\mathcal{B}_{N})$ (for any $N\geq1$), after all, is all
of $\l^{2}CH^{2}(X_{\C})$.%
}

We do this by producing (for any $N$) a cycle in $\text{im}(\mathcal{B}_{N+1})\m\text{im}(\mathcal{B}_{N})$.
Take $p_{1},\ldots,p_{N+1}\in C_{1}$ and $q_{1},\ldots,q_{N+1}\in C_{2}$
to be {}``algebraically independent'' general points, in the sense
that $\bar{\Q}\left((C_{1}\times C_{2})^{\times(N+1)}\right)\cong\bar{\Q}(p_{1},q_{1},\ldots,p_{N+1},q_{N+1})=:k$
(so $trdeg(k/\bar{\Q})=2N+2$); set $\z_{N+1}:=\sum_{j=1}^{N+1}B(p_{j},q_{j})$.

Recall that $V=H^{0}(\O_{X}^{2})^{\vee}$. Write $\{\omega_{\ell}\}\in H^{0}(\O_{X}^{2})$
for a basis, and $\pi_{j}:\, S=(C_{1}\times C_{2})^{\times(N+1)}\to C_{1}\times C_{2}$
for projection to the $j^{\text{th}}$ factor. Clearly one may take
$\bar{\zeta}_{N+1}:=\sum_{j=1}^{N+1}(id_{X}\times\pi_{j})^{*}\Delta_{X}\,\in Z^{2}((C_{1}\times C_{2})^{\times(N+2)})$;
and so \[
\overline{[\zeta_{N+1}]_{2}}\,=\,\,\sum_{j=1}^{N+1}\left(\sum_{\ell}[\omega_{\ell}]^{\vee}\t\pi_{j}^{*}\omega_{\ell}\right)\,\,=\,\,\sum_{\ell}[\omega_{\ell}]^{\vee}\t\left(\sum_{j=1}^{N+1}\pi_{j}^{*}\omega_{\ell}\right)\]
 gives a $nondegenerate$, anti-symmetric bilinear form $\th_{S,\, s_{0}}^{1}\times\th_{S,\, s_{0}}^{1}\to V$
(where $\theta_{S,\, s_{0}}^{1}$ denotes the holomorphic tangent
space). More precisely, writing $J=\left(\begin{array}{cc}
0 & -1\\
1 & 0\end{array}\right)$, its matrix is proportional (by a vector $\in V$) to $\left(\begin{array}{ccc}
J &  & 0\\
 & \ddots\\
0 &  & J\end{array}\right),$ which has nonzero determinant.

Now any $\z_{N}\in\text{im}(\mathcal{B}_{N})$ is defined over a field
$k_{0}$ with $trdeg(k_{0}/\bar{\Q})\leq2N$. Set $\z=\z_{N+1}-\z_{N}$;
we must show that $\overline{[\zeta]_{2}}\neq0$.

Consider the simplest case $k_{0}\subset k$, so that $S\mappia S_{0}$
and $\overline{[\zeta]_{2}}=\overline{[\zeta_{N+1}]_{2}}-\pi_{0}^{*}\overline{[\zeta_{N}]_{2}}.$
Since $\dim(S_{0})<2N+2,$ it is impossible for $\pi_{0}^{*}$ of
a holomorphic $2$-form $S_{0}$ to be nondegenerate on $S$, and
so $\overline{[\zeta_{N+1}]_{2}}\neq\pi_{0}^{*}\overline{[\zeta_{N}]_{2}}$
as desired.

More generally $\zeta$ lives on $X\times\eta_{S'}$, where $S\mappib S'\mappia S_{0}$
(with $\dim S'=:T$), and $\zeta=\pi^{*}\zeta_{N+1}-\pi_{0}^{*}\zeta_{N}.$
One easily shows the existence of a subspace of $\th_{S',\, s_{0}'}^{1}$
of dimension at least $T-N$, on which $\pi_{0}^{*}\overline{[\zeta_{N}]_{2}}$
restricts to the zero form. The maximum dimension of such a subspace
for $\pi^{*}\overline{[\zeta_{N+1}]_{2}}$ is $T-N-1$. So $\overline{[\zeta]_{2}}=\pi^{*}\overline{[\zeta_{N+1}]_{2}}-\pi_{0}^{*}\overline{[\zeta_{N}]_{2}}\neq0$
and Mumford is proved.

\begin{rem*}
It would $not$ have done, to merely exhibit infinitely many linearly
independent $\overline{[\zeta]_{2}}$-classes. this only establishes
that $\l^{2}CH^{2}(X_{\C})$ contains an $\infty$-dimensional $\Q$-vector-space,
which is not at all Mumford's theorem (see the warning in Remark $8.2(ii)$).
So this is why we needed the Prop. $(3.2)$.
\end{rem*}

\section{\textbf{Formulas for the reduced invariants}}

We now produce the currents representing $\chi_{n}(\z)$, beginning
with $n=3$ (from which the general case will be clear). We then turn
to $n=2$ and perform a computation.

Let $p,q,r$ be points in the respective $C_{i}(k)$ {[}$i=1,2,3${]},
and consider the {}``$3$-box'' cycle $\z=B(p,q,r)=$\[
(p,q,r)-(o,q,r)-(p,o,r)-(p,q,o)+(o,o,r)+(o,q,o)+(p,o,o)-(o,o,o)\]
in $Z^{3}\left((C_{1}\times C_{2}\times C_{3})_{k}\right)$ with class
$\left\langle \z\right\rangle \in\l^{3}CH^{3}\left((C_{1}\times C_{2}\times C_{3})_{k}\right).$
As usual we write $\zeta\in Z^{3}\left((X\times\ns)_{\bar{\Q}}\right)$
for its $\bar{\Q}$-spread, where $X=C_{1}\times C_{2}\times C_{3}$
and $\bar{\Q}(S)\cong k$.

Let $\a_{i}^{j=1,\ldots,2g_{i}}$ (with $g_{i}$ the genus) be topological
cycles on $C_{i}$ spanning $H_{1}(C_{i},\Q)$ and all based at $o$.
We take $\hat{\a}_{i}^{j}$ to be a dual basis of {}``markings''
avoiding $o$, and based at another point which the $\a_{i}^{j}$
avoid. Things are set up as usual so that $[\hat{\a}_{i}^{j}]=[\a_{i}^{j\pm g_{i}}]$
(but $\hat{\a}_{i}^{j}\neq\a_{i}^{j\pm g_{i}}$), and $\a_{i}^{j_{1}}\cap\hat{\a}_{i}^{j_{2}}$
is empty unless $j_{1}=j_{2}$ (and then it is a point). The interior
of our fundamental domain (with {}``center'' at $o$) is then given
by $C_{i}\m\cup_{j}\hat{\a}_{i}^{j}$, and we denote its closure by
$\dd_{i}$.

The forms $\{\w_{i}^{k=1,\ldots,g_{i}}\}\in\O^{1}(C_{i})$ integrate
to give functions $z_{i}^{k}$ on $\dd_{i}$ which are zero at $o$.
They may be viewed as discontinuous, complex-valued functions on $C_{i}$
with cuts at the markings, i.e. as $0$-currents like $\log z$ on
$\C^{*}$, so that\[
\di[z_{i}^{k}]=\w_{i}^{k}-\sum_{j=1}^{2g_{i}}\w_{i}^{k}(\a_{i}^{j})\cdot\delta_{\hat{\a}_{i}^{j}}\]
where $\w(\a):=\int_{\a}\w$. Now the spreads of the points $p,\, q,\, r$
give maps $s_{i}:\, S\to C_{i}$, and by composition with the $z_{i}^{k}$,
zero-currents $f_{i}^{k}:=z_{i}^{k}\circ s_{i}\in\DI^{0}(S)$, or
functions on $S$ with branch cuts $T_{i}^{j}=s_{i}^{-1}(\hat{\a}_{i}^{j}).$
We have\[
\di[f_{i}^{k}]=df_{i}^{k}-\sum_{j=1}^{2g_{i}}\w_{i}^{k}(\a_{i}^{j})\cdot\delta_{T_{i}^{j}}\]
where $df_{i}^{k}:=s_{i}^{*}\w_{i}^{k}$. (This is of course the analogue
of the formula $\di[\log f]=\dlog f-2\pi i\cdot\delta_{T_{f}}$.)

Recalling $V:=\left(F^{3}H^{3}(X,\C)\right)^{\vee}\cong\left(\O^{1}(C_{1})\t\O^{1}(C_{2})\t\O^{1}(C_{3})\right)^{\vee}$
and $\Lambda=\text{im}\{ H_{3}(X,\Q)\}\subset V$, the $V[/\Lambda]$-valued
$2$-current representing $\chi_{3}(\z)$ is

\begin{equation} \sum_{k_1,k_2,k_3} [\omega^{k_1}_1 \v \omega^{k_2}_2 \v \omega^{k_3}_3 ]^{\vee} \t \left( \begin{matrix} f^{k_1}_1 df^{k_2}_2 \v df^{k_3}_3 + \sum _{j_1 =1}^{2g_1} (\omega_1^{k_1} (\alpha ^{j_1}_1 )) f^{k_2}_2 df^{k_3}_3 \cdot \delta_{T^{j_1}_1} \\ + \sum _{j_1 ,j_2} (\omega^{k_1}_1 (\alpha ^{j_1}_1 )) ( \omega^{k_2}_2 (\alpha ^{k_2}_2 )) f^{k_3}_3 \cdot \delta _{T^{j_1}_1 \cap T^{j_2}_2 } \end{matrix} \right) . \end{equation}

${}$\\
(An analogy to the Milnor regulator current $\log f_{1}\dlog f_{2}\v\dlog f_{3}+2\pi i\log f_{2}\dlog f_{3}\cdot\delta_{T_{f_{1}}}-4\pi^{2}\log f_{3}\cdot\delta_{T_{f_{1}}\cap T_{f_{2}}},$
where $f_{j}\in\bar{\Q}(S)$ are viewed as functions from $\ns\to\C^{*}$,
is implied here.) Notice that $\di$ of this current is the $\Lambda$-valued
$3$-chain\[
T_{\zeta}:=\sum_{k_{1},k_{2},k_{3}}[\omega_{1}^{k_{1}}\v\omega_{2}^{k_{2}}\v\omega_{3}^{k_{3}}]^{\vee}\t\left(\sum_{j_{1},j_{2},j_{3}}(\omega_{1}^{k_{1}}(\alpha_{1}^{j_{1}}))(\omega_{2}^{k_{2}}(\alpha_{2}^{j_{2}}))(\omega_{3}^{k_{3}}(\alpha_{3}^{j_{3}}))\cdot\delta_{T_{1}^{j_{1}}\cap T_{2}^{j_{2}}\cap T_{3}^{j_{3}}}\right)\]
\[
=\sum_{j_{1},j_{2},j_{3}}\left(\int_{\alpha_{1}^{j_{1}}\times\alpha_{2}^{j_{2}}\times\alpha_{3}^{j_{3}}}(\cdot)\right)\cdot\delta_{T_{1}^{j_{1}}\cap T_{2}^{j_{2}}\cap T_{3}^{j_{3}}}\]
plus the holomorphic ($V$-valued) $3$-form\[
\overline{[\zeta]_{3}}=\sum_{k_{1},k_{2},k_{3}}[\omega_{1}^{k_{1}}\v\omega_{2}^{k_{2}}\v\omega_{3}^{k_{3}}]^{\vee}\t df_{1}^{k_{1}}\v df_{2}^{k_{2}}\v df_{3}^{k_{3}}.\]
If the latter is zero then the current $(14.1)$ gives (by operating
on topological $2$-cycles) a well-defined class $\overline{[AJ\zeta]_{2}}\in\uh^{2}(\ns,V/\Lambda).$

Similarly, for $\z=B(p,q)\in Z^{2}\left((C_{1}\times C_{2})_{k}\right)$
one has for $\overline{[AJ\zeta]_{1}}$

\small \begin{equation} \sum _{k_1, k_2} [\omega^{k_1}_1 \v \omega^{k_2}_2 ]^{\vee} \t \left( f^{k_1}_1 df^{k_2}_2+\sum_{j_1 =1}^{2g_1} (\omega ^{k_1}_1 (\alpha ^{j_1}_1 )) f^{k_2}_2 \cdot \delta _{T^{j_1}_1} \right) \in \uh^1 \left( \ns, \{ F^2 H^2 (C_1 \times C_2 ,\C) \}^{\vee} \left/ \Lambda \right. \right) \end{equation} \normalsize

${}$\\
provided the $df_{1}^{k_{1}}\v df_{2}^{k_{2}}$ vanish, as for instance
is the case when $\dim S=trdeg(k/\bar{\Q})=1.$ We now give an example
of such a case, where $\overline{[AJ\zeta]_{1}}\neq0$ (and thus $\z\ratnequiv0$).

We mentioned in $\S10$ that $\dim_{\Q}\Lambda$ might be larger than
$\dim_{\R}V\,[=g_{1}g_{2}]$; one situation where this does $not$
happen is for $C_{1}=C_{2}=E$ an elliptic curve (defined $/\bar{\Q}$)
with $complex$ $multiplication$. To be concrete, let $E=\mathcal{V}(\mathsf{F})\subset\P^{2}$
be defined by the equation \[
\mathsf{F}(x,y,z)=y^{2}z-x^{3}+5xz^{2}=0,\]
with $\omega=\text{Res}_{E}\left(\frac{z\, dx\v dy-y\, dx\v dz+x\, dy\v dz}{\mathsf{F}(x,y,z)}\right)\in\O^{1}(E)$
(in affine coordinates $\omega=dx/y$). For the $0$-currents resulting
from the integrals $\int_{o}^{*}\omega$ on $\dd$ we will write simply
$z_{1},z_{2}$; these are just the standard {}``plane coordinates''
on the factors $E\times E$. We choose the topological $1$-cycles
$\alpha,\,\beta$ in such a way that $\O_{1}=\O=\int_{\alpha}\omega$
and $\Omega_{2}=i\O=\int_{\beta}\omega$, where $\O\in\R$ is some
transcendental number. We then have $\int_{\alpha\times\beta}dz_{1}\v dz_{2}=(i\O)\O=i\O^{2}$,
$\int_{\beta\times\beta}dz_{1}\v dz_{2}=(i\O)(i\O)=-\O^{2}$, $\int_{\alpha\times\alpha}dz_{1}\v dz_{2}=\O^{2}$,
which give the period lattice $\Q\left\langle \O^{2},i\O^{2}\right\rangle =\Lambda\subset V=\{ F^{2}H^{2}(E\times E)\}^{\vee}=\left\langle dz_{1}\v dz_{2}\right\rangle ^{\vee}=\C.$

Now take a base point $o$ defined over $\bar{\Q}$ (say, the point
at $\infty$), and a general point $p$ defined over $k\cong\bar{\Q}(E)$
(so that $S$ and $E$ will be birationally equivalent). There also
exists a {}``nontorsion'' point $\xi\in E(\bar{\Q})$, which means
that $\int_{o}^{\xi}\omega$ is nonzero in $\C\left/\Q\left\langle \O,i\O\right\rangle \right.$.
(In fact, for this particular curve one can take $o,\,\xi\,\in E(\Q)$.)
We study the cycle\[
\z=B(p,\, p-\xi)=(p,\, p-\xi)-(o,\, p-\xi)-(p,\, o)+(o,\, o),\]
 with class $\left\langle \z\right\rangle \in\l^{2}CH^{2}\left((E\times E)_{k}\right)$
and spread $\bar{\zeta}\subseteq E\times E\times E$, and $\overline{[\zeta]_{2}}=0$
(since $\dim S=1$).

The spreads of $p$ and $p-\xi$ give rise to maps $(s_{1},s_{2})=(id,\, id-\xi):\, S\simeq E\to E\times E$.
We choose our cuts $\hat{\a},$ $\hat{\B}$ so that $z=\int_{o}\w$
$[=z(o)]$ takes values in the square $-\O/2\leq\Im(z)\leq\O/2,\,\,-\O/2\leq\Re(z)\leq\O/2$
centered at $0$, and pick $\a,\,\B$ to have support along the real
and imaginary axes, respectively. Composing $(id,\, id-\xi)$ with
$z$ gives zero-currents $f$ and $g$ on $S\simeq E$, with cuts
$T_{f}^{\a}$, $T_{f}^{\B}$ (and $T_{g}^{\a}$, $T_{g}^{\B}$). Using
$df$, $dg$ (both $=dz$) for the pullbacks of $\w$ (by $s_{1},\, s_{2}$),
we can now write down the basic $1$-current (from $(14.2)$)\[
\overline{[AJ\zeta]_{1}}=[f\, dg\,+\,\O\, g\cdot\delta_{T_{f}^{\a}}\,+\, i\O\, g\cdot\delta_{T_{f}^{\B}}]\,\in\,\hom\left(\uh_{1}(\eta_{E}),\,\C/\Q\left\langle \O^{2},i\O^{2}\right\rangle \right).\]
Now in fact, since $f=z\circ id$, $T_{f}^{\a}=\hat{\a}$ and $T_{f}^{\B}=\hat{\B}$.
So integrating the current term by term over $\a\in H_{1}(E)\cong\uh_{1}(\eta_{E})$,
we have\[
\int_{\a}f\, dg\,=\,\int_{\a}z\, dz\,=\,\int_{-\O/2}^{\O/2}z\, dz\,=0,\]
\[
\sum_{p\in T_{f}^{\a}\cap\a}\O\cdot g(p)\,=\,\sum_{p\in\hat{\a}\cap\a}\O\cdot g(p)\,=\,\O\cdot g(\pm\O/2)\,=\,\pm\frac{\O^{2}}{2}-\xi\O,\]
\[
\sum_{p\in T_{f}^{\B}\cap\a}i\O\cdot g(p)\,=\,0\,\,\,\text{since\, }\hat{\B}\cap\a=\emptyset\,\,\,\text{(they\, are\, parallel).}\]
The only nonzero entry is nontorsion in $\C/\Z\left\langle \O^{2},i\O^{2}\right\rangle $,
which is to say nonzero in $V/\Lambda$, and so $0\neq\overline{[AJ\zeta]_{1}}$,
which $\implies$ $0\neq\left\langle \z\right\rangle \in\l^{2}CH^{2}\left((E\times E)_{k}\right).$
So we have just used basic calculus to show a cycle $\z=B(p,p-\xi)\in\ker(Alb)$
is not rationally equivalent to zero. For purposes of comparison (and
this is easy to show explicitly), the cycle $2\z_{(p,p)}$ $is$ rationally
equivalent to zero.

One can perform a similar computation for $\z=B(\xi,p,q)$ or $B(p,q-p,q-\xi)$
on $X=E\times E\times E$, where $q$ and $p$ are algebraically independent
general points on $E$, using $(14.1)$. The reader who wishes to
attempt this as an exercise should note that for $S=E\times E$, $\uh_{2}(\ns,\Q)\cong H_{2}^{tr}(E\times E,\Q)$
is spanned by just $\alpha\times\alpha-\beta\times\beta$ and $\alpha\times\beta+\beta\times\alpha$
(the Neron-Severi group has rank $4$). Whichever of these one integrates
$\overline{[AJ\zeta]_{2}}$ over, some care with signs is necessary.

In \cite[sec. 6.3]{Ke4} we will treat the more nontrivial example
$\z=B(p,q)\in Z^{2}((E\times E)_{k})$, where $(p,q)$ is the inclusion
of a very general point on the (genus $3$) Fermat quartic curve under
$\mathfrak{{F}}\to E\times E$. This is accomplished by turning $\overline{[AJ(\zeta)]_{1}}$
(in the form above for general $C_{1}\times C_{2}$) into a collection
of iterated integrals.

\section{\textbf{A $\Psi_{2}(\z)$ calculation}}

We now give an elementary proof of the very general result of \cite{RS}
regarding $0$-cycles on products of curves. This essentially contains
the $E\times E$ example above, as well as Nori's examples (which
exclude CM curves) --- but not the Fermat example. Since the cycles
considered are in the Albanese kernel, already $[AJ\zeta_{0}],\,[\zeta]_{0},\,[\zeta]_{1}$
are zero. Since they are of transcendence degree $1$, $[\zeta]_{2}=0$;
hence $[\zeta]=0$ and by Lefschetz (1,1) (see $\S4$) one may arrange
$[\bar{\zeta}]=0$, though this will also be clear from the form of
the cycle. So $[AJ(\zeta)]_{1}$ (or rather the slight quotient $[AJ(\zeta)]_{1}^{tr}$)
may be computed by Griffiths's prescription $\int_{\d^{-1}\bar{\zeta}}(\cdot)$
as a functional in {[}the numerator of{]} $(12.5)$.

Let $C_{1}$ and $C_{2}$ be smooth projective curves defined $/\bar{\Q}$,
both of positive genus. Take $o\in C_{1}$ to be a rational base point
and $p\in C_{1}(k)$ (where $k\cong\bar{\Q}(C_{1})$) to be general.
Choose a $0$-cycle $\W/\bar{\Q}$ with $\left\langle \W\right\rangle \in CH_{hom}^{1}(C_{2})$
and nontorsion image $AJ(\W)\in J^{1}(C_{2})$; write $\d^{-1}\W$
for any choice of bounding chain (connecting the dots). 

\begin{prop}
Set $X=C_{1}\times C_{2}$, and consider the $0$-cycle $\z:=((p)-(o))\times\W\in ker(\text{Alb})\subset Z^{2}(X).$
Then $\Psi_{2}(\z)$ is nontrivial.
\end{prop}
The remainder of the present section is devoted to the proof.

\subsection*{Notation:}

If $\H_{(\Q)}$ is a rational Hodge structure of weight $2m-1$, we
write $J(\H):=\H_{\C}/(F^{m}\H_{\C}+\H_{\Q})$ for the (rational)
Jacobian. (For a curve $C$, $J(C)$ is just $J(H^{1}(C))$.) Injective
and surjective maps of HS induce (resp.) injective and surjective
maps of Jacobians, and $J(\H_{1}\oplus\H_{2})=J(\H_{1})\oplus J(\H_{2})$.
Moreover for $\H_{0}\subseteq\H$, $\H/\H_{0}$ is naturally a HS
and $J(\H/\H_{0})\cong J(\H)/J(\H_{0})$.

We will work with the following (weight $3$) rational HS's, all subHS's
of $\H_{0}$:\[
\H_{0}=H^{1}(C_{1})\t H^{2}(X),\,\,\,\,\F_{0}=H^{1}(C_{1})\t F_{h}^{1}H^{2}(X),\]
\[
\H_{1}=H^{1}(C_{1})\t H^{1}(C_{1})\t H^{1}(C_{2}),\,\,\,\,\F_{1}=\F_{0}\cap\H_{1}=H^{1}(C_{1})\t F_{h}^{1}\left\{ H^{1}(C_{1})\t H^{1}(C_{2})\right\} ,\]
\[
\H_{\Delta}=\Q[\Delta_{C_{1}}]_{1}\t H^{1}(C_{2}).\]
Note in particular that $F_{h}^{1}H^{2}(X)=H_{alg}^{2}(X)$ and $(\H_{0}/\F_{0})=H^{1}(C_{1})\t H_{tr}^{2}(X)$;
and that $[\Delta_{C_{1}}]_{1}\in H^{1}(C_{1})\t H^{1}(C_{1})$ is
the $1^{\text{{st}}}$ K\"unneth component of the fundamental class,
so that $J(\H_{\Delta})\cong\Q[\Delta_{C_{1}}]_{1}\t J(C_{2})$.

\subsection*{The Setup:}

We spread $\z$ out over $S=C_{1}$ (recall $X=C_{1}\times C_{2}$),
with full spread \[
\bar{\zeta}:=(\Delta_{C_{1}}-C_{1}\times\{ o\})\times\W\,\in\, Z^{2}(C_{1}\times(C_{1}\times C_{2})_{/\bar{\Q}}).\]
This has $[\bar{\zeta}]=0$ (because $\W\homeq0$), with bounding
chain $\d^{-1}\bar{\zeta}:=(\Delta_{C_{1}}-C_{1}\times\{ o\})\times\d^{-1}\W.$
The image of the functional\[
\left\{ \int_{\d^{-1}\bar{\zeta}}(\cdot)\right\} \in\frac{\left\{ F^{2}\left(H^{1}(C_{1},\C)\t H^{2}(X,\C)\right)\right\} ^{\vee}}{\text{{im}}\left\{ H_{1}(C_{1},\Q)\t H_{2}(X,\Q)\right\} }\cong J(\H_{0})\]
under the projection to $(12.5)\cong J(\H_{0}/\F_{0})$, is obviously
$[AJ(\zeta)]_{1}^{tr}$ and we must show this is not zero. Consider
the commutative diagram

\begin{equation*} \xymatrix{ & J(\H_0) \ar @{>>} [r]^{p_0} \ar @{>>} [d]^{p_1} & J(\H_0 / \F_0) \ar @{>>} [d]^{\overline{p}_1} \\ J(\H_{\Delta}) \ar @{^{(}->} [r]^{\iota} & J(\H_1) \ar @{>>} [r]^{p_2} & J(\H_1 / \F_1) \ar @{=} [r] & {\frac{J(\H_1)}{J(\F_1)}} \, . } \end{equation*}${}$\\
One easily shows%
\footnote{Viewing $\H_{0},\,\H_{1}$ as their own duals, $p_{1}$ is induced
by restriction of functionals $(F^{2}\H_{0}^{\C})^{\vee}\twoheadrightarrow(F^{2}\H_{1}^{\C})^{\vee}$
and kills $\int_{C_{1}\times\{ o\}\times\d^{-1}\W}(\cdot)$; hence
$p_{1}\left(\int_{\d^{-1}\bar{\zeta}}(\cdot)\right)=p_{1}\left(\int_{\Delta_{C_{1}}\times\d^{-1}\W}(\cdot)\right)$.
Under $(H^{1}(C_{1})^{\t2})^{\vee}\cong H^{1}(C_{1})^{\t2}$, $\int_{\Delta_{C_{1}}}(\cdot)$
identifies with $[\Delta_{C_{1}}]_{1}$; while $\int_{\d^{-1}\W}(\cdot)\in\left(F^{1}H^{1}(C_{2},\C)\right)^{\vee}$
yields $AJ(\W)$.%
} that $\iota\left([\Delta_{C_{1}}]_{1}\t AJ_{C_{2}}(\W)\right)=p_{1}\left(\int_{\d^{-1}\bar{\zeta}}(\cdot)\right)$
in $J(\H_{1})$. By hypothesis $AJ(\W)\neq0$, so this is nonzero;
we are done if we can show $p_{2}$ does not kill it, since then $\overline{p}_{1}\left([AJ(\zeta)]_{1}^{tr}\right)\neq0$.
This is accomplished if $(p_{2}\circ\iota)$ is injective, which would
follow if $J(\F_{1})\oplus J(\H_{\Delta})\hookrightarrow J(\H_{1})$.

\subsection*{The Hodge-theoretic input:}

Thus we have merely to show that\[
\F_{1}\cap\H_{\Delta}=\{0\}\]
in $\H_{1}$ (we actually show $\F_{1}^{\C}\cap\H_{\Delta}^{\C}=\{0\}$).
Let $\{\omega_{i}\}\subseteq\Omega^{1}(C_{1}{}_{/\C})$ be a basis
satisfying $\sqrt{-1}\int_{C_{1}}\omega_{i}\wedge\overline{\omega}_{j}=\delta_{ij}$.
An element $\mathsf{{f}}\in\F_{1}^{\C}$ is of the form $\sum\omega_{i}\t A_{i}+\sum\overline{\omega}_{i}\t B_{i}$
where $A_{i},\, B_{i}$ belong to $\left\{ F_{h}^{1}\left(H^{1}(C_{1})\t H^{1}(C_{1})\right)\right\} _{\C}$
(and are hence of type $(1,1)$). Suppose $\mathsf{{f}}$ also belongs
to $\H_{\Delta}^{\C}$; that is, $\mathsf{{f}}=[\Delta_{C_{1}}]_{1}\t\Gamma$
for some $\Gamma\in H^{1}(C_{2},\C)$. Writing \[
\sum\omega_{i}\t A_{i}+\sum\overline{\omega}_{i}\t B_{i}\,=\, i\left(\sum\omega_{i}\t\overline{\omega}_{i}-\sum\overline{\omega}_{i}\t\omega_{i}\right)\t\Gamma,\]
we have $A_{1}=\overline{\omega}_{1}\t\Gamma$, $B_{1}=\omega_{1}\t\Gamma$
which is impossible in terms of $type$, unless $\Gamma=0$. Hence
$\mathsf{{f}}=0$ and we are done.

\section{\textbf{A surprise from transcendental number theory}}

In this section we consider $\chi_{2}(\z)$ for the cycles of $\S15$,
and ask whether $\chi_{2}$ misses anything that $\Psi_{2}$ picks
up.

Here is another way of putting Proposition $15.1$: let $C_{1},\, C_{2}$
(both of positive genus) be defined $/\bar{\Q}$ with $p\in C_{1}(\C)$
very general, $o\in C_{1}(\bar{\Q})$ and $\W\in Z_{0}^{hom}(C_{2}{}_{/\bar{\Q}})$.
Set $\z=(p-o)\times\W$. Then the following four conditions are equivalent:

\subsubsection*{(I)}

$AJ(\W)$ is nonzero in $J^{1}(C_{2})$ (i.e., nontorsion in the integral
Jacobian).

\subsubsection*{(II)}

$[AJ(\zeta)]_{1}^{tr}\neq0$ in $(12.5)$.

\subsubsection*{(III)}

$\Psi_{2}(\z)\neq0$.

\subsubsection*{(IV)}

$\left\langle \z\right\rangle $ is nonzero in $CH_{0}(C_{1}\times C_{2}\textrm{{}}_{/\C})$.

\begin{proof}
We have shown \emph{(I) $\implies$ (II) $\implies$ (III) $\implies$
(IV)}, and {}``not \emph{(I)}'' $\implies$ {}``not \emph{(IV)}''
is clear.
\end{proof}
Recall that $\chi_{2}(\z)$ is a projected image of $\Psi_{2}(\z)$
(and easier to compute). In the situation at hand (i.e. $\dim X=2$,
$\dim S=1$), it is just \[
\overline{[AJ(\zeta)]_{1}}\in\frac{\left\{ H^{1}(S,\C)\t H^{2,0}(X,\C)\right\} ^{\vee}}{\text{{im}}\left\{ H_{1}(S,\Q)\t H_{2}(X,\Q)\right\} },\]
which can be obtained from $[AJ(\zeta)]_{1}^{tr}\in(12.5)$ by forgetting
the $H_{tr}^{1,1}(X)$-part of the functional. But it does surprisingsly
well at detecting cycles as far as we can compute:

\begin{prop}
With the notation above, and $C_{i}=E_{i}$ \emph{(}$i=1,2$\emph{)}
\textbf{elliptic} curves $/\bar{\Q}$,$\vspace{2mm}$\\
(V). $\chi_{2}(\z)\neq0$$\vspace{2mm}$\\
is equivalent to the conditions (I) thru (IV), assuming a conjecture
of Waldschmidt for the case where $E_{1}$ and $E_{2}$ are non-isogenous
non-CM curves.
\end{prop}
In particular, we have not been able to find a $0$-cycle $\z$ with
$\chi_{2}(\z)=0$ but $\Psi_{2}(\z)\neq0$.

The rest of the section is the proof. We will show that for $\z$
of the above form, assuming $\overline{[AJ(\zeta)]_{1}}=0$ and $AJ_{E_{2}}(\W)\neq0$
leads to a contradiction. (Here as usual $AJ(\W)$ is the image in
the $rational$ Jacobian.) Hence \emph{(I)} $\implies$ \emph{(V)}
(the converse, of course, is clear).

Note first that the rational Neron-Severi group $NS(E_{1}\times E_{2})\t\Q$
has rank $2,\,3,\,4$ according as $E_{1}$ and $E_{2}$ are nonisogenous,
isogenous without CM, isogenous with CM. In the latter (rank $4$)
case, $H_{tr}^{1,1}(E_{1}\times E_{2})=0$ and so $[AJ(\zeta)]_{1}^{tr}$,
$\overline{[AJ(\zeta)]_{1}}$ are the same invariant. So in this case
we have immediately $\overline{[AJ(\zeta)]_{1}}=0$ $\implies$ $AJ(\W)=0$.

\subsection*{Notation:}

We may write the elliptic curves $E_{j}$ ($j=1,2$) in Weierstrass
form $\{\overline{y^{2}=4x^{3}-(g_{2})_{j}x-(g_{3})_{j}}\}$ with
$(g_{2})_{j}$, $(g_{3})_{j}\,\in\,\bar{\Q}$, $\omega_{j}=$ pullback
to $E_{j}$ of $dx/y$, and $o=$ point at $\infty$ (and $p\in E_{1}(\C)$
very general). Pick $\alpha_{j},\,\beta_{j}\,\in\, Z_{1}^{top}(E_{j})$
spanning $H_{1}(E_{j},\Q)$ such that $\alpha_{j}\cdot\beta_{j}=+1$,
and write periods $\Omega_{j\alpha}=\int_{\alpha_{j}}\omega_{j}$,
$\Omega_{j\beta}=\int_{\beta_{j}}\omega_{j}$, $\tau_{j}:=\Omega_{j\beta}/\Omega_{j\alpha}$
(with $\Im(\tau_{j})>0$). Up to $\ratequiv$, $\W$ is just $(q)-(o)$
for some $q\in E_{2}(\bar{\Q})$, and we write $\d^{-1}\W=\overrightarrow{o.q}$
for some fixed choice of path. Integrating along it yields $\xi:=\int_{o}^{q}\omega_{2}\in\C$,
which projects to $AJ_{E_{2}}(\W)\,\in\,\C/\Q\left\langle \Omega_{2\alpha},\Omega_{2\beta}\right\rangle $.
Finally, let $A_{j}$ and $B_{j}$ be Poincar\'e duals ($C^{\infty}$
$1$-forms) for $\alpha_{j},\,\beta_{j}$; and $\wp_{j}$ denote the
Weierstrass $\wp$-function on $\left(\C/\Z\left\langle \Omega_{j\alpha},\Omega_{j\beta}\right\rangle \right)\cong E_{j}(\C)$.

\subsection*{Computations:}

We compute $\chi_{2}(\z)$, for $\z=((p)-(o))\times((q)-(o))$. If
we write $\d^{-1}\bar{\zeta}=\left(\Delta_{E_{1}}-E_{1}\times\{ o\}\right)\times\overrightarrow{o.q}$,
evaluating $\int_{\d^{-1}\bar{\zeta}}(\cdot)\in\left\{ H^{1}(E_{1})\t H^{2,0}(E_{1}\times E_{2})\right\} ^{\vee}$
on the basis $\mathfrak{{B}}=\{ A_{1}\t\omega_{1}\wedge\omega_{2},\, B_{1}\t\omega_{1}\wedge\omega_{2}\}$
yields a vector\[
\vec{\nu}_{\z}=\left(\Omega_{1\alpha}\xi,\,\Omega_{1\beta}\xi\right)\,\in\,\C^{2}.\]
(The $E_{1}\times\{ o\}$ term in $\d^{-1}\bar{\zeta}$ plays no role.)
The projection of $\vec{\nu}_{\z}$ to $\C^{2}/\mathbb{{L}}$ gives
$\overline{[AJ(\zeta)]_{1}}$, where $\mathbb{{L}}$ denotes the $\Q$-lattice
generated by the vectors produced by evaluating a basis of $H_{1}(E_{1},\Q)\t H_{2}(E_{1}\times E_{2},\Q)$
on $\mathfrak{{B}}$. Obviously a basis for $H_{1}(E_{1},\Q)^{\otimes2}\otimes H_{1}(E_{2},\Q)$
will suffice; taking this to be\[
\mathfrak{{C}}\,=\,\{\alpha_{1}\t\alpha_{1}\t\alpha_{2},\,\alpha_{1}\t\alpha_{1}\t\beta_{2},\,\alpha_{1}\t\beta_{1}\t\alpha_{2},\,\alpha_{1}\t\beta_{1}\t\beta_{2},\]
\[
\mspace{50mu}\beta_{1}\t\alpha_{1}\t\alpha_{2},\,\beta_{1}\t\alpha_{1}\t\beta_{2},\,\beta_{1}\t\beta_{1}\t\alpha_{2},\,\beta_{1}\t\beta_{1}\t\beta_{2}\},\]
we have the following generators for $\mathbb{{L}}$: \[
(0,\Omega_{1\alpha}\Omega_{2\alpha}),\,(0,\Omega_{1\alpha}\Omega_{2\beta}),\,(0,\Omega_{1\beta}\Omega_{2\alpha}),\,(0,\Omega_{1\beta}\Omega_{2\beta}),\]
\[
(\Omega_{1\alpha}\Omega_{2\alpha},0),\,(\Omega_{1\alpha}\Omega_{2\beta},0),\,(\Omega_{1\beta}\Omega_{2\alpha},0),\,(\Omega_{1\beta}\Omega_{2\beta},0).\]

\subsection*{Main argument:}

Assume $\overline{[AJ(\zeta)]_{1}}=0$ and $AJ_{E_{2}}((q)-(o))\neq0$
(in the rational Jacobian). The second condition means that $\xi\notin\Q\left\langle \Omega_{2\alpha},\Omega_{2\beta}\right\rangle $.
The first translates to $\vec{\nu}_{\z}\in\mathbb{{L}}$, hence to
the two equations

\begin{equation} \xi = a_1\Omega_{2\alpha} + a_2 \Omega_{2\beta} +a_3 \tau_1 \Omega_{2\alpha} + a_4 \tau_1 \Omega_{2\beta} \end{equation}
\begin{equation} \xi = b_1 \tau_1^{-1} \Omega_{2\alpha} + b_2 \tau_1^{-1} \Omega_{2\beta} + b_3 \Omega_{2\alpha} + b_4 \Omega_{2\beta} , \end{equation}${}$\\
for some $a_{k},b_{k}\in\Q$. We cannot have $a_{3}$ and $a_{4}$
both zero; nor can $b_{1}$ and $b_{2}$ both be zero. Without loss
of generality, we may assume $b_{3}=b_{4}=0$. (This is seen by adding
topological $1$-cycles to $\overrightarrow{o_{2}.q}$, which effects
translations of $\xi$ and $\vec{\nu}_{\z}$ by respective elements
of $\Q\left\langle \Omega_{2\alpha},\Omega_{2\beta}\right\rangle $
and $\mathbb{{L}}$.) Setting the right-hand sides of $(16.1)$ and
$(16.2)$ equal gives

\begin{equation} \frac{b_1 \Omega_{2\alpha}+b_2 \Omega_{2\beta}}{\tau_1} = a_1 \Omega_{2\alpha} + a_2 \Omega_{2\beta} +(a_3 \Omega_{2\alpha} + a_4 \Omega_{2\beta})\tau_1, \end{equation}and
multiplying by $\tau_{1}/\Omega_{2\alpha}$ and rearranging yields

\begin{equation} \tau_1^2 (a_3 + a_4 \tau_2) + \tau_1 (a_1+a_2\tau_2) - (b_1+b_2\tau_2) =0. \end{equation}${}$\\
Recall that (since $E_{i}$ is defined $/\bar{\Q}$) each $\tau_{i}$
is either a quadratic irrationality ($\Leftrightarrow$ $E_{i}$ has
CM) or transcendental.$\vspace{2mm}$\\
\emph{$\underline{\text{{Case\,1}}}$: Assume $E_{1},\, E_{2}$ non-isogenous
non-CM (the {}``general'' case).}$\vspace{2mm}$\\
Set $u_{1}=\Omega_{2\alpha}$, $u_{2}=\Omega_{2\beta}$, $u_{3}=\frac{b_{1}\Omega_{2\alpha}+b_{2}\Omega_{2\beta}}{\tau_{1}}(\neq0)$.
Then $\tau_{1}=\frac{b_{1}u_{1}+b_{2}u_{2}}{u_{3}}$ and $(16.3)$
says\[
u_{3}=a_{1}u_{1}+a_{2}u_{2}+(a_{3}u_{1}+a_{4}u_{2})\left(\frac{b_{1}u_{1}+b_{2}u_{2}}{u_{3}}\right);\]
multiplying by $u_{3}$ gives\[
u_{3}^{2}=a_{1}u_{1}u_{3}+a_{2}u_{2}u_{3}+a_{3}b_{1}u_{1}^{2}+a_{4}b_{2}u_{2}^{2}+(a_{3}b_{2}+a_{4}b_{1})u_{1}u_{2}.\]
This is a nontrivial algebraic relation amongst the $\{ u_{i}\}$
(the coefficient of $u_{3}^{2}$ is $1$). Now assume the following,
a more general version of which appears as Conjecture $4.1$ in \cite{W}:

\subsection*{Conjecture:}

\emph{Let $\wp$ denote a Weierstrass $\wp$-function with algebraic
invariants} ($g_{2},g_{3}\in\bar{\Q}$). \emph{Let} $u_{i}\in\C$
($i=1,\ldots,n$) \emph{be such that} (\emph{for each $i$}) \emph{}$u_{i}$
\emph{is a pole of $\wp$ or} $\wp(u_{i})\in\bar{\Q}$. \emph{Then
$\{ u_{i}\}$ algebraically dependent} $\implies$ $\{ u_{i}\}$ \emph{linearly
dependent over}\[
K_{\tau}:=\left\{ \begin{array}{c}
\Q\,\,\,\,\,\,\,\,\tau\notin\bar{\Q}\\
\Q(\tau)\,\,\,\,\,\,\,\,\tau\in\Q\end{array}\right..\]

\begin{rem*}
If true, this (apparently standard) conjecture absorbs the result
of Schneider quoted under the {}``transcendental input'' heading
below. 
\end{rem*}
Applying this with $\wp=\wp_{2}$ (and $K=\Q$ since $E_{2}$ doesn't
have CM) yields ($A,B,C\in\Q$ such that)\[
0=Au_{1}+Bu_{2}+Cu_{3}\]
\[
=A\Omega_{2\alpha}+B\Omega_{2\beta}+Cb_{1}(\Omega_{2\alpha}/\tau_{1})+Cb_{2}(\Omega_{2\beta}/\tau_{1});\]
multiplying by $\tau_{1}/\Omega_{2\alpha}$ gives\[
0=A\tau_{1}+B\tau_{1}\tau_{2}+Cb_{1}+Cb_{2}\tau_{2}\]
\[
\implies\,\,\tau_{2}=-\frac{A\tau_{1}+Cb_{1}}{B\tau_{1}+Cb_{2}}.\]
The latter formula implies the existence of an isogeny $E_{1}\simeq E_{2}$,
contradicting our assumption.$\vspace{2mm}$\\
\emph{$\underline{\text{{Case\,2}}}$: Assume $E_{1}\simeq E_{2}$
non-CM.}$\vspace{2mm}$\\
Then for $A,B,C,D\in\Q$,\[
\tau_{2}=\frac{A+B\tau_{1}}{C+D\tau_{1}}\]
\[
\implies\,\,\,0=A+B\tau_{1}-C\tau_{2}-D\tau_{1}\tau_{2}.\]
 Multiplying by $\Omega_{2\alpha}$ gives

\begin{equation} 0= A\Omega_{2\alpha} + B \tau_1 \Omega_{2\alpha} - C \Omega_{2\beta} - D\tau_1 \Omega_{2\beta} . \end{equation}${}$\\
\emph{Case 2a:} $D\neq0$.$\vspace{1mm}$\\
Multiply $(16.5)$ by $a_{4}/D$ and add to $(16.1)$ to get\[
(16.1')\mspace{160mu}\xi=a_{1}'\Omega_{2\alpha}+a_{2}'\Omega_{2\beta}+a_{3}'\tau_{1}\Omega_{2\alpha}.\mspace{170mu}\]
\[
(16.4')\mspace{130mu}a_{3}'\tau_{1}^{2}+(a_{1}'+a_{2}'\tau_{2})\tau_{1}-(b_{1}+b_{2}\tau_{2})=0.\mspace{140mu}\]
Substituting for $\tau_{2}$ and multplying by $(C+D\tau_{1})$ gives\[
0=a_{3}'\tau_{1}^{2}(C+D\tau_{1})+\tau_{1}\{ a_{1}'(C+D\tau_{1})+a_{2}'(A+B\tau_{1})\}-\{ b_{1}(C+D\tau_{1})+b_{2}(A+B\tau_{1})\}\]
\[
=\tau_{1}^{3}(a_{3}'D)+\tau_{1}^{2}(a_{3}'C+a_{1}'D+a_{2}'B)+\tau_{1}(a_{1}'C+a_{2}'A-b_{1}D-b_{2}B)-(b_{1}C+b_{2}A).\]
This polynomial in $\tau_{1}\notin\bar{\Q}$ must vanish, hence $0=a_{3}'D$.
But $D=0$ contradicts the hypothesis of \emph{Case 2a}, while $a_{3}'=0$
together with $(16.1')$ contradicts $\xi\notin\Q\left\langle \Omega_{2\alpha},\Omega_{2\beta}\right\rangle .$$\vspace{1mm}$\\
\emph{Case 2b:} $D=0$. $\vspace{1mm}$\\
Then we rewrite $\tau_{2}=A+B\tau_{1}$ (absorbing $C$). Substituting
this into $(16.4)$ gives \[
0=\tau_{1}^{2}\{(a_{3}+a_{4}A)+a_{4}B\tau_{1}\}+\tau_{1}\{(a_{1}+a_{2}A)+a_{2}B\tau_{1}\}-\{(b_{1}+b_{2}A)+b_{2}B\tau_{1}\}\]
\[
=(a_{4}B)\tau_{1}^{3}+(a_{3}+a_{4}A+a_{2}B)\tau_{1}^{2}+(a_{1}+a_{2}A-b_{2}B)\tau_{1}-(b_{1}+b_{2}A).\]
Again the polynomial must vanish since $\tau_{1}\notin\bar{\Q}$,
so in particular $b_{1}=-b_{2}A$. By $(16.2)$ (with $b_{3}=b_{4}=0$,
of course), and $\tau_{1}=\frac{\tau_{2}-A}{B}$, \[
\xi=\tau_{1}^{-1}(b_{1}\Omega_{2\alpha}+b_{2}\Omega_{2\beta})=\left(\frac{\Omega_{2\alpha}B}{\tau_{2}-A}\right)(b_{1}+b_{2}\tau_{2})\]
\[
=\left(\frac{\Omega_{2\alpha}B}{\tau_{2}-A}\right)(-b_{2}A+b_{2}\tau_{2})=\Omega_{2\alpha}Bb_{2}\in\Q\left\langle \Omega_{2\alpha},\Omega_{2\beta}\right\rangle ,\]
again a contradiction.$\vspace{2mm}$\\
\emph{$\underline{\text{{Case\,3}}}$: $E_{1}$ or $E_{2}$ has CM.}$\vspace{2mm}$\\
If $E_{2}$ has CM then $\tau_{2}\in\bar{\Q}$ and $(16.4)\,\implies\,\tau_{1}\in\bar{\Q}$,
\emph{hence $E_{1}$ has CM.} If $E_{2}\simeq E_{1}$ then $\text{{rank}}(NS)=4$
leads to a contradiction as indicated above; so $E_{1}\,\,\nsimeq E_{2}$
$\implies$ $\Q\left\langle 1,\tau_{1}\right\rangle \neq\Q\left\langle 1,\tau_{2}\right\rangle $.

Hence this case reduces to: $E_{1}$ of CM, $\tau_{1}\notin\Q\left\langle 1,\tau_{2}\right\rangle $.
Returning to $(16.1)$, choose $M\in\Z^{+}$ sufficiently large that
$Ma_{i}\in\Z$ for $i=1,2,3,4$; and set $\lambda_{0}:=Ma_{3}\Omega_{2\alpha}+Ma_{4}\Omega_{2\beta}\in\Z\left\langle \Omega_{2\alpha},\Omega_{2\beta}\right\rangle $.
Then we have $M\xi\equiv\tau_{1}\lambda_{0}$ in $\C/\Z\left\langle \Omega_{2\alpha},\Omega_{2\beta}\right\rangle .$
Pick $N\in\Z^{+}$ such that $\lambda:=\frac{\lambda_{0}}{N}\notin\Z\left\langle \Omega_{2\alpha},\Omega_{2\beta}\right\rangle $.
Since $E_{2}/\bar{\Q}$, torsion points are algebraic and $\wp_{2}(\lambda)\in\bar{\Q}$.

\subsection*{The transcendental input:}

It is common knowledge that the periods $\Omega_{\alpha},\,\Omega_{\beta}$
of an elliptic curve defined $/\bar{\Q}$ are transcendental. (For
this and the result that follows, it is essential that $g_{2},\, g_{3}\in\bar{\Q}$
and $\omega[=dz]$ be defined $/\bar{\Q}$ {[}e.g., $dx/y${]} ---
rescaling is not allowed.) 

Less well-known (but quite useful) is a result of T. Schneider \cite[Satz 16 (p. 61)]{Sc}
which can be read as follows for $E/\bar{\Q}$ with associated Weierstrass
$\wp(z)$, writing $K_{\tau}=\Q(\tau)$ if $E$ has CM and $K_{\tau}=\Q$
otherwise:$\vspace{2mm}$\\
\emph{If $u\in\C\m\Z\left\langle \Omega_{\alpha},\Omega_{\beta}\right\rangle $
is such that $\wp(u)\in\bar{\Q}$, and $\mathsf{{b}}\in\bar{\Q}\m K_{\tau}$
is such that $\mathsf{{b}}\cdot u\notin\Z\left\langle \Omega_{\alpha},\Omega_{\beta}\right\rangle $,
then $\wp(\mathsf{{b}}\cdot u)$ is transcendental.}$\vspace{2mm}$\\
We apply this to $\wp=\wp_{2}$ with $u=\lambda$ and \textbf{$\mathsf{{b}}=N\tau_{1}$};
note that since $K_{\tau_{2}}\subseteq\Q\left\langle 1,\tau_{2}\right\rangle $
(whether or not $E_{2}$ has CM) and $\tau_{1}\notin\Q\left\langle 1,\tau_{2}\right\rangle $,
$N\tau_{1}\notin K_{\tau_{2}}$. Moreover, since $E_{1}$ has CM,
$N\tau_{1}\in\bar{\Q}$. So the requirements are met and $\wp_{2}(N\tau_{1}\cdot\lambda)=\wp_{2}(M\xi)\notin\bar{\Q}$;
since the group law is algebraic $\wp_{2}(\xi)\notin\bar{\Q}$. But
this is the $x$-coordinate of $q$, contradicting $q\in E_{2}(\bar{\Q})$.
So the proof is finished.

\begin{rem}
We emphasize that the obstruction to the existence of a cycle (of
the above form) in $\ker(\chi_{2})$ but not $\ker(\Psi_{2})$, is
not Hodge-theoretic in nature. One can repeat all the constructions
of this paper (except the expected injectivity of $\Psi$ and its
consequences) with the spread field $\bar{\Q}$ replaced by some $\bar{L}(\supseteq\bar{\Q})$.
This yields (weaker) analogues of $\chi_{2}$ and $\Psi_{2}$ which
are defined on $\ker(Alb)$. If we take $\z=((p)-(o))\times((q)-(o))$
with $q\in E_{2}(\bar{L})$ and $p$ not defined $/\bar{L}$, the
computations proceed as above, with the same Hodge-theoretic target
spaces. Moreover, in this case there is no contradiction at the end
of the above proof.

For example, take $E_{1},\, E_{2}/\bar{\Q}$ with $\tau_{1}=\sqrt{-1}$
(so $E_{1}$ has CM) and pick $q=\left(\wp_{2}(\sqrt{-1}\Omega_{2\alpha}),\wp_{2}'(\sqrt{-1}\Omega_{2\alpha})\right)\in E_{2}(\C)$
so that $AJ_{E_{2}}\left((q)-(o)\right)\equiv\sqrt{-1}\Omega_{2\alpha}$.
This $q$ is actually defined over some $L$ with $\text{{trdeg}}(L/\bar{\Q})=1$,
and we can tensor $E_{1}$ and $E_{2}$ up to $\bar{L}$. Then we
take $p\in E_{1}(\C)$ very general and not defined $/\bar{L}$, and
spread $\z$ (defined as above) out over $\bar{L}$. Essentially by
$\S15$ $\Psi_{2}^{\bar{L}}(\z)\neq0$, while $\chi_{2}^{\bar{L}}(\z)$
identifies with $\vec{\nu}_{\z}=$$(\Omega_{1\alpha}\cdot\sqrt{-1}\Omega_{2\alpha},\Omega_{1\beta}\cdot\sqrt{-1}\Omega_{2\alpha})=$
$(\tau_{1}\Omega_{1\alpha}\cdot\Omega_{2\alpha},-\frac{1}{\tau_{1}}\Omega_{1\beta}\cdot\Omega_{2\alpha})=$
$(\Omega_{1\beta}\Omega_{2\alpha},-\Omega_{1\alpha}\Omega_{2\alpha})\in\mathbb{{L}}$
and hence is $0$ in $\C^{2}/\mathbb{{L}}$.
\end{rem}

\section{\textbf{$0$-cycles on special Kummer surfaces}}

Let $\EE\mappi\mathbb{{P}}^{1}$ be an elliptic fibration (defined
$/\bar{\Q}$) where $\EE$ is an exceptional $K3$ surface (i.e. $\text{rank}(Pic(\EE))=20$);
the Fermat quartic is the example we have in mind. The following questions
arise, if we take a smooth {}``algebraic'' fiber $E_{q}=\pi^{-1}(q)$
($q\in\bar{\Q}$), with general point $p$ and algebraic base point
$o/\bar{\Q}$. Since $H^{0}(\O_{\EE}^{1})=0$, \[
\z_{\EE}:=(p)-(o)\in ker(\text{Alb}).\]
Is $\Psi_{2}(\z_{\EE})$ nontrivial? Is $\z_{\EE}\ratnequiv0$?

One approach to the first question is to try and compute it by the
differential character approach, i.e. the membrane integrals $(11.1)$.
Since $trdeg(\z/\bar{\Q})=1$, $[\zeta]_{2}=0$ and so $\Psi_{2}(\z)=[AJ\zeta]_{1}$;
moreover since $H_{tr}^{1,1}(\EE)=0$, $\chi_{2}(\z)=[AJ(\zeta)]_{1}^{tr}$.
So this is a rare instance (like the self-product of a CM elliptic
curve in $\S14$) where the membrane integrals derived from loops
on $S\,[=E_{q}\text{\, here}]$ compute (essentially) the entire invariant.

Let's be concrete about these: spreading $\{ p\}$ merely gives the
image of $\Delta_{E_{q}}$ under $E_{q}\times E_{q}\mapidi E_{q}\times\EE$.
Take cycles $\alpha,\,\beta\in H_{1}(E_{q})$ and let $\mu_{\a},\,\mu_{\B}$
be bounding membranes $\in C_{2}(\EE)$ (possible since $H_{1}(\EE)=0$).
One merely has to show that $\int_{\mu_{\a}}$ and $\int_{\mu_{\B}}$
do not both give torsion classes in $\{ H^{2,0}(\EE,\C)\}^{\vee}/im\{ H_{2}(\EE,\Z)\}$.
Moreover $H^{2,0}(\EE)$ is generated by some $\O$ which has only
$2$ nontrivial periods. Now for instance the Fermat quartic has {}``CM''
(in this case, an order-$4$ automorphism) and it follows that the
period ratio is $1:i$. So we merely have to show that \textbf{no
integral multiple of $\int_{\mu_{\a}}\O\left/\int_{\mu_{\B}}\O\right.$
is a Gaussian integer}. (At least this would be $sufficient$.)

We venture a guess that it is transcendental, but pursue an alternate
route (and answer only the second question above). Consider the fact
that the Fermat quartic is a \emph{special Kummer surface} (\cite{PS},
\cite{PL}). That is, it may be constructed as a blow-up of a quotient
of a product of elliptic curves. Namely, if $E_{1}(\C)\cong\C/\Z\left\langle 1,\frac{i}{2}\right\rangle $
and $E_{2}(\C)\cong\C/\Z\left\langle 1,i\right\rangle $ then the
action of $(-id_{1},-id_{2}):\, E_{1}\times E_{2}\to E_{1}\times E_{2}$
induces a quotient morphism $\theta:\, E_{1}\times E_{2}\to\hat{\EE}$.
Products of $2$-torsion points give rise to $16$ singular points;
blowing these up yields $b:\,\EE\to\hat{\EE}$, with $\EE$ a $K3$
surface fibering over $E_{2}/\left\langle -id_{2},id_{2}\right\rangle \,\cong\,\P^{1}.$

For example, if $E_{2}=\overline{\{ y^{2}=x^{3}-5x\}}$ then $\{ x,y\}\mapsto\{ x\}\in\P^{1}$
gives the map $E_{2}\maprho\mathbb{{P}}^{1}$. The singular fibers
of $\EE\mappi\mathbb{{P}}^{1}$ lie over the $\rho-images$ of the
$2$-torsion points, namely $\{\pm\sqrt{5},0,\infty\}.$ They consist
of (connected) configurations of rational curves. So any two points
on a singular fiber are rationally equivalent. On the other hand,
the restriction of $\theta$ to $E_{1}\times(E_{2}\m\{2\text{-torsion}\})$
gives an \'etale double cover of the smooth fibers $\pi^{-1}\left(\mathbb{{P}}^{1}\m\{\pm\sqrt{5},0,\infty\}\right)$. 

Now take points $o=\rho^{-1}(\infty)$ and $\xi\in E_{2}(\bar{\Q})$
such that the class $(\xi)-(o)$ in nontorsion, and write $q:=\rho(\xi)\in\mathbb{{P}}^{1}(\bar{\Q})\m\{\pm\sqrt{5},0,\infty\}$.
Also write $o\in E_{1}(\bar{\Q})$ for a base point and $p\in E_{1}(\C)$
for a very general point. As above, we denote the push-forward of
$(p)-(o)$ under the composition $E_{1}=E_{1}\times\{\xi\}\rTo_{b^{-1}\circ\theta}^{\cong}\pi^{-1}(q)=:E_{q}\subseteq\EE$
by $\z_{\EE}$. Let $\W$ be another $0$-cycle on $\EE$, consisting
of a difference of two points in $\pi^{-1}(\infty)$ such that $b_{*}(\W)=\theta_{*}((o,o)-(p,o))$.
Finally, define $\z:=(p,\xi)-(o,\xi)-(p,o)+(o,o)\in Z_{0}(E_{1}\times E_{2}{}_{/\C})$.

From Proposition $15.1$, we know $\z\ratnequiv0$; while $\W\ratequiv0$
is obvious. Moreover, $\z$ is the $2$-box $B(p,\,\xi)$ with $\theta^{*}\theta_{*}\z=B(p,\,\xi)+B(-p,\,-\xi)$;
explicit rational equivalences show $B(p,\,\xi)\ratequiv B(-p,\,-\xi)$,
and so $\theta^{*}\theta_{*}\z\ratequiv2\z\ratnequiv0$. It follows
that $0\ratnequiv\theta_{*}\z=b_{*}(\z_{\EE}+\W)$, so that $\z_{\EE}+\W\ratnequiv0$;
hence $\z_{\EE}\ratnequiv0$, as we wanted to show.

Of course, this is not surprising in light of the isomorphism between
the Albanese kernels of an abelian variety and its associated Kummer
surface, proved by Bloch in \cite[appendix]{BKL}.

\address{\small  \noun{Department of Mathematics, University of Chicago,
Chicago, IL 60637 USA}}

\email{matkerr@math.uchicago.edu}
\end{document}